\documentclass[12pt]{iopart} 
\usepackage{latexsym}  
\usepackage{amsfonts}  
\usepackage{iopams}  
\usepackage[dvips]{graphicx}  
\usepackage{color}  

\usepackage{times}  
 
\usepackage{amssymb}  
  
\font\mat=msbm10  
\def\reel{\mat{R}}  
\def\R{\mbox{\reel}}  
\def\Z{\mbox{$\mathbb{Z}$}}  
\def\P{\mbox{$\mathbb{P}$}}  
\def\I{\mbox{$\mathbb{I}$}}

\def\u{\mbox{$\mathbf{u}$}}  
\def\v{\mbox{$\mathbf{v}$}} 
\def\j{\mbox{$\mathbf{j}$}}  
\def\k{\mbox{$\mathbf{k}$}}  
\def\n{\mbox{$\mathbf{n}$}}  
\def\x{\mbox{$\mathbf{x}$}}  
\def\iint{\int\!\!\!\!\int} 
  
\begin{document}  
\title[Divergence-free Wavelets for Navier-Stokes]{Towards a divergence-free wavelet method for the simulation of 2D/3D turbulent flows}  
\author{Erwan Deriaz\dag\  
\footnote[3]{To 
whom correspondence should be addressed (Erwan.Deriaz@imag.fr)} 
and Val\'erie Perrier\dag} 
 
\address{\dag\ Laboratoire de Mod\'elisation et Calcul de l'IMAG,  
BP 53 - 38 041 Grenoble Cedex 9, France}

\begin{abstract} 
In this paper, we investigate the  use of  
compactly supported {\it divergence-free wavelets} for the  
representation  
of the Navier-Stokes solution. After reminding the theoretical construction of  
divergence-free wavelet vectors, we present in detail the bases and  
corresponding fast algorithms for 2D and 3D incompressible flows.  
In order to compute the nonlinear term, we propose a new method which  
provides in practice with the Hodge decomposition of any flow: this  
decomposition  
enables us to separate the incompressible part of the flow from its  
orthogonal complement, which corresponds to the gradient component of  
the flow.  
Finally we show numerical tests to validate our approach.  
  
\end{abstract}  
 
\submitto{Journal of Turbulence} 
 
\section{Introduction}  
  
The prediction of fully-developed turbulent flows represents an  
extremely challenging field of research in scientific computing. The  
{\it Direct Numerical Simulations} (DNS) of turbulence  requires  
the integration in time of the nonlinear Navier-Stokes equations, which assumes the  
computation of all scales of motion.  
However,  
at large Reynolds number, turbulent flows generate increasingly small  
scales: to be realistic, the discretization in space (and correlatively  
in time) ought to handle a huge number of degrees of freedom, that is  
in 3D  out of the reach of  
available computers.  
  
 Many tentatives have been done or are underway to overcome this  
problem: one can cite the {\it Vortex Methods} which are able to generate  
very thin scales, or {\it Large Eddy Simulation (LES)} and subgrid-scale  
techniques which separate the flow into large scales, that are  
explicitly computed, from the small scales, that are parametrized or  
statistically computed.  
  
\  
  
In that context, wavelet bases offer an intermediate decomposition   to  
suitably represent  
the intermittent spatial structure of turbulent flows, with only  
few degrees of freedom: this property is mainly due to the  
good localization, both in physical and frequency domains, of the  
basis functions.  
The wavelet decomposition was introduced in the beginning of the  
90s for the analysis of turbulent flows \cite{F92,M91,L93}. Wavelet based methods for the resolution of the  
Navier-Stokes equations appear later \cite{CP96, FS96, FK96, KGKFS98,  
GK00}. They have also been used to define LES-type methods such as the  
CVS method \cite{FS01}. Most of the  
works cited below use a Galerkin or a Petrov-Galerkin approach for the 2D vorticity  
formulation with periodic boundary  
conditions. However, if we want to turn to the 3D case, with non  
periodic boundary conditions, these approaches are no more available.

An alternative was at the same period firstly considered by K. Urban 
and after 
investigated by several authors: 
they proposed to use the {\it divergence-free wavelet bases}   
originally designed by Lemari\'e-Rieusset  
\cite{L92}.  Divergence-free wavelet vectors have been implemented and  
used to analyze 2D turbulence flows \cite{AUDRL02,KKR00,U02}, as well as to compute  
the 2D/3D Stokes solution for the  driven cavity problem  
\cite{U94,U96}. Seeing that divergence-free wavelets are  
constructed from standard compactly supported biorthogonal wavelet bases, they allow to incorporate  
boundary conditions in their construction \cite{DKU96,U00}.

This research direction is of great interest, since divergence-free 
wavelets provide with bases suitable to represent the  
incompressible Navier-Stokes solution, in two and three dimension. Our 
objective is now to investigate their feasibility and amenability for 
such problem. The first point lies in avoiding the pressure by 
projecting the equations  
onto the  
space of divergence-free vectors. This (orthogonal) projection is the 
well-known Leray projector, and it can be computed explicitly in 
Fourier space, for periodic boundary conditions.  
Unfortunately, as already noted by K. Urban \cite{U00}, 
if we want  to explicit the Leray operator in terms of  
divergence-free wavelets,  since they form   
biorthogonal bases (and not orthogonal), they would not give rise, in a  
simple way, to the  
orthogonal projection onto the space of divergence-free vectors.

\  
  
Nevertheless, we propose in the present paper to investigate the  
use of divergence-free wavelets for the simulation of turbulent flows.  
Firstly, we remind the basic ingredients of the theory of compactly  
supported  divergence-free wavelet vectors, developed by Lemari\'e-Rieusset  
\cite{L92}. In  section 3, we present in detail the bases we proposed to  
implement in space dimensions 2 and 3. We will see that the choice of the complement wavelet basis is not  
unique from this construction, and it induces the values of 
divergence-free coefficients, for compressible flows. We discuss the algorithmic 
implementation of divergence-free wavelet coefficients in
dimensions 2 and 3, leading to fast  
algorithms (in O($N$) operations  
where $N$ is the number of grid points).  
  
Section 4 is devoted to the Hodge decomposition of a compressible  
field, in a wavelet formulation: the method we present uses both the 
biorthogonal projectors 
on divergence-free, and  on curl-free wavelets: our method is  
an iterative procedure, and  
we will experimentally prove that it converges.  
The last section presents numerical tests to validate our approach:  
nonlinear compression of 2D and 3D incompressible turbulent flows, and  
Hodge decomposition of well chosen examples, such as the nonlinear 
term of the Navier-Stokes equations.  
  
  
\section{Theory of divergence-free wavelet bases}\label{basic}  
  
In this section, we review briefly the relevant properties of wavelet  
bases, that will be used for the construction of divergence-free  
wavelets. Compactly supported divergence-free vector wavelets were  
originally designed by Lemari\'e-Rieusset, in the context of  
biorthogonal Multiresolution Analyses. We illustrate the construction  
with the explicit example of splines of degree 1 and 2. For more  
details, we refer to \cite{L92,Dau92,KL98,U96}.  
  
\subsection{Multiresolution Analyses (MRA)}  
  
Multiresolution Analyses (MRA) are approximation spaces allowing the  
construction of wavelet bases, introduced by S. Mallat \cite{M98}. We  
begin with the one-dimensional case of functions defined on the real line.  
  
\  
  
\noindent  
\textbf{Definition (MRA)}: {\it  
A Multiresolution Analysis of $L^2(\mathbb{R})$ is a sequence of closed subspaces $(V_j)_{j\in \mathbb{Z}}$  
verifying:  
\begin{description}  
\item[(1)] $\forall j,~~V_j \subset V_{j+1}, \quad \bigcap_{j\in \mathbb{Z}} V_j=\{0\}, \quad  
\bigcup_{j\in \mathbb{Z}} V_j \mbox{ is dense in }  L^2(\mathbb{R})$ \qquad  
\item[(2)] {\rm (Dilation invariance)} \quad $f \in V_j \quad \iff \quad f(2 .) \in V_{j+1}$ \qquad  
  
\item[(3)] {\rm (Shift-invariance)} \quad There exists a function 
$\phi\in V_0 $ such that the family $\{\phi(.-k)~;~k\in\mathbb{Z}\}$ form a  
(Riesz) basis of $V_0$. \\  
 \end{description} }

\noindent 
The function $\phi$ in (3) is called a {\it scaling function} of the MRA.\\ 
Here, $j$ denotes the level of refinement. By virtue of the dilation 
invariance property (2) above, we can deduce that each space  
$V_j$ is spanned by $\{\phi_{j,k}~;~k\in\Z\}$ where $\phi_{j,k}(x)=2^{j/2}\phi(2^j x-k)$.

\  
  
\noindent 
Wavelets appear as bases of complementary spaces  $W_j$:  
 \begin{equation}  
 \label{Wj}  
 V_{j+1}=V_j \oplus W_j   
 \end{equation}  
 where the sum is direct, but not necessarily orthogonal. In this context (called the  
  biorthogonal case), the choice of  
 spaces $W_j$ is not unique. The problem of constructing the spaces $W_j$ means to find a function  
$\psi$, called {\it wavelet} such that the system $\{\psi(.-k)~~;k\in\mathbb{Z}\}$ spans  
$W_0$. Repeated decomposition of $V_j$ yields the multiresolution analysis of $V_j$ with the wavelet spaces:  
$$V_{j}=V_0 \bigoplus_{\ell=0}^{j-1} W_{\ell} $$  
which leads, when $j\to +\infty$, to the wavelet decomposition of the whole space:  
$$L^2(\mathbb{R})=V_0 \bigoplus_{\ell=0}^{+\infty} W_{\ell} $$  
As a result, we can write any  
function $f\in L^2(\mathbb{R})$ in the basis  
$\{\phi_k,~\psi_{j,k}~;~j \geq 0,k \in \mathbb{Z}\}$,  
with $\phi_k=\phi\cdot -k)$ and $\quad \psi_{j,k}=2^{j/2}\psi(2^j\cdot-k)$:  
\begin{equation}  
\label{decomp}  
f=\sum_{k\in\mathbb{Z}}c_k~\mathbf{\phi}_k+\sum_{j\geq  
0}\sum_{k\in\mathbb{Z}}~d_{j,k}~\mathbf{\psi}_{j,k}  
\end{equation}  
  
\  
  
\noindent  
{\textbf{Dual bases}}:  
 Let a  pair $(\phi,\psi)$ of scaling function and wavelet, arising from
a Multiresolution Analysis, be given, then we can associate a unique dual pair $(\phi^*,\psi^*)$,  
such that the following biorthogonality (in space $L^2$) relations are fulfilled: for all $k\in\mathbb{Z}$ and $j\geq0$,  
\begin{equation}  
\label{dual}  
<\phi|\phi^*_{k}>=\delta_{k,0}, ~ <\phi|\psi^*_{j,k}>=0, ~ 
<\psi|\psi^*_{j,k}>=\delta_{j,0}\delta_{k,0}, ~ <\psi|\phi^*_{k}>=0  
\end{equation}  
Moreover the dual scaling functions $\phi^*_{k}$ and the dual wavelets $\psi^*_{j,k}$ have the same  
structure as above: $\phi^*_k=\phi^*(\cdot -k)$ and $\psi^*_{j,k}=2^{j/2}\psi^*(2^j\cdot-k)$.  
  
\  
  
\noindent  
\textbf{Scaling equations and filter design}: Since the function $\frac{1}{\sqrt{2}}\phi(\frac{\cdot}{2})$ lives in $V_0$,  
there exists a sequence $(h_k)$ (also called the {\it low pass filter}) verifying:  
\begin{equation}  
\label{hk}  
\frac{1}{\sqrt{2}}\phi(\frac{x}{2})=\sum_{k\in\mathbb{Z}}h_k ~\phi(x-k)  
\end{equation}  
By applying the Fourier transform\footnote{The Fourier transform of a function$f$ is defined by  
$\hat f (\xi)=\int_{-\infty}^{+\infty} f(x) ~e^{-ix\xi} dx$}, (\ref{hk}) rewrites:  
$$  
\hat{\phi}(2\xi)=m_0(\xi)\hat{\phi}(\xi)  
$$  
where $m_0(\xi)=\frac{1}{\sqrt{2}}\sum_{k\in\mathbb{Z}}h_k e^{-ik\xi}$ is the transfer function of the filter $(h_k)$.\\  
 
\noindent   
Again, because of $W_{-1}\subset V_0$, the wavelet  satisfies a two-scale equation:  
\begin{equation}  
\label{gk}  
\frac{1}{\sqrt{2}}\psi(\frac{x}{2})=\sum_{k\in\mathbb{Z}}g_k ~\phi(x-k)  
\end{equation}  
where the coefficients $(g_k)$ are called the {\it high pass filter}. Again the Fourier transform  of  
$\psi$ expresses with the transfer function $n_0$ of filter $g_k$ as:  
$$  
\hat{\psi}(2\xi)=n_0(\xi)\hat{\phi}(\xi)  
$$  
In the same way, the dual functions satisfy scaling equations:  
\begin{equation}  
\label{hk*}  
\begin{array}{ccc}  
\frac{1}{\sqrt{2}}\phi^*(\frac{x}{2})=\sum_{k\in\mathbb{Z}}h^*_k ~\phi^*(x-k), & & \hat{\phi^*}(2\xi)=m^*_0(\xi)\hat{\phi^*}(\xi)\\  
\frac{1}{\sqrt{2}}\psi^*(\frac{x}{2})=\sum_{k\in\mathbb{Z}}g^*_k ~\phi^*(x-k), & & \hat{\psi^*}(2\xi)=n^*_0(\xi)\hat{\phi^*}(\xi)  
\end{array}  
\end{equation}  
Following \cite{KL98,Dau92},  the biorthogonality conditions (\ref{dual}) for the scaling functions imply:  
$$  
m_0(\xi)\overline{m^*_0(\xi)}+m_0(\xi+\pi)\overline{m^*_0(\xi+\pi)}=1\quad  
$$  
while one can choose, for example, as  transfer functions for the associated wavelets: 
$$  
\begin{array}{ccc}  
n_0(\xi)=e^{-i\xi}~\overline{m^*_0(\xi+\pi)}, & & n^*_0(\xi)=e^{-i\xi}~\overline{m_0(\xi+\pi)}  
\end{array}  
$$  
which corresponds to: 
$$g_k=(-1)^{1-k}~ h^*_{1-k},~~~g^*_k=(-1)^{1-k} ~h_{1-k},~~\forall k $$  
In practice, the filter coefficients  $h_k$ and $g_k$ are all what is needed to compute the wavelet  
decomposition (\ref{decomp}) of a given function. Notice that these filters are finite if and only if the  
functions $\psi$ and $\phi$ are compactly supported.  
  
\  
  
\begin{figure} [!hb] 
\begin{center}  
\begin{tabular}{llll}  
\hspace{-0.8cm}\includegraphics[angle=-90,scale=0.15]{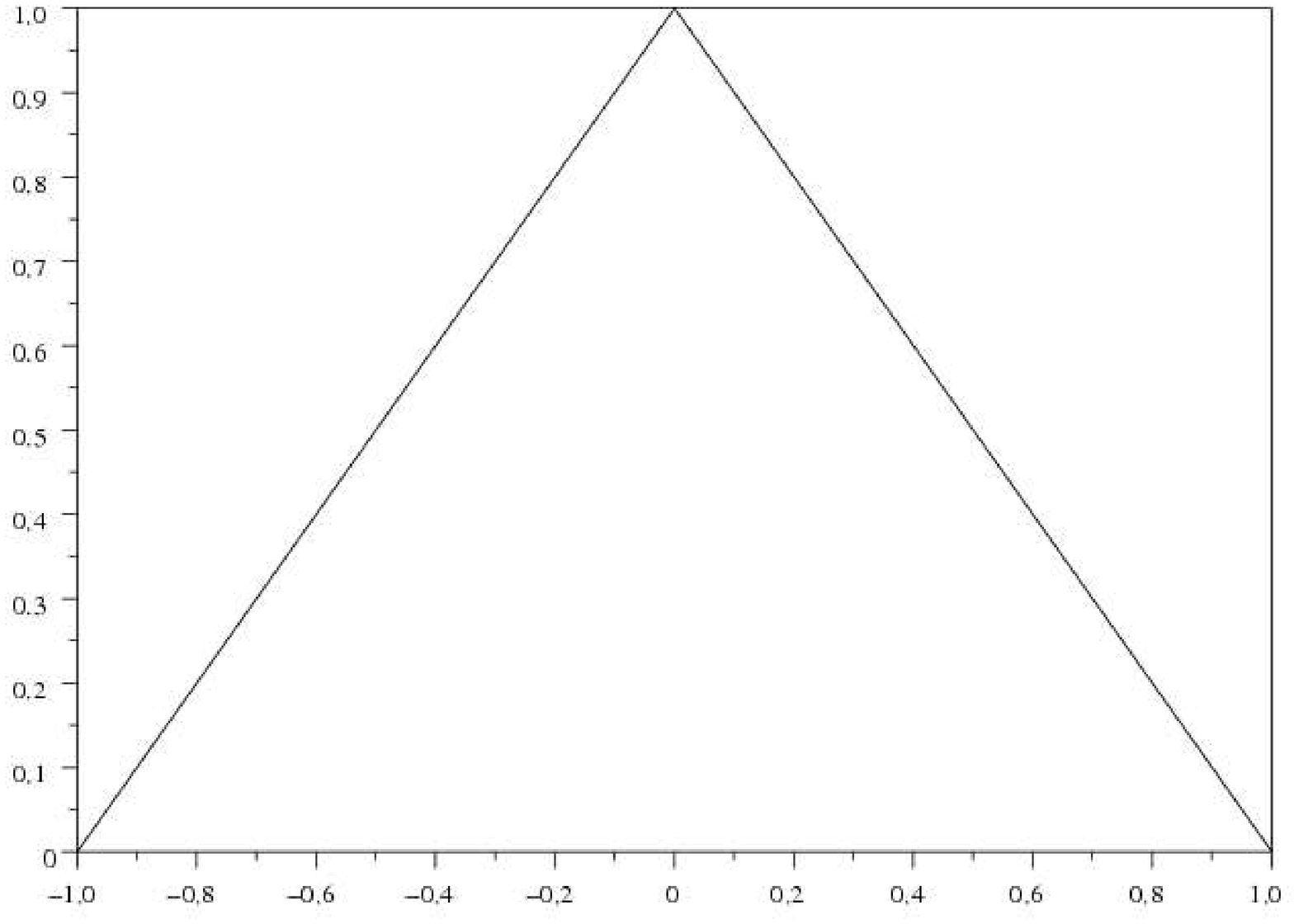} & \includegraphics[angle=-90,scale=0.15]{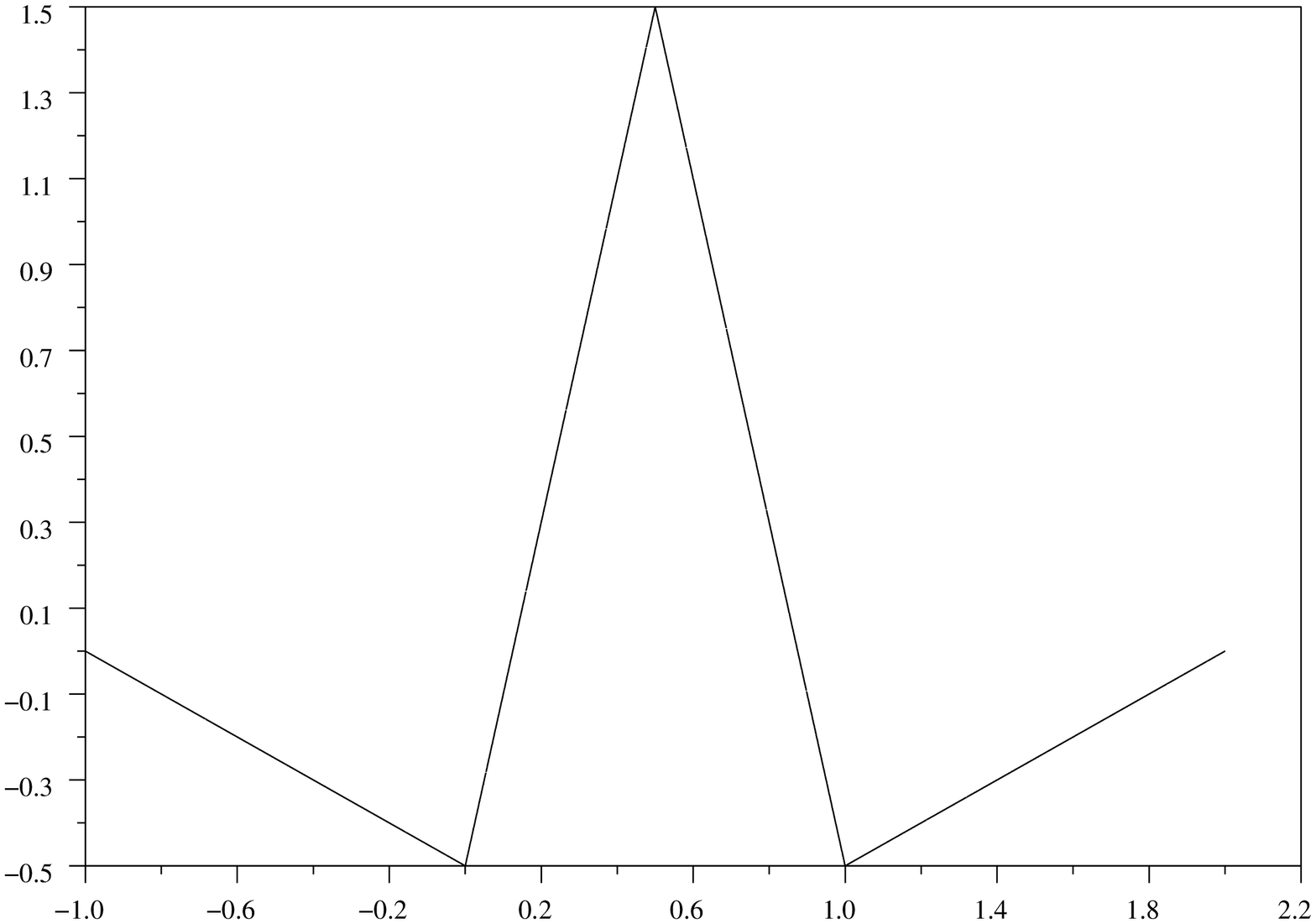} &  
\includegraphics[angle=-90,scale=0.15]{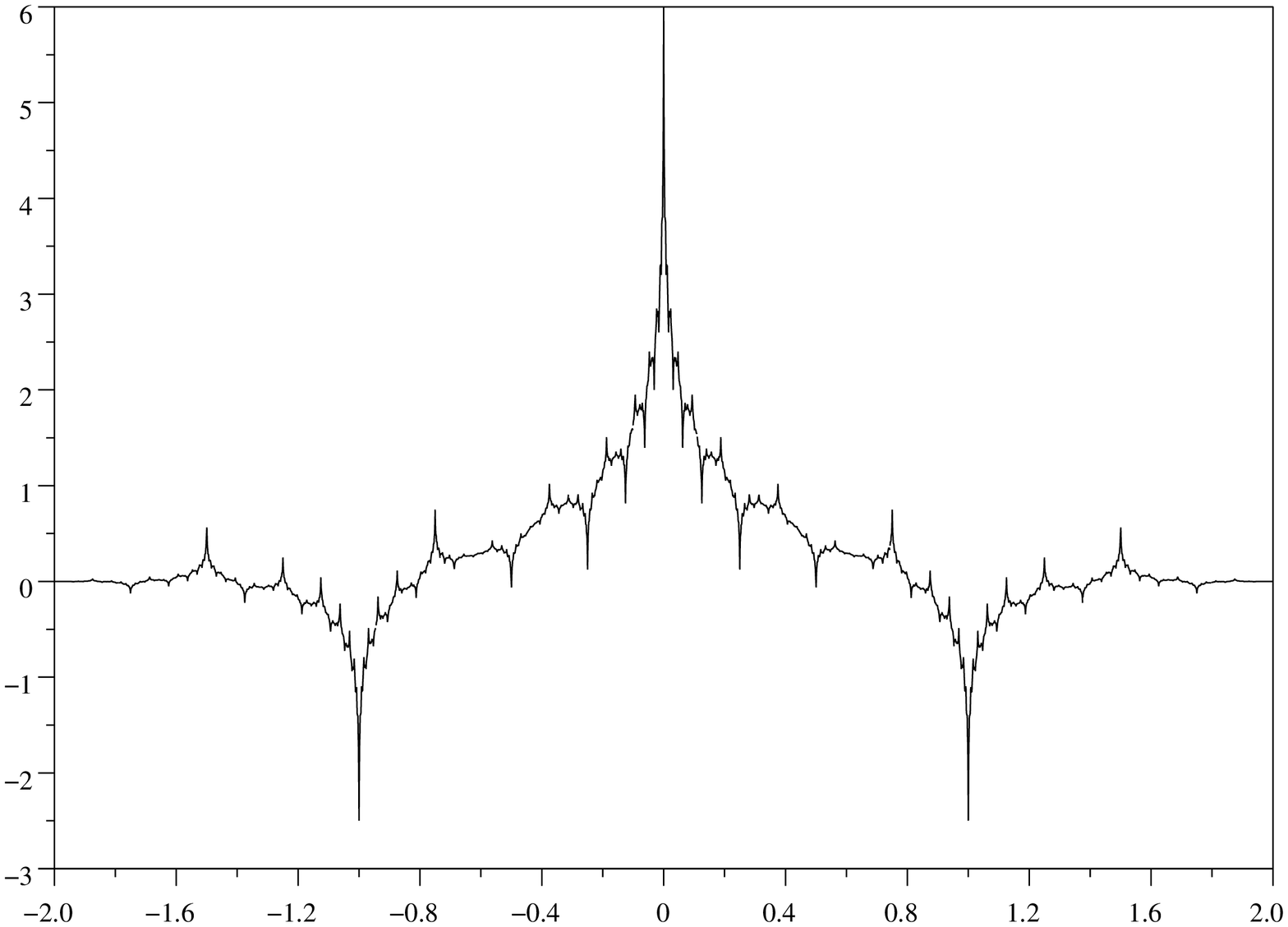} & \includegraphics[angle=-90,scale=0.15]{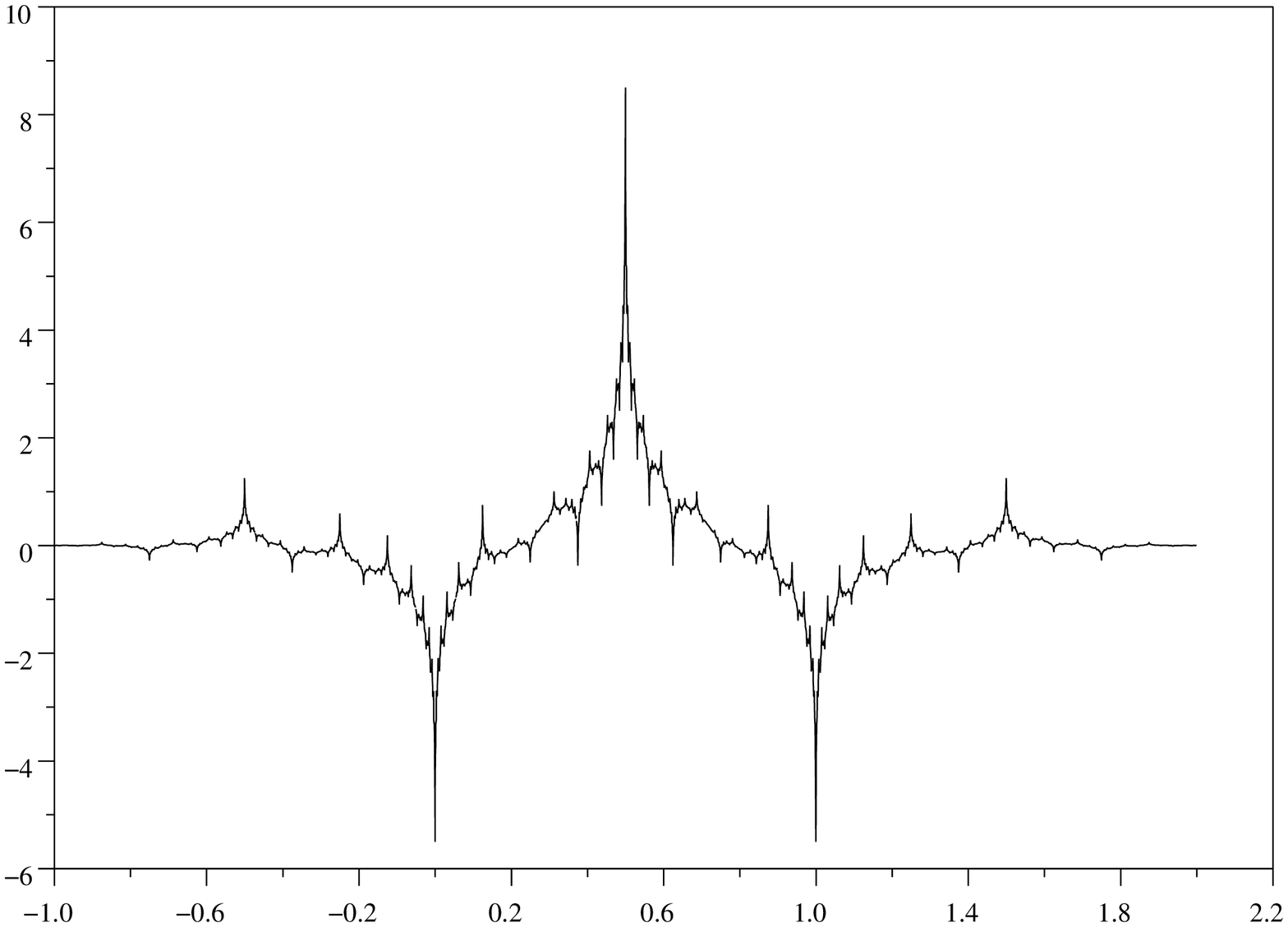}  
\end{tabular}  
\end{center}   
\caption{\label{spline1} From left to right: the scaling function $\phi$ with its associated symmetric wavelet  
with shortest support, and their duals: the dual scaling function $\phi^*$ and the dual wavelet $\psi^*$.}  
 
\end{figure}

\noindent  
\textbf{Example: symmetric biorthogonal splines of degree 1}\\  
A simple example for spaces $V_j$ are the spaces of continuous functions, which are  
piecewise linear on the intervals $[k2^{-j},(k+1)2^{-j}]$, for  
$k\in\Z$. In this case we can choose as scaling function the hat  
function  
$\phi(x)=\max(0,1-|x|)$. Its transfer function is given by   
\begin{equation}  
\label{m0} 
m_0(\xi)=e^{i\xi}~\left(\frac{1+e^{-i\xi}}{2}\right)^2  
\end{equation}  
The shortest even dual scaling function associated with $\phi$ is associated with the filter:  
\begin{equation}  
\label{m0*} 
m^*_0(\xi)=e^{i\xi}~\left(\frac{1+e^{-i\xi}}{2}\right)^2(2-\cos \xi)  
\end{equation}  
The corresponding values of filters $(h_k)$ and $(h^*_k)$ are given in table \ref{filtres}.  
\noindent  
Figure \ref{spline1} displays the scaling functions and their associated wavelets in this case.

\subsection{Decomposition-recomposition algorithm and useful example}  
\label{splines}
  
In the context of a biorthogonal Multiresolution Analysis, the wavelet decomposition of a given function  
$f\in L^2(\R)$ (equation (\ref{decomp})), is obtained through the now well   
known Fast Wavelet Transform \cite{M98}. We briefly review here the formula that will be useful for the  
following.  
  
In practice we begin with an approximation $f_J$ of $f$ in some space $V_J$ of the MRA. This approximation  
may be the oblique projection of $f$ onto $V_J$ following the 
direction perpendicular to $V^*_J$ (also called the  
biorthogonal projection on $V_J$); but usually $f_J$ means an interpolating function of $f$, associated to  
nodes $\{k2^{-J}~;~k\in \Z\}$ (see in last section some examples of procedure).

This  approximation $f_J$ of  function $f$ may be expanded in terms of 
the scaling basis  
$\phi_{J,k}=2^{J/2}\phi(2^J x -k)$ of $V_J$:  
$$  
f_J (x)=2^{J/2}\sum_{k=-\infty}^{+\infty}c_{J,k}~\phi(2^Jx-k)  
$$  
The wavelet decomposition of $f_J$ corresponds to a truncated sum in equation (\ref{decomp}) up to the  
level $J-1$, meaning that  
$2^{-J}$ is the finest scale of approximation:  
\begin{equation}  
\label{decompJ}  
f_J=\sum_{k\in\mathbb{Z}}c_k~\mathbf{\phi}_k+\sum_{j=0}^{J-1}\sum_{k\in\mathbb{Z}}~d_{j,k}~\mathbf{\psi}_{j,k}  
\end{equation}  
The wavelet coefficients $d_{j,k}$ are then computed recursively on the level $j$, from $j=J-1$ to $j=1$  
using the decomposition of spaces (\ref{Wj}):  
$$  
f_{j+1}=\sum_{k=-\infty}^{+\infty}c_{j+1,k}~\phi_{j+1,k} =\sum_{k=-\infty}^{+\infty}c_{j,k}~\phi_{j,k}  
+\sum_{k=-\infty}^{+\infty}d_{j,k}~\psi_{j,k}  
$$  
(Here $f_{j+1}$ denotes the biorthogonal projection of $f_J$ onto $V_{j+1}$).  
By biorthogonality (\ref{dual}) one has:  
$$  
c_{j,k}=<f_{j+1}|\phi^*_{j,k}>~,~~d_{j,k}=<f_{j+1}|\psi^*_{j,k}>~,~~c_{j+1,k}=<f_{j+1}|\phi^*_{j+1,k}>  
$$  
which yield the decomposition formula: for all $j=0,\dots J-1$, 
$$\left\{  
\begin{array}{c}  
c_{j,k}=\sum_{\ell}h^*_{\ell}~c_{j+1,\ell+2k}\\  
\vspace{0.1cm}  
\\  
d_{j,k}=\sum_{\ell}g^*_{\ell}~c_{j+1,\ell+2k}  
\end{array}\right.  
$$  
where the filters $h^*_{\ell}$ and $g^*_{\ell}$ arise from the scaling equations (\ref{hk*})  
of the dual basis functions.\\  
In the same way, we obtain the reconstruction formula:  
$$  
c_{j+1,k}=\sum_{\ell}\left(h_{k-2\ell}~c_{j,\ell}  
+g_{k-2\ell}~d_{j,\ell}\right)  
$$  
where $h_{k}$ and $g_k$ are the filters provided by the scaling equations (\ref{hk}, \ref{gk})  
of the primal basis functions.\\  
  
The computing cost for the whole wavelet decomposition (\ref{decompJ}) (as well as for the recomposition)  
is about $C 2^J$ operations, where $2^J$ is the number of point values $f(k2^{-J})$ we start with, and $C$  
means the length of the filters ($h^*_k$ for the decomposition, $h_k$ for the synthesis).  
  
\  
  
\begin{figure}[!h]  
\begin{center}  
\begin{tabular}{cc}  
\includegraphics[width=3.5cm,height=7cm,angle=-90]{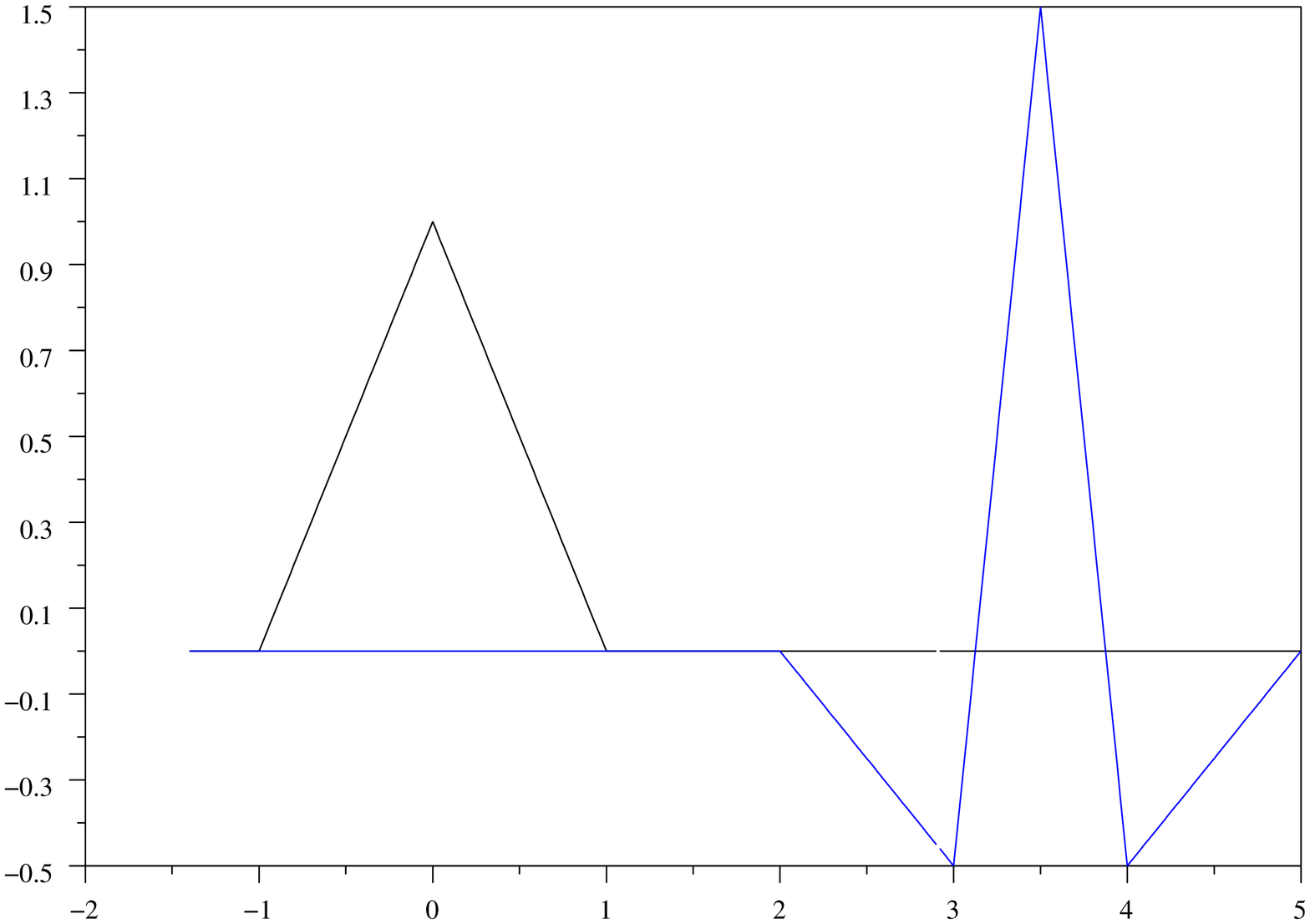} & 
\includegraphics[width=3.5cm,height=7cm,angle=-90]{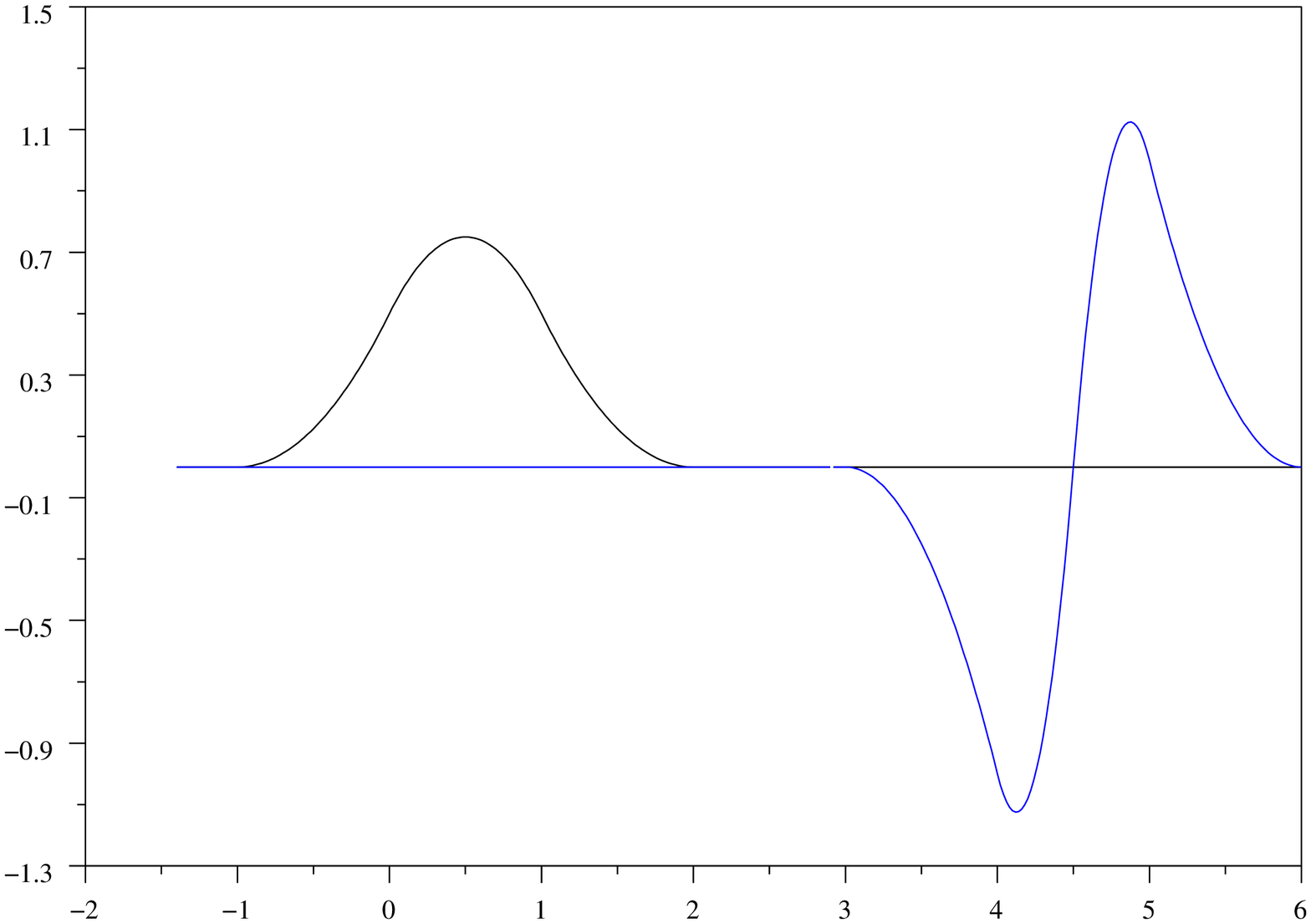}\\  
$\phi_0\qquad \qquad \qquad \psi_0$ & $\phi_1\qquad \qquad \qquad \psi_1$  
\end{tabular}  
\end{center} 
\caption{\label{V0V1} Scaling functions and associated even and odd wavelets with shortest support, for splines of degree 1 (left) and 2 (right).}  
 \end{figure}

\noindent  
\textbf{Example:  spline wavelets of degree 1 and 2}: Biorthogonal splines provide with wavelet bases which are regular,  
compactly-supported and easy to implement.   
The scaling functions of the associated MRA are standard B-spline bases, and the wavelets are constructed easily, by linear  
combinations of translated B-splines. We focus here on two examples of wavelet bases, which will be useful  
for the construction of divergence-free wavelets: splines of degree 1 ($V^0_j$ MRA spaces) and splines of  
degree 2 ($V^1_j$ MRA spaces). In both cases we draw the scaling functions and the associated wavelets  
with shortest support (figure \ref{V0V1}).  
In both cases the filters are easy to compute (see \cite{KL98,Dau92}). Since the support of basis function  
is very short, the filters have a few non zero coefficients. The values of the decomposition and  
reconstruction filters are given on Table \ref{filtres}.  
 
\begin{table}[ht] 
\begin{center}  
$$ 
\begin{array}{|c|c|c|c|c|c|c|}  
\hline  
\ell & -2 & -1 & 0 & 1 & 2 & 3 \\  
\hline  
\frac{1}{\sqrt{2}}h^{*0}_{\ell} & -1/8 & 1/4 & 3/4 & 1/4 & -1/8 &0\\  
\hline  
\frac{1}{\sqrt{2}}g^{*0}_{\ell} & 0 & 0 & -1/4 & 1/2 & -1/4&0 \\  
\hline  
\frac{1}{\sqrt{2}}h_{\ell}^{0}&0 & 1/4 & 1/2 & 1/4 & 0 & 0 \\  
\hline  
\frac{1}{\sqrt{2}}g_{\ell}^{0}&0 & -1/8 & -1/4 & 3/4 & -1/4 & -1/8 \\  
\hline  
\end{array}  
~~~~~~  
\begin{array}{|c|c|c|c|c|}  
\hline  
\ell &-1 & 0 & 1 & 2 \\  
\hline  
\frac{1}{\sqrt{2}}h^{*1}_{\ell}  & -1/4 & 3/4 & 3/4 & -1/4\\   
\hline  
\frac{1}{\sqrt{2}}g^{*1}_{\ell}  & 1/8 & -3/8 & 3/8 & -1/8\\   
\hline  
\frac{1}{\sqrt{2}}h_l^{1} & 1/8 & 3/8 & 3/8 & 1/8 \\  
\hline  
\frac{1}{\sqrt{2}}g_l^{1} & -1/4 & -3/4 & 3/4 & 1/4  \\  
\hline  
\end{array}  
$$  
\end{center}  
\caption{\label{filtres} Decomposition filter ($h^*_k,~g^*_k$) and reconstruction filter ($h_k,~g_k$)  
coefficients, associated to piecewise linear splines (left) and piecewise quadratic splines (right),  
verifying (\ref{derb}) (see hereafter) with  
shortest supports.}  
\end{table}  
 
\subsection{Multivariate wavelets}  
\label{multi} 
  
\begin{figure}[!hb] 
\input{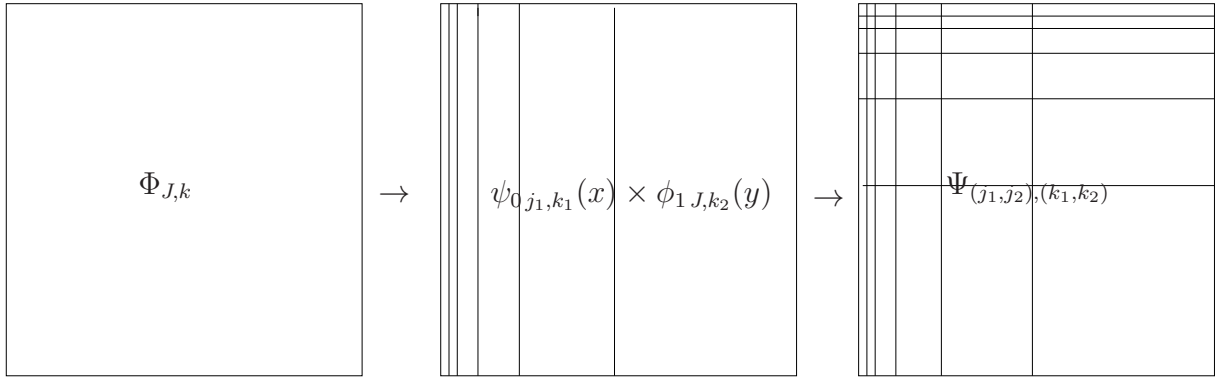}  
\caption{\label{fig:aniso2D} Anisotropic 2D wavelet transform.}  
\end{figure}

The above considerations can be extended to multi-D. The simplest way to obtain multivariate wavelets is to employ anisotropic or  
isotropic tensor products of one-dimensional functions.\\  
To be more precise, we focus on the two dimensional case: let $(V^0_j)$ and $(V^1_j)$ two   
multiresolution analyses of $L^2(\R)$ be given, associated with scaling functions and wavelets $(\phi_0,  
\psi_0)$ and $(\phi_1, \psi_1)$; the two dimensional tensor product space $V^0_J\otimes V^1_J$ is  
generated by the scaling basis $\{\phi_{0,J,k_1}(x)\phi_{1,J,k_2}(y);(k_1,k_2)\in\Z^2\}$, where  
$\phi_{0,J,k_1}(x)=2^{J/2}\phi_0(2^Jx-k_1)$ and similarly for $\phi_{1,J,k_2}$. Then each function $f_J$  
in $V^0_J\otimes V^1_J$ can be written:  
\begin{equation}  
\label{UJ}  
f_J(x,y)=\sum_{k_1=-\infty}^{\infty}\sum_{k_2=-\infty}^{\infty}c_{J,k_1,k_2}~2^J\phi_0(2^Jx-k_1)~\phi_1(2^Jy-k_2)  
\end{equation}  
The anisotropic 2D wavelets are constructed with  tensor products of wavelets at different scales  
$\{\psi_{0,j_1,k_1}(x)\psi_{1,j_2,k_2}(y)\}$. For certain choices of $j_1, j_2$, the support of the  
functions may be very lengthened. In this case the wavelet decomposition of $f_J$ writes:  
\begin{eqnarray*}  
f_J(x,y) & =&\sum_{(k_1,k_2)\in\Z^2} c_{k_1,k_2}~\phi_0(x-k_1)~\phi_1(y-k_2)\\  
&&+\sum_{j_1=0}^{J-1}\sum_{j_2=0}^{J-1}2^{(j_1+j_2)/2} \sum_{(k_1,k_2)\in\Z^2} d_{j_1,j_2,k_1,k_2}~\psi_0(2^{j_1}x-k_1)~\psi_1(2^{j_2}y-k_2)  
\end{eqnarray*}  
The anisotropic decomposition is the most easy way to compute a multi-dimensional wavelet transform, as it  
corresponds to apply  one-dimensional wavelet decompositions in each direction. In the 2D case, this is  
schematized in figure \ref{fig:aniso2D}.   
  
\begin{figure}[!h]  
\input{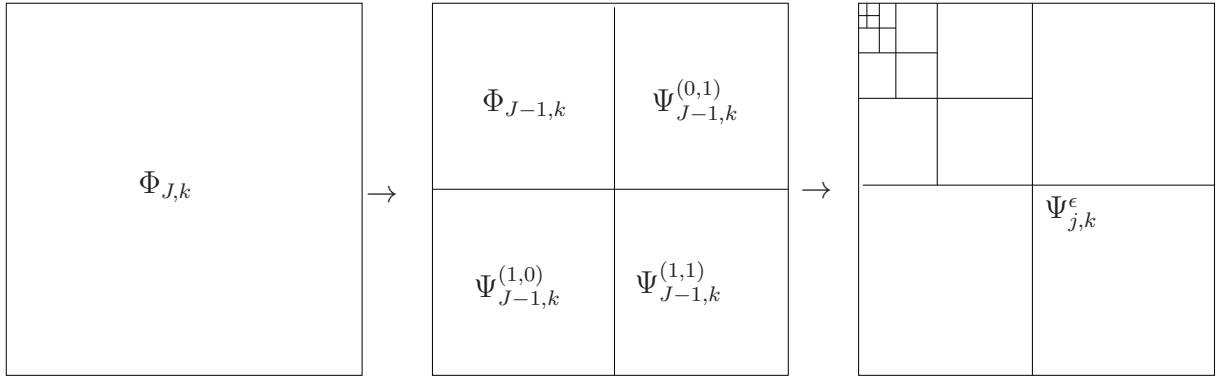}  
\caption{\label{fig:iso2D} Isotropic 2D wavelet transform.}  
\end{figure}

In the isotropic case, the 2D wavelets are obtained through tensor products of wavelets  
and scaling functions or wavelets at the {\it same} scale. This produces the following  
decomposition for $f_J$:  
\begin{eqnarray*}  
f_J(x,y) & =&\sum_{(k_1,k_2)\in\Z^2} c_{k_1,k_2}~\phi_0(x-k_1)~\phi_1(y-k_2)\\  
&&+\sum_{j_=0}^{J-1}\left( \sum_{k_1,k_2}d_{j,k_1,k_2}^{(1,0)}~\psi_{0,j,k_1}(x)~\phi_{1,j,k_2}(y)+  
 \sum_{k_1,k_2}d_{j,k_1,k_2}^{(1,1)}~\psi_{0,j,k_1}(x)~\psi_{1,j,k_2}(y)\right.\\  
 && \left.\quad +\sum_{k_1,k_2}d_{j,k_1,k_2}^{(0,1)}~\phi_{0,j,k_1}(x)~\psi_{1,j,k_2}(y)\right)  
\end{eqnarray*}  
As one can see, this decomposition involves three kinds of wavelets, 
one following the direction $x$: 
$\Psi^{(1,0)}(x,y)=\psi_{0}(x)~\phi_{1}(y)$, one  
following the direction $y$: 
$\Psi^{(0,1)}(x,y)=\phi_{0}(x)~\psi_{1}(y)$ 
 and one in both directions: 
$\Psi^{(1,1)}(x,y)=\psi_{0}(x)~\psi_{1}(y)$. The interest of this basis remains in the  
fact that the size of their support is proportional to $2^{-j}$ in each direction, {\it i.e.} the basis  
functions are rather isotropic. The principle of the associated decomposition algorithm is illustrated by  
figure \ref{fig:iso2D}.  
  
\subsection{Theoretical ground of the divergence-free wavelet vectors}  
 \label{theory} 
  
Let introduce $$\vec{\textbf{H}}_{\textrm{div}}(\R^n)=\{f\in (L^2(\mathbb{R}^n))^n / \textrm{div}f\in L^2(\mathbb{R}^n)  
,\quad \textrm{div}f=0\}$$  
 the space of divergence-free vector functions in $\R^n$.\\  
The construction of compactly divergence-free wavelets in $(L^2(\mathbb{R}^n))^n$, which will correspond  
to Riesz bases of $\vec{\textbf{H}}_{\textrm{div}}(\R^n)$,  was originally derived by Lemari\'e-Rieusset  
\cite{L92}: it is based on the following proposition, which relates two different multiresolution analyses  
of $L^2(\mathbb{R})$ by differentiation and integration:\\  
  
\noindent \textbf{Proposition}\label{theoreme-Lemarie}: {\it 
Let $(V^0_j)$ a  
one-dimensional MRA with a derivable  
scaling function $\phi_1$ and a  
wavelet $\psi_1$, be given. Then, we can build a MRA $(V^1_j)$ with a scaling function $\phi_0$  
and a wavelet $\psi_0$ verifying:  
$$  
V_0^0=\textrm{span}\{\phi_0(x-k),k\in\Z\} \qquad  
V_0^1=\textrm{span}\{\phi_1(x-k),k\in\Z\}  
$$  
and  
\begin{equation}  
\label{derb}  
\phi_1'(x)=\phi_0(x)-\phi_0(x-1) \qquad 
\psi_1'(x)=4~\psi_0(x)  
\end{equation}  
For the refinement polynomials it can be traduced by:  
$$  
m_0(\xi)=\frac{2}{1+e^{-i\xi}}~m_1(\xi)  
$$  
Equation (\ref{derb}) rewrites for the dual functions $\phi_0^*$, $\psi_0^*$, $\phi_1^*$, and 
$\psi_1^*$:  
\begin{equation}  
\label{derb*}  
{\phi_0^*}'(x)=\phi_1^*(x+1)-\phi_1^*(x)\\  
 {\psi_0^*}'(x)=-4~\psi_1^*(x)  
\end{equation}  
which induces: 
$$  
m_0^*(\xi)=\frac{1+e^{i\xi}}{2}~m_1^*(\xi)  
$$  
}  
 
\noindent  
\textbf{Example:} As an example of functions fulfilling the above proposition, one shall cite the piecewise  
linear spline functions $\phi_0$, $\psi_0$, associated with the piecewise quadratic spline functions  
$\phi_1$, $\psi_1$ introduced in section 2.2, and plotted on figure \ref{V0V1}.  
  
\  
  
Then, the divergence-free wavelets are explicitly constructed by combining suitable tensor products of  
these functions. For instance, in the 2D case we may have the following basis \cite{L92}:  
  
\  
  
\noindent  
\textbf{Example:} The 2D {\it divergence-free vector scaling function} 
takes the form:  
$$  
\begin{array}{llll}  
\Phi_{\textrm{div}}(x_1,x_2) =\left| \begin{array}{l} \phi_1(x_1){\phi_1}'(x_2) \\  
- {\phi_1}'(x_1)\phi_1(x_2) \end{array} \right. 
& =& 
\left| \begin{array}{l} \phi_1(x_1)~[\phi_0(x_2)-\phi_0(x_2-1)] \\  
- [\phi_0(x_1)-\phi_0(x_1-1)]  ~\phi_1(x_2) \end{array} \right. 
\end{array}  
$$  
and the corresponding {\it isotropic} vector wavelets are given by the system:  
$$  
\begin{array}{ll}  
\Psi_{\textrm{div}}^{(1,0)}(x_1,x_2) =\left| \begin{array}{l} -\frac{1}{4}\psi_1(x_1)[\phi_0(x_2)-\phi_0(x_2-1)] \\  
\psi_0(x_1)\phi_1(x_2) \end{array} \right. \\  
\vspace{1mm} & \vspace{1mm} \\  
\Psi_{\textrm{div}}^{(0,1)}(x_1,x_2) =\left| \begin{array}{l} \phi_1(x_1)\psi_0(x_2) \\  
-\frac{1}{4}[\phi_0(x_1)-\phi_0(x_1-1)]\psi_1(x_2) \end{array} \right. \\  
\vspace{1mm} & \vspace{1mm} \\  
\Psi_{\textrm{div}}^{(1,1)}(x_1,x_2) =\left| \begin{array}{l} \psi_1(x_1)\psi_0(x_2) \\  
-\psi_0(x_1)\psi_1(x_2)\end{array}\right.   
\end{array}  
$$  
  
\begin{figure}[!h]  
\begin{center}  
\begin{tabular}{ccc}  
\includegraphics[width=5cm,height=5cm,angle=-90]{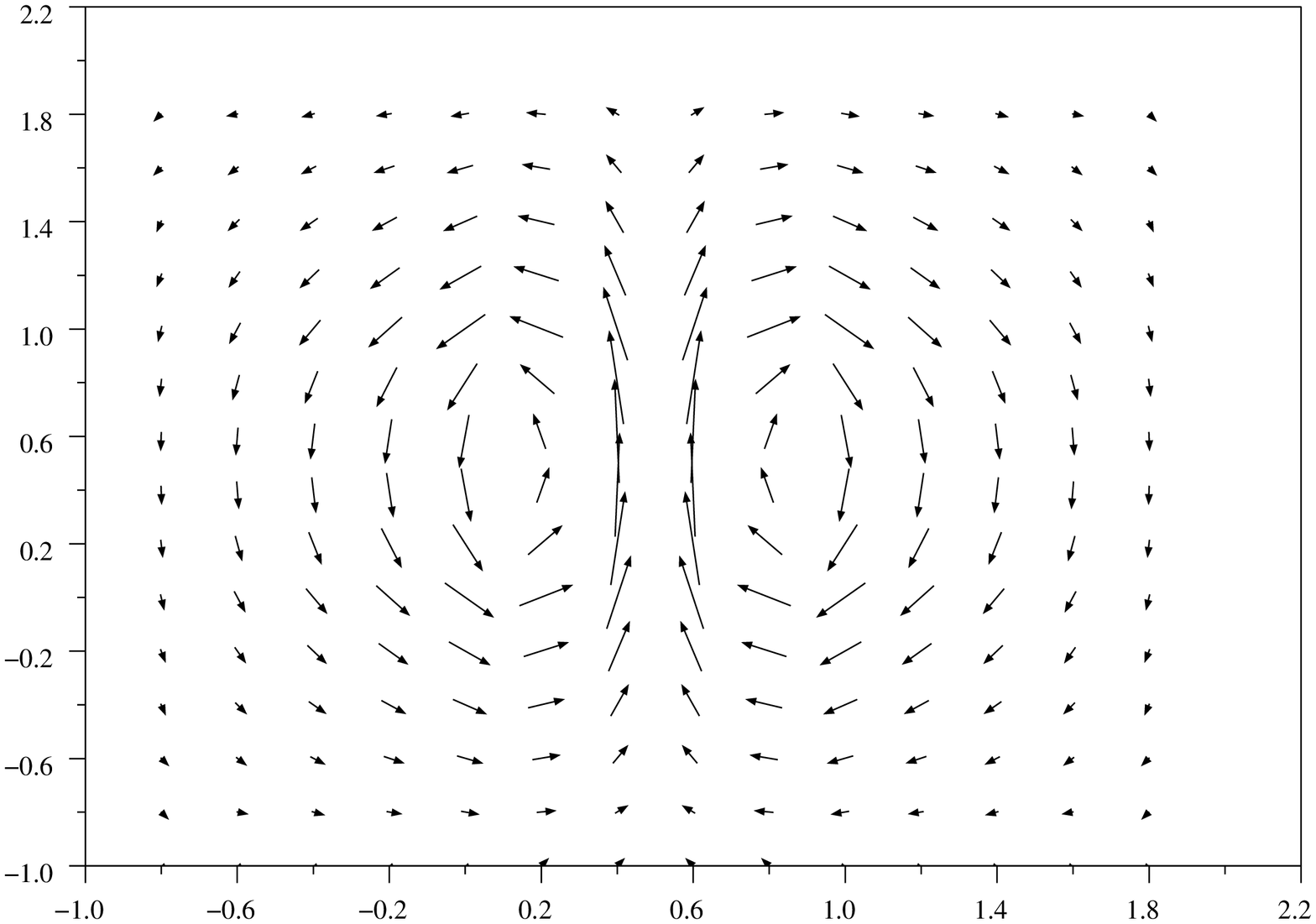} &  
\includegraphics[width=5cm,height=5cm,angle=-90]{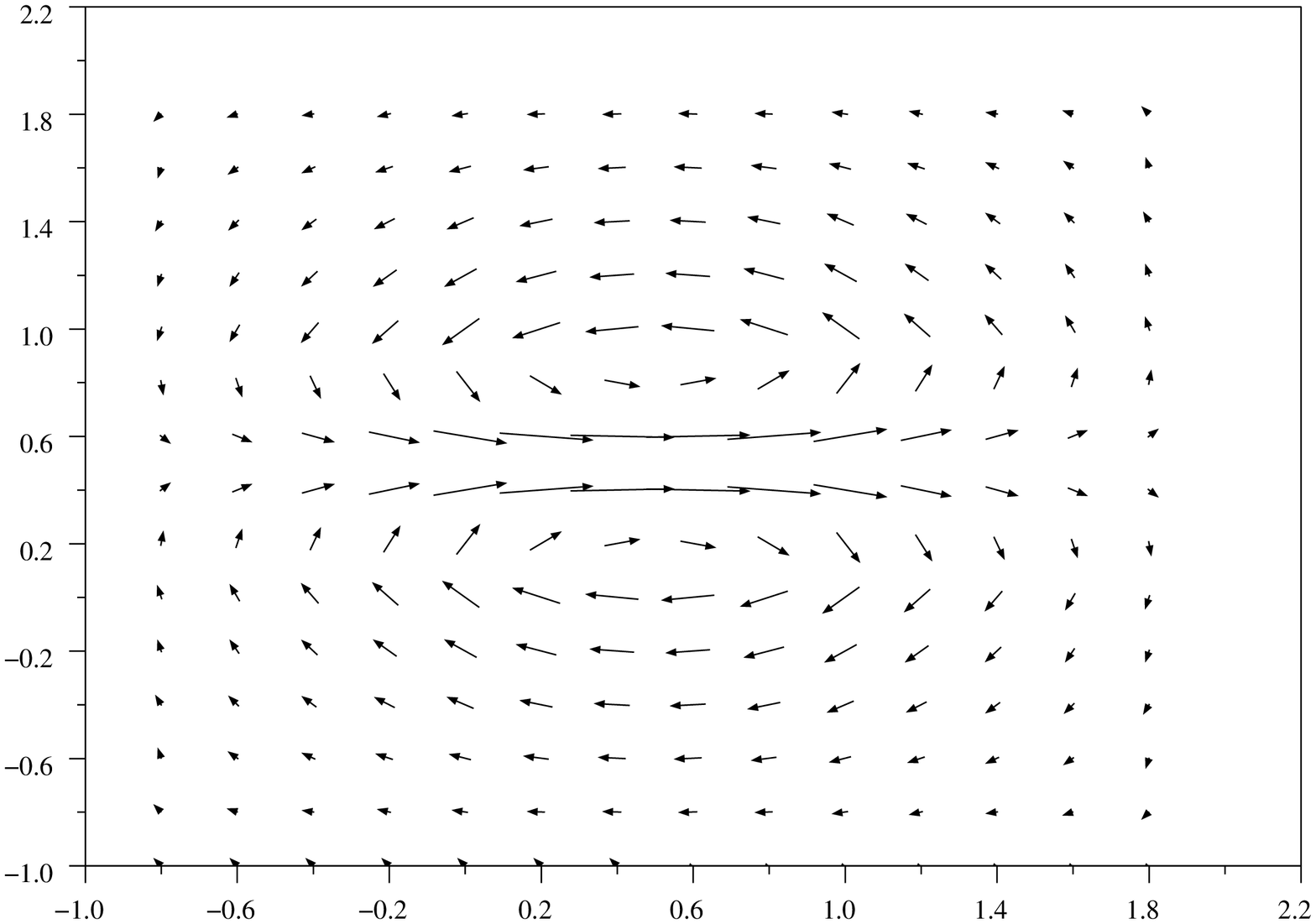} &  
\includegraphics[width=5cm,height=5cm,angle=-90]{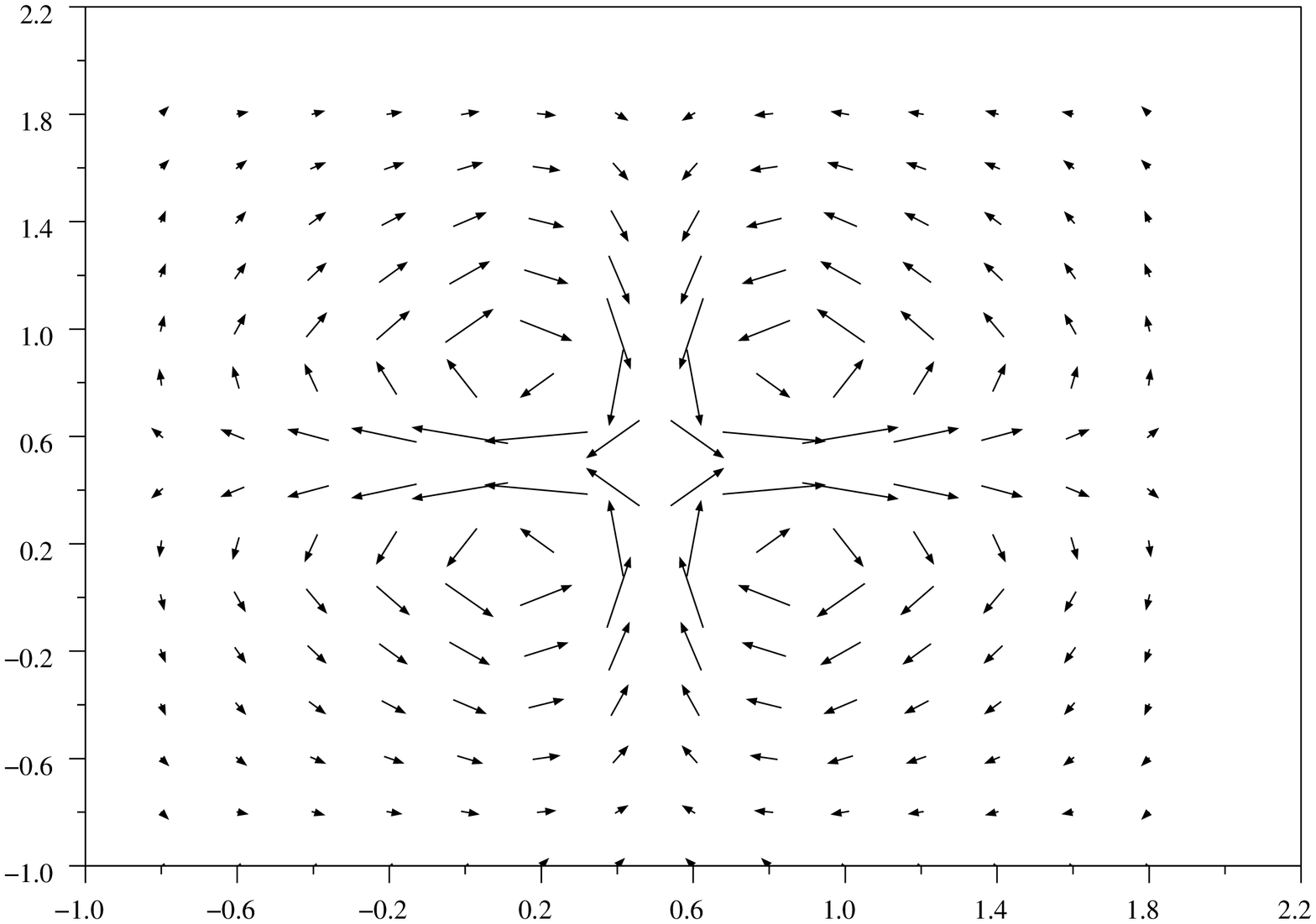}  
\end{tabular}  
\end{center}  
\caption{\label{fig:wave2D} Isotropic 2D generating divergence-free wavelets $\Psi_{\textrm{div}}^{(1,0)}$ (left),  
 $\Psi_{\textrm{div}}^{(0,1)}$ (center) and  $\Psi_{\textrm{div}}^{(1,1)}$ (right).}  
\end{figure}

It can be easily  seen that the dilated and translated functions $\Psi_{\textrm{div},j,\k}^{\varepsilon}=2^j\Psi_{\textrm{div}}^{\varepsilon}(2^jx_1-k_1,2^jx_2-k_2)$  
with   
$j,k_1,k_2\in \Z$ and $\varepsilon\in\{0,1\}^2\setminus(0,0)$ span the space  
$\vec{\textbf{H}}_{\textrm{div}}(\R^2)$ of divergence-free vector 
functions in $\R^2$. We represent in figure  
\ref{fig:wave2D} the three generating functions in the case of   
spline generators of degree 1 and 2 of  figure \ref{V0V1}.  
  
More generally, the construction of divergence-free wavelets in $\R^n$ 
is carried out  by suitable combinations of  
tensor products of functions $\phi_0$, $\psi_0$, and $\phi_1$, 
$\psi_1$, satisfying the above proposition  
(see \cite{L92, U96, U00}). These allow  
to state the following theorem of existence of {\it isotropic} divergence-free  
wavelet bases in the general case \cite{L92}:\\  
  
\noindent \textbf{Theorem}: {\it  
There exist $(n-1)(2^n-1)$ vector functions  
{\rm $\Psi_{\textrm{div},i}^{\varepsilon}\in \vec{\textbf{H}}_{\textrm{div}}(\R^n)$} ($\varepsilon\in\Omega_n^*$ of cardinal $(2^n-1)$, $1\leq i\leq n-1$) compactly  
supported, such that every vector function  
$\u\in \vec{\textbf{H}}_{\textrm{div}}(\R^n)$ can be expanded in a 
unique way:  
{\rm 
$$  
\u=\sum_{j\in \Z}\sum_{\varepsilon\in \Omega_n^*}\sum_{\k\in 
\Z^n}d_{\mbox{div},i,j,\k}^{\varepsilon}~ \Psi_{\textrm{div},i,j, \k}^{\varepsilon} 
$$ } 
and one has, for a constant $C>0$ independent from $\u=(u_1, u_2, \dots, u_n)$:  
{\rm $$  
\frac{1}{C}\|\u\|_{\vec L^2}\leq \bigg\{  
\sum_{j\in \Z}\sum_{\varepsilon\in \Omega_n}\sum_{\k\in \Z^n}|d_{\mbox{div},i,j,\k}^{\varepsilon}|^2\bigg\}^{1/2}  
\leq C\|\u\|_{\vec L^2}  
$$ } 
where $\|\u|_{\vec L^2}^2=\sum_{i=1}^n \| u_i \|_{L^2}^2$.  
} 
\\

These wavelets have already been studied by several authors for 2D analyses of turbulent flows  
\cite{AUDRL02,KKR00}, and also to solve the Stokes problem in two and three dimensions \cite{U94, U96, U00}.  
From now on, we will focus on the 2D and 3D case, and we will see 
that the expansion of compressible flow in terms of divergence-free 
wavelet bases is not uniquely given. 
We will also present in a practical way, the  associated fast algorithms.  
  
\section{Practical implementation of divergence-free wavelets}  
  
In this section, we present in detail the two and three-dimensional divergence-free wavelet bases, and we  
study how to compute the associated fast wavelet transforms. We present the constructions in the {\it  
isotropic} multivariate wavelet case (see section 2.3),  and also in 
the {\it anisotropic} one, which differs somewhat from previous studies.\\  
  
In the following, we are given two 1D multiresolution analyses $(V^0_j)$ and $(V^1_j)$ satisfying the  
proposition of section 2.4. We note $\phi_0$, $\psi_0$ and $\phi_1$, $\psi_1$ their  
associated (one-dimensional) scaling functions and wavelets.  
\subsection{Isotropic divergence-free wavelet transforms}  
  
\subsubsection{The 2D case} 
\label{isotropic2D}  
  
The starting point of the construction lies in considering as 2D multiresolution analysis of $L^2(\R^2)^2$  
the vector space of tensor-products  $(V^1_j\otimes V^0_j)\times 
(V^0_j\otimes V^1_j)$.  
In the isotropic case, the {\it 2D  
scaling functions} $\Phi_1$, $\Phi_2$ and the {\it 2D wavelets}  
$\Psi_1^{\varepsilon}$, $\Psi_2^{\varepsilon}$  
of this MRA are given by:  
$$  
\begin{array}{llcll}  
\Phi_1(x_1,x_2)= & \left| \begin{array}{l} \phi_1(x_1)\phi_0(x_2) \\ 0 \end{array} \right.& &  
\Phi_2(x_1,x_2)= & \left| \begin{array}{l} 0 \\ \phi_0(x_1)\phi_1(x_2) \end{array} \right.\\  
\vspace{0.3mm} & \vspace{0.3mm} & \vspace{0.3mm} & \vspace{0.3mm} \\  
\Psi_1^{(1,0)}(x_1,x_2)= & \left| \begin{array}{l} \psi_1(x_1)\phi_0(x_2) \\ 0 \end{array} \right.& &  
\Psi_2^{(1,0)}(x_1,x_2)= & \left| \begin{array}{l} 0 \\ \psi_0(x_1)\phi_1(x_2) \end{array} \right.\\  
\vspace{0.3mm} & \vspace{0.3mm} & \vspace{0.3mm} & \vspace{0.3mm} \\  
\Psi_1^{(0,1)}(x_1,x_2)= & \left| \begin{array}{l} \phi_1(x_1)\psi_0(x_2) \\ 0 \end{array} \right.& &  
\Psi_2^{(0,1)}(x_1,x_2)= & \left| \begin{array}{l} 0 \\ \phi_0(x_1)\psi_1(x_2) \end{array} \right.\\  
\vspace{0.3mm} & \vspace{0.3mm} & \vspace{0.3mm} & \vspace{0.3mm} \\  
\Psi_1^{(1,1)}(x_1,x_2)= & \left| \begin{array}{l} \psi_1(x_1)\psi_0(x_2) \\ 0 \end{array} \right.& &  
\Psi_2^{(1,1)}(x_1,x_2)= & \left| \begin{array}{l} 0 \\ \psi_0(x_1)\psi_1(x_2) \end{array} \right.  
\end{array}  
$$  
As explained in section 2.3, the functions  
\begin{equation*}
\bigg\{ \Psi_{i,j,\k}^{\varepsilon}(x_1, x_2)=2^j \Psi_{i}^{\varepsilon}(2^jx_1-k_1,2^jx_2-k_2)\bigg\}\\ 
\end{equation*} 
with $j\in\Z, \k=(k_1,k_2)\in \Z^2, \varepsilon\in\{(0,1),(1,0), (1,1)\}, i=1,2$, form a basis of $(L^2(\mathbb{R}^2))^2$. Then, a velocity field  
$\u$ in $(L^2(\mathbb{R}^2))^2$ has the following wavelet decomposition:  
\begin{eqnarray}  
\label{V}  
\u &= ~\sum_{j\in \Z}\sum_{\k\in \Z^2} &\left(d_{1,j,\k}^{(1,0)}~\Psi_{1,j,\k}^{(1,0)}+  
d_{2,j,\k}^{(1,0)}~\Psi_{2,j,\k}^{(1,0)}+d_{1,j,\k}^{(0,1)}~\Psi_{1,j,\k}^{(0,1)}\right. 
\\ 
\nonumber 
& & + \left. d_{2,j,\k}^{(0,1)}~\Psi_{2,j,\k}^{(0,1)}+d_{1,j,\k}^{(1,1)}~\Psi_{1,j,\k}^{(1,1)}+  
d_{2,j,\k}^{(1,1)}~\Psi_{2,j,\k}^{(1,1)}\right)  
\end{eqnarray}  
  
By a simple linear change of basis we are able to find out a divergence-free wavelet basis, and its complement:  
$$  
\begin{array}{ccl}  
\left\{ \begin{array}{c} \Psi_1^{(1,0)} \\ \Psi_2^{(1,0)} \end{array} \right.& \longrightarrow &  
\left\{ \begin{array}{l} \Psi_{\textrm{div}}^{(1,0)} =\Psi_2^{(1,0)} -\frac{1}{4}[\Psi_1^{(1,0)} -\Psi_1^{(1,0)}(x_1,x_2-1)] \\  
\Psi_{n}^{(1,0)} =\Psi_{1}^{(1,0)}  \end{array} \right. \\  
\vspace{0.3mm} & \vspace{0.3mm} & \vspace{0.3mm} \\  
\left\{ \begin{array}{c} \Psi_1^{(0,1)}  \\ \Psi_2^{(0,1)}  \end{array} \right.& \longrightarrow &  
\left\{ \begin{array}{l} \Psi_{\textrm{div}}^{(0,1)} =\Psi_1^{(0,1)} -\frac{1}{4}[\Psi_2^{(0,1)} -\Psi_2^{(0,1)}(x_1-1,x_2)] \\  
\Psi_{n}^{(0,1)} =\Psi_{2}^{(0,1)}  \end{array} \right. \\  
\vspace{0.3mm} & \vspace{0.3mm} & \vspace{0.3mm} \\  
\left\{ \begin{array}{c} \Psi_1^{(1,1)}  \\ \Psi_2^{(1,1)}  \end{array} \right.& \longrightarrow &  
\left\{ \begin{array}{l} \Psi_{\textrm{div}}^{(1,1)} =\Psi_1^{(1,1)} -\Psi_2^{(1,1)}  \\  
\Psi_{n}^{(1,1)} =\Psi_1^{(1,1)} +\Psi_2^{(1,1)}  \end{array} \right.  
\end{array}  
$$

The first functions $\Psi_{\textrm{div}}^{\varepsilon}$ yield a divergence-free  
basis (already shown in example of section 2.4),  and the second ones $\Psi_{n}^{\varepsilon}$  are the complement functions  
corresponding to non  
divergence-free part of the data. Remark that the functions $\Psi_{\textrm{div}}^{\varepsilon}$ and  
$\Psi_{n}^{\varepsilon}$ are not orthogonal. Moreover, the choice of the functions $\Psi_{n}^{\varepsilon}$  
is not unique, and this choice has an influence on the values of all the coefficients, when applying this  
transform to a compressible flow. The choice of 
$\Psi_{n}^{\varepsilon}$ we made in this work, leads to very simple formula 
to obtain the divergence-free coefficients.\\  
Now, the expansion (\ref{V}) of a given vector function $\u$ can be rewritten:  
\begin{eqnarray}  
\nonumber  
\u &=& \sum_{j\in \Z}\sum_{\k\in \Z^2}\left(d_{\textrm{div},j,\k}^{(1,0)}~\Psi_{\textrm{div},j,\k}^{(1,0)}  
+d_{\textrm{div},j,\k}^{(0,1)}~\Psi_{\textrm{div},j,\k}^{(0,1)}+  
d_{\textrm{div},j,\k}^{(1,1)}~\Psi_{\textrm{div},j,\k}^{(1,1)}\right)  
\\  
&&+ \sum_{j\in \Z}\sum_{\k\in \Z^2}\left(d_{n,j,\k}^{(1,0)}~\Psi_{n,j,\k}^{(1,0)} +d_{n,j,\k}^{(0,1)}~\Psi_{n,j,\k}^{(0,1)}+d_{n,j,\k}^{(1,1)}~\Psi_{n,j,\k}^{(1,1)}\right)  
\end{eqnarray}  
where the new  coefficients are directly expressed from the original ones by:

\begin{eqnarray}  
\label{ddiv} \hspace{-1cm} 
(d_{div})~ \left\{  
\begin{array}{l}  
d_{\textrm{div},j,\k}^{(1,0)}=d_{2,j,\k}^{(1,0)} \\  
\vspace{-0.6mm} \\  
d_{\textrm{div},j,\k}^{(0,1)}=d_{1,j,\k}^{(0,1)} \\  
\vspace{-0.6mm} \\  
d_{\textrm{div},j,\k}^{(1,1)}=\frac{1}{2}d_{1,j,\k}^{(1,1)}-\frac{1}{2}d_{2,j,\k}^{(1,1)}  
\end{array} \right.  
& \hspace{0cm}  
(d_n) ~\left\{  
\begin{array}{l}  
d_{n,j,\k}^{(1,0)}=d_{1,j,\k}^{(1,0)}+\frac{1}{4}d_{2,j,\k}^{(1,0)}-\frac{1}{4}d_{2,j,k_1,k_2-1}^{(1,0)} \\  
\vspace{-0.6mm} \\  
d_{n,j,\k}^{(0,1)}=d_{2,j,\k}^{(0,1)}+\frac{1}{4}d_{1,j,\k}^{(0,1)}-\frac{1}{4}d_{1,j,k_1-1,k_2}^{(0,1)} \\  
\vspace{-0.6mm} \\  
d_{n,j,\k}^{(1,1)}=\frac{1}{2}d_{1,j,\k}^{(1,1)}+\frac{1}{2}d_{2,j,\k}^{(1,1)}  
\end{array}\right. &   
\end{eqnarray}  
As one can see, the computation of divergence-free wavelet coefficients $d_{\textrm{div},j,\k}^{\varepsilon}$  
is reduced to a very simple linear combination of the standard wavelet coefficients  
$d_{i,j,\k}^{\varepsilon}$ provided by the biorthogonal fast wavelet transform, which makes programming easy.  
  
\

\noindent \textbf{Remark}: $\textrm{div}(\Psi_{n}^{(1,0)})$, $\textrm{div}(\Psi_{n}^{(0,1)})$ and  
$\textrm{div}(\Psi_{n}^{(1,1)})$ are generating functions of the 
scalar space $V^0\otimes V^0$.  
Moreover, for all $\u$, we have  
$$  
\textrm{div} ~\u = \sum_{j\in \Z}\sum_{\k\in \Z^2} \left(  
d_{n,j,\k}^{(1,0)}~\textrm{div}(\Psi_{n,j,\k}^{(1,0)})+  
d_{n,j,\k}^{(0,1)}~\textrm{div}(\Psi_{n,j,\k}^{(0,1)})+  
d_{n,j,\k}^{(1,1)}~\textrm{div}(\Psi_{n,j,\k}^{(1,1)}) \right)  
$$  
Then, the incompressibility condition $\textrm{div} ~\u= 0$ is equivalent to  
$d_{n,j,\k}^{\varepsilon}=0$, for all $j,\k, \varepsilon$.\\  
  
For incompressible flows, since the biorthogonal projectors  onto the spaces $(V^1_j\otimes V^0_j)\times 
(V^0_j\otimes V^1_j)$ commute with partial derivatives \cite{L92},   
the coefficients $d_{\textrm{div},j,\k}^{\varepsilon}$ are uniquely determined, by the formula $(d_{div})$ in  
equation (\ref{ddiv}). We present in section 5.2, numerical experiments on 2D incompressible turbulent  
flows.\\  
Difficulties arise when we want to compute the divergence-free part of a compressible flow. Because of the  
non-orthogonality between the divergence-free basis  
$(\Psi_{\textrm{div}}^{\varepsilon})$ and its complement  $(\Psi_{n}^{\varepsilon})$, the values of the  
divergence-free wavelet coefficients depend on the choice of the complement basis. We address  
 this problem in the specific section 4, on the Hodge decomposition.   
  
\subsubsection{The isotropic 3D case}
\label{iso3D} 
  
The construction of the 3D divergence-free wavelet bases may be obtained in a similar fashion as for the  
2D bases. Again, the key-point is to start with a vector multiresolution analysis of $(L^2(\R^3))^3$, of  
the type 
$$(V^1_j\otimes V^0_j\otimes V^0_j)\times (V^0_j\otimes V^1_j\otimes 
V^0_j)\times (V^0_j\otimes V^0_j\otimes V^1_j)$$ 
This MRA provides naturally with 3 generating {\it 3D-vector scaling functions}:  
  
$$  
\begin{array}{lclcl}  
\Phi_{1}(x_1,x_2,x_3) =  
\left| \begin{array}{l}  
\phi_1(x_1)\phi_0(x_2)\phi_0(x_3) \\ 0 \\ 0 \end{array} \right. & &  
\Phi_{2} =\left| \begin{array}{l}  
0 \\ \phi_0\phi_1\phi_0 \\ 0 \end{array} \right. & &  
\Phi_{3} =\left| \begin{array}{l}  
0 \\ 0 \\ \phi_0\phi_0\phi_1 \end{array} \right.  
\end{array}  
$$  
and 21 generating {\it 3D-vector wavelets}:  
 $$
 \left\{\Psi_{i}^{\varepsilon}~|~ i=1,2,3~,~\varepsilon=(\varepsilon_1, \varepsilon_2,\varepsilon_3)~\mbox{with}~\varepsilon_i=0,1~\mbox{and}~\varepsilon\neq(0,0,0)\right\}
 $$
For example, we give below the expressions of wavelets $\Psi_{i}^{(1,0,0)}$, $\Psi_{i}^{(1,0,0)}$, and $\Psi_{i}^{(1,0,0)}$:  
$$  
\begin{array}{lclcl}  
\Psi_{1}^{(1,0,0)}(x_1,x_2,x_3) =\left| \begin{array}{l}  
\psi_1(x_1)\phi_0(x_2)\phi_0(x_3) \\ 0 \\ 0 \end{array} \right. & &  
\Psi_{2}^{(1,0,0)} =\left| \begin{array}{l}  
0 \\ \psi_0\phi_1\phi_0 \\ 0 \end{array} \right. & &  
\Psi_{3}^{(1,0,0)}=\left| \begin{array}{l}  
0 \\ 0 \\ \psi_0\phi_0\phi_1 \end{array} \right.\\  
\vspace{1mm} & \vspace{1mm} & \vspace{1mm} & \vspace{1mm} & \vspace{1mm} \\  
\Psi_{1}^{(1,1,0)}(x_1,x_2,x_3) =\left| \begin{array}{l}  
\psi_1(x_1)\psi_0(x_2)\phi_0(x_3) \\ 0 \\ 0 \end{array} \right. & &  
\Psi_{2}^{(1,1,0)} =\left| \begin{array}{l}  
0 \\ \psi_0\psi_1\phi_0 \\ 0 \end{array} \right. & &  
\Psi_{3}^{(1,1,0)} =\left| \begin{array}{l}  
0 \\ 0 \\ \psi_0\psi_0\phi_1 \end{array} \right.\\  
\vspace{1mm} & \vspace{1mm} & \vspace{1mm} & \vspace{1mm} & \vspace{1mm} \\  
\Psi_{1}^{(1,1,1)}(x_1,x_2,x_3) =\left| \begin{array}{l}  
\psi_1(x_1)\psi_0(x_2)\psi_0(x_3) \\ 0 \\ 0 \end{array} \right. & &  
\Psi_{2}^{(1,1,1)} =\left| \begin{array}{l}  
0 \\ \psi_0\psi_1\psi_0 \\ 0 \end{array} \right. & &  
\Psi_{3}^{(1,1,1)} =\left| \begin{array}{l}  
0 \\ 0 \\ \psi_0\psi_0\psi_1 \end{array} \right.\\  
\end{array}  
$$
and it goes similarly for $(0,1,0)$, $(0,0,1)$, $(1,0,1)$ and $(0,1,1)$.\\  
Let introduce $\Omega_3^*=\{\varepsilon\in \{0,1\}^3\setminus(0,0,0)\}$. The {\it isotropic} wavelet expansion of a
given function $\u$ writes:  
\begin{equation}  
\label{V3d}  
\u= \sum_{j\in \Z}\sum_{\k\in \Z^3}\sum_{\varepsilon\in \Omega_3^*} \left(  
d_{1,j,\k}^{\varepsilon}~\Psi_{1,j,\k}^{\varepsilon}+d_{2,j,\k}^{\varepsilon}~\Psi_{2,j,\k}^{\varepsilon}  
+d_{3,j,\k}^{\varepsilon}~\Psi_{3,j,\k}^{\varepsilon} \right)  
\end{equation}  
 Following theorem of section \ref{theory}, there exist 14 kinds of isotropic divergence-free wavelets, with arbitrary possible 
choices concerning the privileged direction of each basis function. In the following  
we do not detail all the expressions we choose, but only some typical ones:  
 
$$  
\begin{array}{l}  
\Psi_{\textrm{div},1}^{(1,0,0)}(x_1,x_2,x_3) =\left| \begin{array}{l}  
-\frac{1}{4}\psi_1(x_1)[\phi_0(x_2)-\phi_0(x_2-1)]\phi_0(x_3) \\  
\psi_0(x_1)\phi_1(x_2)\phi_0(x_3) \\ 0 \end{array} \right. \\  
\vspace{1mm} \\  
\Psi_{\textrm{div},2}^{(1,0,0)}(x_1,x_2,x_3) =\left| \begin{array}{l}  
-\frac{1}{4}\psi_1(x_1)\phi_0(x_2)[\phi_0(x_3)-\phi_0(x_3-1)] \\  
0 \\ \psi_0(x_1)\phi_0(x_2)\phi_1(x_3)\end{array} \right. \\  
\vspace{1mm} \\  
\Psi_{\textrm{div},1}^{(1,1,0)}(x_1,x_2,x_3) =\left| \begin{array}{l}  
\psi_1(x_1)\psi_0(x_2)\phi_0(x_3) \\  
-\psi_0(x_1)\psi_1(x_2)\phi_0(x_3) \\ 0 \end{array} \right. \\  
\vspace{1mm} \\  
\Psi_{\textrm{div},2}^{(1,1,0)}(x_1,x_2,x_3) =\left| \begin{array}{l}  
-\frac{1}{8}\psi_1(x_1)\psi_0(x_2)[\phi_0(x_3)-\phi_0(x_3-1)] \\  
-\frac{1}{8}\psi_0(x_1)\psi_1(x_2)[\phi_0(x_3)-\phi_0(x_3-1)] \\  
\psi_0(x_1)\psi_0(x_2)\phi_1(x_3) \end{array} \right. \\  
\vspace{1mm} \\  
\Psi_{\textrm{div},1}^{(1,1,1)}(x_1,x_2,x_3) =\left| \begin{array}{l}  
-\psi_1(x_1)\psi_0(x_2)\psi_0(x_3) \\ 0 \\  
\psi_0(x_1)\psi_0(x_2)\psi_1(x_3) \end{array}\right.  
\vspace{1mm} \\  
\Psi_{\textrm{div},2}^{(1,1,1)}(x_1,x_2,x_3) =\left| \begin{array}{l}  
0 \\ \psi_0(x_1)\psi_1(x_2)\psi_0(x_3) \\  
-\psi_0(x_1)\psi_0(x_2)\psi_1(x_3) \end{array}\right.  
\end{array}  
$$  
It goes similarly for all basis function: for each $\varepsilon\in\Omega_3^*$ given, two divergence-free wavelets $\Psi_{\textrm{div},i}^{\varepsilon}$, $i=1,2$ are carried out by linear combination of  
$\Psi_{1}^{\varepsilon}$, $\Psi_{2}^{\varepsilon}$, $\Psi_{3}^{\varepsilon}$, in order to satisfy the divergence-free condition.  
The complement wavelet $\Psi_{\textrm{n}}^{\varepsilon}$ is constructed in order to take care of the symmetry. For example, we consider: 
$$\left\{  
\begin{array}{l}  
\Psi_{\textrm{div},1}^{(1,0,0)}=\Psi_{2}^{(1,0,0)}-\frac{1}{4}(\Psi_{1}^{(1,0,0)}(.,.,.)-\Psi_{1}^{(1,0,0)}(.,.-1,.)) \\  
\vspace{1mm} \\  
\Psi_{\textrm{div},2}^{(1,0,0)}=\Psi_{3}^{(1,0,0)}-\frac{1}{4}(\Psi_{1}^{(1,0,0)}(.,.,.)-\Psi_{1}^{(1,0,0)}(.,.,.-1)) \\  
\vspace{1mm} \\  
\Psi_{\textrm{n}}^{(1,0,0)}=\Psi_{1}^{(1,0,0)}   
\end{array}\right.  
$$  
and similarly for $\Psi_{\textrm{div},i}^{(0,1,0)}$ and $\Psi_{\textrm{div},i}^{(0,0,1)}$, $i=1,2$.  
$$\left\{  
\begin{array}{ll}  
\Psi_{\textrm{div},1}^{(1,1,0)}=&\Psi_{1}^{(1,1,0)}-\Psi_{2}^{(1,1,0)} \\  
\vspace{1mm} & \vspace{1mm} \\  
\Psi_{\textrm{div},2}^{(1,1,0)}=&\Psi_{3}^{(1,1,0)}-\frac{1}{8}(\Psi_{1}^{(1,1,0)}(.,.,.)-\Psi_{1}^{(1,1,0)}(.,.,.-1))
\\ 
 &-\frac{1}{8}(\Psi_{2}^{(1,1,0)}(.,.,.)-\Psi_{2}^{(1,1,0)}(.,.,.-1)) \\  
\vspace{1mm} & \vspace{1mm} \\  
\Psi_{\textrm{n}}^{(1,1,0)}=&\Psi_{1}^{(1,1,0)}+\Psi_{2}^{(1,1,0)}   
\end{array}\right.  
$$  
and similarly for $\Psi_{\textrm{div},i}^{(0,1,1)}$ and $\Psi_{\textrm{div},i}^{(1,0,1)}$, $i=1,2$.  
$$\left\{  
\begin{array}{l}  
\Psi_{\textrm{div},1}^{(1,1,1)}=\Psi_{3}^{(1,1,1)}-\Psi_{1}^{(1,1,1)} \\  
\vspace{1mm} \\  
\Psi_{\textrm{div},2}^{(1,1,1)}=\Psi_{2}^{(1,1,1)}-\Psi_{3}^{(1,1,1)} \\  
\vspace{1mm} \\  
\Psi_{\textrm{n}}^{(1,1,1)}=\Psi_{1}^{(1,1,1)}+\Psi_{2}^{(1,1,1)}+\Psi_{3}^{(1,1,1)}  
\end{array}\right.  
$$  
Now we can rewrite (\ref{V3d}):  
$$  
\u= \sum_{j\in \Z}\sum_{\k\in \Z^3}\sum_{\varepsilon\in \Omega_3^*} \left( 
d_{\textrm{div},1,j,\k}^{\varepsilon}\Psi_{\textrm{div},1,j,\k}^{\varepsilon}  
+d_{\textrm{div},2,j,\k}^{\varepsilon}\Psi_{\textrm{div},2,j,\k}^{\varepsilon}  
+d_{\textrm{n},j,\k}^{\varepsilon}\Psi_{\textrm{n},j,\k}^{\varepsilon} \right) 
$$  
where the divergence-free wavelets are simply obtained from the standard ones, for example:  
$$ \begin{array}{ll} 
\left\{\begin{array}{l}  
d_{\textrm{div},1}^{(1,0,0)}=d_{2}^{(1,0,0)} \\  
\vspace{0mm} \\  
d_{\textrm{div},2}^{(1,0,0)}=d_{3}^{(1,0,0)}\end{array}\right. 
& 
\left\{\begin{array}{l} 
d_{\textrm{div},1}^{(1,1,0)}=\frac{1}{2}(d_{1}^{(1,1,0)}-d_{2}^{(1,1,0)}) \\  
\vspace{0mm} \\  
d_{\textrm{div},2}^{(1,1,0)}=d_{3}^{(1,1,0)} \end{array}\right.  
\\
\left\{\begin{array}{l} 
d_{\textrm{div},1}^{(1,1,1)}=\frac{1}{3}(-2d_{1}^{(1,1,1)}+d_{2}^{(1,1,1)}+d_{3}^{(1,1,1)}) \\  
\vspace{0mm} \\  
d_{\textrm{div},2}^{(1,1,1)}=\frac{1}{3}(-d_{1}^{(1,1,1)}+2d_{2}^{(1,1,1)}-d_{3}^{(1,1,1)})  
\end{array}\right. & 
\end{array} 
$$ 
The complement coefficients are in this case: 
$$(d_{\mbox{n}}) ~\left\{ \begin{array}{l} 
d_{\textrm{n},k}^{(1,0,0)}=d_{1,k_1,k_2,k_3}^{(1,0,0)}+\frac{1}{4}(d_{2,k_1,k_2,k_3}^{(1,0,0)}-d_{2,k_1,k_2-1,k_3}^{(1,0,0)})  
+\frac{1}{4}(d_{3,k_1,k_2,k_3}^{(1,0,0)}-d_{3,k_1,k_2,k_3-1}^{(1,0,0)}) \\  
\vspace{0mm} \\  
d_{\textrm{n},k}^{(1,1,0)}=\frac{1}{2}(d_{1}^{(1,1,0)}+d_{2}^{(1,1,0)})  
+\frac{1}{8}(d_{3,k_1,k_2,k_3}^{(1,1,0)}-d_{3,k_1,k_2,k_3-1}^{(1,1,0)}) \\  
\vspace{0mm} \\  
d_{\textrm{n}}^{(1,1,1)}=\frac{1}{3}(d_{1}^{(1,1,1)}+d_{2}^{(1,1,1)}+d_{3}^{(1,1,1)}) \end{array} \right. 
$$  
As for the two-dimensional case, the computation of divergence-free  wavelet coefficients of any 3D vector field lies in a short linear combination of standard biorthogonal wavelet coefficients, arising from the fast wavelet transform. 
  
\subsection{Anisotropic divergence-free wavelet transforms}  
  
In this section we will construct anisotropic wavelets that are divergence-free. Since  
the one-dimensional wavelets verify ${\psi_1}' = 4 \psi_0$, we derive easily divergence-free
wavelet bases by tensor products of one-dimensional wavelets. We detail in the following the
construction of such bases in the two and three dimensional cases. 
  
\subsubsection{The anisotropic 2D case} 
 \label{div2Dan} 
  
The 2D anisotropic divergence-free wavelets are given by: 
$$  
\Psi_{\textrm{div},\j,\k}^{\textrm{an}}(x_1,x_2)=\left| \begin{array}{l}  
2^{j_2}\psi_1(2^{j_1}x_1-k_1)\psi_0(2^{j_2}x_2-k_2) \\  
-2^{j_1}\psi_0(2^{j_1}x_1-k_1)\psi_1(2^{j_2}x_2-k_2)\end{array}\right.  
$$  
where $\j=(j_1,j_2)\in\Z^2$ is the scale parameter (which is different
in both directions), and  
$\k=(k_1,k_2)\in\Z^2$ is the position parameter. When the indices $\k$ and $\j$ vary
in $\Z^2$, the family
$\{\Psi_{\textrm{div},\j,\k}^{\textrm{an}}\}$ forms a basis of  
$\vec{\textbf{H}}_{\textrm{div}}(\R^2)$. \\ 
We introduced the complement functions:  
$$  
\Psi_{\textrm{n},\j,\k}^{\textrm{an}}(x_1,x_2)=\left| \begin{array}{l}  
2^{j_1}\psi_1(2^{j_1}x_1-k_1)\psi_0(2^{j_2}x_2-k_2) \\  
2^{j_2}\psi_0(2^{j_1}x_1-k_1)\psi_1(2^{j_2}x_2-k_2)\end{array}\right.  
$$  
  
The anisotropic divergence-free wavelet transform works similarly as the isotropic one but with fewer
elements to be computed. The decomposition of a given vector function $\u$ begins with the anisotropic wavelet decomposition
associated to the MRA  $(V^1_j\otimes V^0_j)\times (V^0_j\otimes V^1_j)$ (see section \ref{multi}); 
$$  
\u= \sum_{\j\in \Z^2}\sum_{\k \in \Z^2}\left( d_{1,\j,\k}^{\textrm{an}}~\Psi_{1,\j,\k}^{\textrm{an}}+  
d_{2,\j,\k}^{\textrm{an}}~\Psi_{2,\j,\k}^{\textrm{an}} \right) 
$$  
with:  
$$  
\begin{array}{ll}  
\Psi_{1,\j,\k}^{\textrm{an}}(x_1,x_2)= & \left| \begin{array}{l} \psi_1(2^{j_1}x_1-k_1)\psi_0(2^{j_2}x_2-k_2) \\ 0  
\end{array} \right.\\  
\vspace{1mm} & \vspace{1mm} \\  
\Psi_{2,\j,\k}^{\textrm{an}}(x_1,x_2)= & \left| \begin{array}{l} 0 \\ \psi_0(2^{j_1}x_1-k_1)\psi_1(2^{j_2}x_2-k_2) \end{array} \right.  
\end{array}  
$$  
for $\j,\k\in \Z^2$. Remark that for more simplicity, the dilated
functions are not 
normalized, in $L^2$-norm.\\ 
$\u$ can be expanded onto the new basis: 
\begin{equation} 
\label{vdiv} 
\u= \sum_{\j\in \Z^2}\sum_{\k\in \Z^2}\left( d_{\textrm{div},\j,\k}^{\textrm{an}}~\Psi_{\textrm{div},\j,\k}^{\textrm{an}}+  
d_{n,\j,\k}^{\textrm{an}}~\Psi_{n,\j,\k}^{\textrm{an}} \right) 
\end{equation} 
with the corresponding coefficients:  
\begin{equation} 
\label{ddiv}  
\begin{array}{l}  
d_{\textrm{div},\j,\k}^{\textrm{an}}=
\frac{2^{j_2}}{2^{2j_1}+2^{2j_2}}~d_{1,\j,\k}^{\textrm{an}}-\frac{2^{j_1}}{2^{2j_1}+2^{2j_2}}~d_{2,\j,\k}^{\textrm{an}}  \\  
\vspace{-2mm} \\  
d_{\mbox{n},\j,\k}^{\textrm{an}}=
\frac{2^{j_1}}{2^{2j_1}+2^{2j_2}}~d_{1,\j,\k}^{\textrm{an}}+\frac{2^{j_2}}{2^{2j_1}+2^{2j_2}}~d_{2,\j,\k}^{\textrm{an}}
\end{array}  
\end{equation} 
  
\subsubsection{The anisotropic 3D case} 
  
In the same way, the (non normalized) anisotropic 3D divergence-free wavelets take the form:  
$$  
\begin{array}{c}  
\Psi_{\textrm{div},1,\j,\k}^{\textrm{an}}(x_1,x_2,x_3)=\left| \begin{array}{l}  
2^{j_2}\psi_1(2^{j_1}x_1-k_1)\psi_0(2^{j_2}x_2-k_2)\psi_0(2^{j_3}x_3-k_3) \\   
-2^{j_1}\psi_0(2^{j_1}x_1-k_1)\psi_1(2^{j_2}x_2-k_2)\psi_0(2^{j_3}x_3-k_3) \\ 0 \end{array}\right.  \\  
\vspace{0cm} \\  
\Psi_{\textrm{div},2,\j,\k}^{\textrm{an}}(x_1,x_2,x_3)=\left| \begin{array}{l} 0 \\  
2^{j_3}\psi_0(2^{j_1}x_1-k_1)\psi_1(2^{j_2}x_2-k_2)\psi_0(2^{j_3}x_3-k_3) \\  
-2^{j_2}\psi_0(2^{j_1}x_1-k_1)\psi_0(2^{j_2}x_2-k_2)\psi_1(2^{j_3}x_3-k_3) \end{array}\right. \\ \vspace{0cm} 
\\ 
\Psi_{\textrm{div},3,\j,\k}^{\textrm{an}}(x_1,x_2,x_3)=\left| \begin{array}{l}  
-2^{j_3}\psi_1(2^{j_1}x_1-k_1)\psi_0(2^{j_2}x_2-k_2)\psi_0(2^{j_3}x_3-k_3) \\ 0 \\  
2^{j_1}\psi_0(2^{j_1}x_1-k_1)\psi_0(2^{j_2}x_2-k_2)\psi_1(2^{j_3}x_3-k_3)  
\end{array}\right. 
\end{array}  
$$  
with $\j=(j_1,j_2,j_3),~\k=(k_1,k_2,k_3)\in \Z^3$.\\  
The 3D divergence-free basis is carried out by considering only two types of functions among the three above.  
As complement basis we introduce a function which is the most as possible orthogonal to the previous
ones: \\ 
$$  
\begin{array}{c}  
\Psi_{\textrm{n},\j,\k}^{\textrm{an}}(x_1,x_2,x_3)=\left| \begin{array}{l}  
2^{j_1}\psi_1(2^{j_1}x_1-k_1)\psi_0(2^{j_2}x_2-k_2)\psi_0(2^{j_3}x_3-k_3) \\   
2^{j_2}\psi_0(2^{j_1}x_1-k_1)\psi_1(2^{j_2}x_2-k_2)\psi_0(2^{j_3}x_3-k_3) \\ 
2^{j_3}\psi_0(2^{j_1}x_1-k_1)\psi_0(2^{j_2}x_2-k_2)\psi_1(2^{j_3}x_3-k_3) \end{array}\right.  
\end{array}  
$$  
The operations to compute divergence-free coefficients and complement coefficients are similar to the 2D case.  
  
  
\section{An iterative algorithm to compute the Hodge wavelet decomposition} 
\label{hodge}  
  
\subsection{Principle of the Hodge decomposition}  
  
The Hodge decomposition consists in splitting a vector function $\u\in (L^2(\R^n))^n$ into  
its divergence-free component $\u_{\mbox{div}}$ and a  gradient vector.  More precisely, there exist a pressure
$p$ and a stream-function $\psi$ such that: 
 
\begin{equation}  
\label{hodge} 
\u=\u_{\mbox{div}} + \vec \nabla p~~~\mbox{and}~~~ \u_{\mbox{div}}=\vec{\mbox{curl}}~ \psi 
\end{equation} 
Moreover, the functions $\vec{\mbox{curl}} ~\psi$ and $\vec \nabla p$ are orthogonal in $(L^2(\R^n))^n$. The
stream-function $\psi$ and the pressure $p$ are unique, up to an additive constant.\\ 
In $\R^2$, the stream-function is a scalar valued function, whereas in $\R^3$ it is a 3D vector function.\\ 
This decomposition may be viewed as the following orthogonal space splitting: 
$$(L^2(\R^n))^n=\vec{\textbf{H}}_{\textrm{div}}(\R^n)\oplus \vec{\textbf{H}}_{\textrm{curl}}(\R^n)$$ 
where we note 
$$\vec{\textbf{H}}_{\textrm{div}}(\R^n)=\{\v\in (L^2(\mathbb{R}^n))^n / \textrm{div}~\v\in L^2(\mathbb{R}^n)  
,\quad \textrm{div}~\v=0\}$$  
the space of divergence-free vector functions, and 
$$\vec{\textbf{H}}_{\textrm{curl}}(\R^n)=\{\v\in (L^2(\mathbb{R}^n))^n / \vec{\textrm{curl}}~\v\in
(L^2(\mathbb{R}^n))^n  
,\quad \vec{\textrm{curl}}~\v=\vec 0 \}$$  
the space of curl-free vector functions (if $n=2$ we have to replace $\vec{\textrm{curl}}~\v\in
(L^2(\mathbb{R}^n))^n$ by $\textrm{curl}~\v\in L^2(\mathbb{R}^2)$ in the definition). 
For the whole space $\R^n$, the proofs of the above decompositions can be derived easily, by mean of the Fourier
transform. In more general domains, we refer to \cite{GR79}.\\ 
 
The objective now is to provide with a wavelet Hodge decomposition. Since in the previous sections we have
constructed wavelet bases of $\vec{\textbf{H}}_{\textrm{div}}(\R^n)$, we have to work analogously to carry out
wavelet bases of $\vec{\textbf{H}}_{\textrm{curl}}(\R^n)$. 
  
\subsection{Construction of a  gradient wavelet basis}  
\label{gradient} 
  
A definition of wavelet bases for the space $\vec{\textbf{H}}_{\textrm{curl}}(\R^n)$ $(n=2, 3)$ has already been
provided by K. Urban in the isotropic case \cite{U00}. We will focus here on the construction of
{\it anisotropic}
curl-free vector wavelets in the 2D case (it goes similarly in the $n$-dimensional case).   
  
This construction is very similar to the divergence-free wavelet construction, despite some crucial differences.
The starting point here is to search wavelets in the MRA 
 $(V^0_J\otimes V^1_J)\times (V^1_J\otimes V^0_J)$ instead of  
$(V^1_J\otimes V^0_J)\times (V^0_J\otimes V^1_J)$, where the one-dimensional spaces $V_0$ and $V_1$ are related
by differentiation and integration (proposition of section \ref{theory}). \\ 
Since $\vec{\textbf{H}}_{\textrm{curl}}(\R^2)$ is the space of gradient functions in $L^2(\R^2)$, we construct
gradient wavelets by taking the gradient  
of a 2D wavelet basis of the MRA $(V^1_{j}\otimes V^1_{j})$.
If we avoid the $L^2$-normalization, the anisotropic gradient wavelets are defined by:  
$$  
\Psi_{\textrm{curl},\j,\k}^{\textrm{an}}(x_1,x_2)=\frac{1}{4}~ \vec\nabla \left(\psi_1(2^{j_1}x_1-k_1)\psi_1(2^{j_2}x_2-k_2)\right) 
=\left| \begin{array}{l}  
2^{j_1}\psi_0(2^{j_1}x_1-k_1)\psi_1(2^{j_2}x_2-k_2) \\ \vspace{0cm} 
\\ 
2^{j_2}\psi_1(2^{j_1}x_1-k_1)\psi_0(2^{j_2}x_2-k_2)\end{array}\right.  
$$  
Thus, when $\j=(j_1,j_2)$, $\k=(k_1,k_2)$ vary in $\Z^2$, the family
$\{\Psi_{\textrm{curl},\j,\k}^{\textrm{an}}\}$ forms a wavelet basis of $\vec{\textbf{H}}_{\textrm{curl}}(\R^2)$.

The decomposition algorithm on curl-free wavelets $\{\Psi_{\textrm{curl},\j,\k}^{\textrm{an}}\}$ works similarly
as the decomposition algorithm on anisotropic divergence-free wavelets. A vector function $\v$ is firstly
approximated in a space $(V^0_J\otimes V^1_J)\times (V^1_J\otimes V^0_J)$ by: 
$$  
\v^{\#}_J=\left| \begin{array}{l}  
v^{\#}_1=\sum_{\k\in \Z^2}c_{1,J,\k}~\phi_0(2^Jx_1-k_1)\phi_1(2^Jx_2-k_2) \\  
v^{\#}_2=\sum_{\k\in \Z^2}c_{2,J,\k}~\phi_1(2^Jx_1-k_1)\phi_0(2^Jx_2-k_2)  
\end{array} \right.  
$$  
By applying to $v^{\#}_1$ the standard anisotropic wavelet transform of  
$(V^0_J\otimes V^1_J)$  and to $v^{\#}_2$ the one of $(V^1_J\otimes V^0_J)$, it rewrites:  
$$  
\v_J^{\#}(x_1,x_2)=\left| \begin{array}{l}  
v^{\#}_1(x_1,x_2)=\sum_{j_1,j_2<J}\sum_{\k\in  
\Z^2}d^{\textrm{an}}_{1,j,\k}~\psi_0(2^{j_1}x_1-k_1)~\psi_1(2^{j_2}x_2-k_2) \\  
v^{\#}_2(x_1,x_2)=\sum_{j_1,j_2<J}\sum_{\k\in  
\Z^2}d^{\textrm{an}}_{2,j,\k}~\psi_1(2^{j_1}x_1-k_1)~\psi_0(2^{j_2}x_2-k_2)  
\end{array} \right.  
$$  
Let us introduce the vector wavelets: 
$$  
\begin{array}{llll}  
\Psi_{1,\j,\k}^{\textrm{an,\#}}(x_1,x_2)= & \left| \begin{array}{l}  
\psi_0(2^{j_1}x_1-k_1)\psi_1(2^{j_2}x_2-k_2) \\ 0 \end{array} \right.&  
\Psi_{2,\j,\k}^{\textrm{an,\#}}(x_1,x_2)= & \left| \begin{array}{l} 0 \\  
\psi_1(2^{j_1}x_1-k_1)\psi_0(2^{j_2}x_2-k_2) \end{array} \right.  
\end{array}  
$$  
then  
$$ 
\v_J^{\#}=\sum_{j_1,j_2<J}\sum_{\k\in \Z^2} \left( 
d^{\textrm{an}}_{1,\j,\k} ~\Psi_{1,\j,\k}^{\textrm{an,\#}}  
+d^{\textrm{an}}_{2,\j,\k} ~\Psi_{2,\j,\k}^{\textrm{an,\#}} \right) 
$$  
Thus, to compute the expansion of $P_J^{\#}\v$ in terms of the gradient vector wavelets, we have to perform the change of basis:  
$$  
\begin{array}{ccl}  
\left\{ \begin{array}{c} \Psi_{1,\j,\k}^{\textrm{an,\#}} \\ \Psi_{2,\j,\k}^{\textrm{an,\#}} \end{array}  
\right.& \longrightarrow &  
\left\{ \begin{array}{ll} \Psi_{\textrm{curl},\j,\k}^{\textrm{an}}  
&=2^{j_1} \Psi_{1,\j,\k}^{\textrm{an,\#}}+2^{j_2} \Psi_{2,\j,\k}^{\textrm{an,\#}} \\  
\Psi_{N,\j,\k}^{\textrm{an}} &=2^{j_2} \Psi_{1,\j,\k}^{\textrm{an,\#}}-2^{j_1} \Psi_{2,\j,\k}^{\textrm{an,\#}}  \end{array} \right.  
\end{array}  
$$  
which leads to:  
\begin{equation} 
\label{vcurl}   
\v_J^{\#}=\sum_{j_1,j_2<J}\sum_{\k\in \Z^2} \left( 
d^{\textrm{an}}_{\textrm{curl},\j,\k} ~\Psi_{\textrm{curl},\j,\k}^{\textrm{an}}  
+d^{\textrm{an}}_{N,\j,\k} \Psi_{N,\j,\k}^{\textrm{an}} \right) 
\end{equation}  
where the curl-free wavelet coefficients are obtained from the standard ones by: 
\begin{equation} 
\label{dcurl} 
d_{\textrm{curl},\j,\k}^{\textrm{an}}=  
\frac{2^{j_1}}{2^{2j_1}+2^{2j_2}}~d_{1,\j,\k}^{\textrm{an}}+\frac{2^{j_2}}{2^{2j_1}+2^{2j_2}}~d_{2,\j,\k}^{\textrm{an}}   
\end{equation} 
associated to complement coefficients: 
\begin{equation} 
\label{dN} 
d_{N,\j,\k}^{\textrm{an}}=  
\frac{2^{j_2}}{2^{2j_1}+2^{2j_2}}~d_{1,\j,\k}^{\textrm{an}}-\frac{2^{j_1}}{2^{2j_1}+2^{2j_2}}~d_{2,\j,\k}^{\textrm{an}}  
\end{equation}  
  
\ 
  
\subsection{Implementation of the Hodge decomposition in the wavelet context}  

\label{Hodge}
  
From now on, our objective is to compute the wavelet decomposition of a given vector function $\v$: this means
to find a divergence-free component $\v_{\mbox{div}}$ and an orthogonal curl-free  
component $\v_{\mbox{curl}}$ such that: 
$$  
\v=\v_{\mbox{div}} + \v_{\mbox{curl}}  
$$  
where:  
$$  
\v_{\mbox{div}}=\sum_{\j,\k} d_{\textrm{div},\j,\k} \Psi_{\textrm{div},\j,\k} \qquad 
\v_{\mbox{curl}}=\sum_{\j,\k} d_{\textrm{curl},\j,\k} \Psi_{\textrm{curl},\j,\k} 
$$  
are the wavelet expansions onto div-free and curl-free wavelet bases constructed previously (section \ref{div2Dan}
and \ref{gradient}). For more simplicity, we will focus on 2D anisotropic wavelet bases (and we will omit the superscript
''an'' in the notation of the basis functions). 
 
\ 
 
To provide with such decomposition, we have to overcome two problems:  
 
- The first one lies in the fact that div-free wavelets and curl-free wavelets form {\it biorthogonal} bases in
their respective spaces, and as already noticed by K. Urban \cite{U00}, they would not give rise, in a simple way,
to the orthogonal projections $\v_{\mbox{div}}$ and $\v_{\mbox{curl}}$ of $\v$. As solution, we propose to
construct, in wavelet spaces, two sequences $(\v_{\mbox{div}}^p)$ and $(\v_{\mbox{curl}}^p)$ that will converge to
$\v_{\mbox{div}}$ and $\v_{\mbox{curl}}$. 
 
- The second difficulty is that div-free wavelets leave in spaces of the form $(V^1_J\otimes V^0_J)\times
(V^0_J\otimes V^1_J)$, whereas curl-free wavelets arise from $(\tilde V^0_J\otimes \tilde V^1_J)\times (\tilde
V^1_J\otimes \tilde V^0_J)$, where $V^0_0, V^1_0$ and $\tilde V^0_0, \tilde V^1_0$ are couples of spaces related by
differentiation and integration. These spaces are different, and in order to construct our approximations
$(\v_{\mbox{div}}^p)$ and $(\v_{\mbox{curl}}^p)$, we have to define a precise interpolation procedure between the
two kinds of spaces. In particular, the spaces $\tilde V^0_0, \tilde V^1_0$ can be suitably chosen from $V^0_0,
V^1_0$. 
 
\subsubsection{Iterative construction of the div-free and curl-free parts of a flow} 
\label{iterhdg}

Let $\v=(v_1,v_2)$ a vector function be given, and suppose that $\v$ is periodic in
both directions, and known on $2^J\times 2^J$ grid points that are not necessarily
the same for $v_1$  and $v_2$.
In the following, we
will note:

- $\I_J \v$ an approximation of $\v$ in the space $(V^1_J\otimes V^0_J)\times
(V^0_J\otimes V^1_J)$, given by some interpolating process.

- $\I_J^{\#} \v$ an approximation of $\v$ in the space $(\tilde V^0_J\otimes \tilde
V^1_J)\times (\tilde V^1_J\otimes \tilde V^0_J)$, also given by some interpolating
process.

\

\noindent We now define the sequences $\v_{\mbox{div}}^p \in (V^1_J\otimes V^0_J)\times
(V^0_J\otimes V^1_J)$ satisfying $\mbox{div}~\v_{\mbox{div}}^p=0$, and
$\v_{\mbox{curl}}^p \in (\tilde V^0_J\otimes \tilde V^1_J)\times (\tilde
V^1_J\otimes \tilde V^0_J)$ satisfying $\mbox{curl}~\v_{\mbox{curl}}^p=0$, as
follows:

- We begin with $\v^0=\I_J \v$ and we compute  $\v_{\mbox{div}}^0$, the
divergence-free wavelet decomposition of $\v^{0}$, and its complement $\v_{n}^0$,
by formula (\ref{vdiv},\ref{ddiv}):
$$
\I_J \v = \v_{\mbox{div}}^0 +
\v_{n}^0=\sum_{\j,\k}d_{\textrm{div},\j,\k}^{0}~\Psi_{\textrm{div},\j,\k}+
\sum_{\j, \k}
d_{n,\j,\k}^{0}~\Psi_{n,\j,\k}
$$
Then we compute at grid points the difference $\v-\v_{\mbox{div}}^0$.\\
Secondly we consider  $\I_J^{\#}(\v-\v_{\mbox{div}}^0)$, and we apply the curl-free wavelet decomposition
(\ref{vcurl},\ref{dcurl},\ref{dN}), leading to a curl-free part and its complement:
$$
\I_J^{\#} (\v-\v_{\mbox{div}}^0)=\v_{\mbox{curl}}^0+\v_{N}^0 = \sum_{\j,
 \k}
d^{0}_{\textrm{curl},\j,\k} ~\Psi_{\textrm{curl},\j,\k}
+ \sum_{\j, \k} d^{0}_{N,\j,\k} \Psi_{N,\j,\k} 
$$
Finally  we define {\it pointwise}:  $\v^{1}=\v- \v_{\mbox{div}}^0 -\v_{\mbox{curl}}^0$.

\

- At step $p$, by knowing $\v^{p}$ at grid points, we are able to construct a divergence
free-part $\v_{\mbox{div}}^p$ of $\I_J \v^{p}$ by (\ref{vdiv}), and $\v_{\mbox{curl}}^p$, the
curl-free component of
$ \I_J^{\#} (\v^{p} - \v_{\mbox{div}}^p)$ by (\ref{vcurl}) ($\v^{p} -
\v_{\mbox{div}}^p$ being computed at grid points). The next term of
the sequence is again defined {\it pointwise}:
\begin{equation}
\label{vp}
\v^{p+1}=\v^{p}-\v_{\mbox{div}}^p -\v_{\mbox{curl}}^p
\end{equation}

\

We iterate this process until $\|\v^{P}\|_{\ell^2}<\epsilon$, and we
obtain:
\begin{eqnarray*}
\v&\approx_{\epsilon} & \sum_{p=1}^P \v_{\mbox{div}}^p+\sum_{p=1}^P \v_{\mbox{curl}}^p \\
&= & \sum_{\j,\k} \left(\sum_{p=1}^P
d_{\textrm{div},\j,\k}^{p}\right)~\Psi_{\textrm{div},\j,\k} +
\sum_{\j, \k} \left(\sum_{p=1}^P
d^{p}_{\textrm{curl},\j,\k} \right) ~\Psi_{\textrm{curl},\j,\k}
\end{eqnarray*} 
where the right hand side is an approximation of $\v$, which interpolates the data
up to an error $\epsilon$ ($\epsilon$ being given). 

For the moment, we are not able to prove theoretically the convergence to $0$ of the
sequence
$(\v^{p})$: we will prove it experimentally in section \ref{hodge-exp}, on
arbitrary fields.
Nevertheless, we can outline some remarks:

- The convergence rate depends on the choice of complement functions
$\Psi_{n,\j,\k}$, $\Psi_{N,\j,\k}$. The more the $L^2$-scalar products  
$<\Psi_{\textrm{div},\j,\k},\Psi_{\textrm{n},\j',\k'}>$ and
$<\Psi_{\textrm{curl},\j,\k},\Psi_{\textrm{N},\j',\k'}>$ are small, the faster the
sequence converges.

- Ideally, we would like to choose the interpolating operators $\I_J$ and $\I_J^{\#}$ 
such that the convergence doesn't depend on
this choice. We propose below a choice for these operators, based on
spline-quasi interpolation, which is satisfactory at relatively-slow convergence rate.

\subsubsection{Hodge-adapted spline interpolation}

In this part, we will detail our choice of operators $\I_J $ and
$\I_J^{\#}$, in the context of spline spaces of degree 1 ($V^0_j$) and
2 ($V^1_j$) that we
have introduced at the beginning. 

\begin{figure} [!h]
\begin{center} 
\input{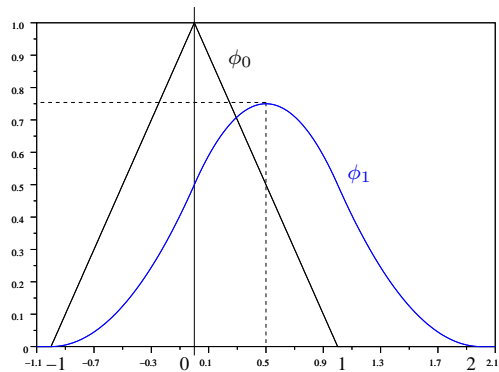} 
\caption{\label{g0g1} The two scaling functions of  $V_0$ and $V_1$,
and their symmetry centers}
\end{center} 
\end{figure}

Let us suppose the components $v_1$ and $v_2$ of 
a velocity field $\v$ be known respectively at  knot points 
$2^{-J}(k_1+\frac{1}{2},k_2)$ and $2^{-J}(k_1,k_2+\frac{1}{2})$, for
$k_1,k_2=0. 2^J-1$. This
choice of grid is induced by the symmetry centers  of scaling functions
$\phi_0$ of $V_0$ and $\phi_1$ of $V_1$ (see Figure \ref{g0g1}).

For $J$ given, $\I_J $ is chosen as an operator of quasi-interpolation (similarly to section
\ref{sec:quasi})  in the spline space $(V^1_J\otimes V^0_J)\times (V^0_J\otimes V^1_J)$
 
$$ 
\mathbb{I}_J \v=\sum_{\k}c_{\k}^1~\Phi_{\textrm{1},J,\k}+\sum_{\k}c_{\k}^2~\Phi_{\textrm{2},J,\k}
$$ 
where $\Phi_{\textrm{1}}$ and $\Phi_{\textrm{2}}$ are the vector scaling functions
introduced in section \ref{isotropic2D}.

The second operator  $\I_J ^{\#}$ provides again with
a quasi-interpolation  of vector functions onto a new  
spline space $(\tilde V^0_J\otimes \tilde V^1_J)\times (\tilde V^1_J\otimes \tilde
V^0_J)$. Under interpolating considerations, we define:
$$\tilde V^0=\{\v~;~\v(x-1/2)\in
V_0\}=span\{\phi_0(x-1/2-k)~;~k\in\Z\}$$
$$\tilde V^1=\{\v~;~\v(x-1/2)\in V_1\}=span\{\phi_1(x-1/2-k)~;~k\in\Z\}$$
Hence we can write: 
$$ 
\mathbb{I}^{\#}_J\v=\sum_{\k} c_{\k}^{\# 1}~\tilde \Phi_{\textrm{1},J,\k}+\sum_{\k}c_{\k}^{\# 2}~\tilde \Phi_{\textrm{2},J,\k}
$$ 
where $\tilde \Phi_{\textrm{1},J,\k}$ and $\tilde \Phi_{\textrm{2},J,\k}$
are the 2D anisotropic vector scaling functions of \ref{isotropic2D}
built from $\tilde\phi_0=\phi_0$ and $\tilde\phi_1=\phi_1$.


\section{Numerical experiments} 
 
In this section, we present our numerical results concerning the application of divergence-free wavelet 
decomposition, for analyzing several data. We begin with the analyses of periodic, numerical, incompressible 
velocity fields in dimensions two and three, arising from pseudo-spectral codes. First, we have to take care of 
the initial interpolation of such fields, in order not to break the incompressible condition satisfied in Fourier 
space. Then, after the vizualisation of the divergence-free wavelet coefficients, we study the compression 
obtained through the wavelet decomposition. In the last part, we investigate and numerically prove  the 
convergence of the algorithm presented in section \ref{hodge}, which provides  with the wavelet Hodge 
decomposition of any flow.  As an example, we compute the div-free
component of the nonlinear term of
the Navier-Stokes equations, and we extract the associated pressure,
directly in wavelet space.
In all the experiments, we will use divergence-free wavelets constructed with splines of degrees 1 and 2.  
  
\subsection{Approximation of the velocity in spline spaces} 
 
Usually, the data are provided by point values of the velocity field. The first step
of the wavelet decomposition consists 
in interpolating the velocity coordinates on the suitable B-spline space. The arising problem is 
that this approximation may not conserve the divergence free condition that is verified 
in Fourier space, when velocity arise from a spectral code. \\ 
The spline approximation of data, obtained through spectral methods,  introduces a slight error 
for the divergence free condition. This difference may not be neglectable. 
For the turbulent fields we studied (2D and 3D) the error is about $1 \%
$ of the $L^2$-norm (i.e. 0.01 \% of the energy). 
 
Thus we propose two ways to overcome the problem.
The first way is to interpolate the velocity in the Fourier domain and to compute exactly its
biorthogonal projection on wavelet spaces.
The second way is to interpolate on the divergence free 
B-spline spaces with a Hodge decomposition made by wavelet decompositions as it is proposed 
in the part \ref{iterhdg} . This can be applied to any compressible flow. 
 
\subsubsection{By quasi-interpolation}
\label{sec:quasi} 
 The spline quasi-interpolation is a good compromise when we have to
deal simultaneously with 
spline approximations of degree even and odd. In this context, the
order of  approximation is $n+1$, by using B-splines 
of degree $n$ \cite{Bo}. An advantage of the procedure is that it  may be applied in any case
of boundary conditions.
 
Let $b$  be a B-spline scaling function ($b=\phi_0$ or $\phi_1$). Given the sampling $f(k/N)$ ($N=2^J$), we want to compute scaling
coefficients $c_k$, of a spline function $f_N$, that will nearly interpolate
the values $f(k/N)$: 
\begin{equation}
\label{quasi}
f_N(x)=\sum_{k\in\Z}c_k~b(Nx-k) 
\end{equation}
$f_N$ is an interpolating function if:
$$ 
\sum_{k\in\Z}c_k ~b(\ell-k)=f(\frac{\ell}{N}) \quad \forall \ell\in \Z 
$$ 
For example, if we consider $b=\phi_1$ (spline of degree 2), the
previous condition implies: 
$$ 
f_N(\frac{\ell}{N})=\frac{1}{2}(c_{\ell-1}+c_\ell)=f(\frac{\ell}{N}) \quad \forall \ell\in \Z 
$$ 
In order to avoid the inversion of a linear system, the
quasi-interpolation introduces, instead of $c_\ell$:
\begin{eqnarray*}
\widetilde{c}_\ell=\frac{5}{8}[f(\frac{\ell}{N})+f(\frac{\ell+1}{N})]-\frac{1}{8}[f(\frac{\ell-1}{N})+f(\frac{\ell+2}{N})] 
\quad \forall \ell\in \Z 
\end{eqnarray*}
By replacing $c_\ell$ by $\widetilde{c}_\ell$ in (\ref{quasi}), we
obtain the following error at each grid point:
\begin{eqnarray*} 
\frac{1}{2}(\widetilde{c}_{\ell-1}+\widetilde{c}_\ell)-f(\frac{\ell}{N})&=&\frac{1}{16}[-f(\frac{\ell-2}{N}) \\
& &+4f(\frac{\ell-1}{N})-6f(\frac{\ell}{N})+4f(\frac{\ell+1}{N})-f(\frac{\ell+2}{N})]\\ 
&=&-\frac{1}{48N^4}f^{(4)}(\theta)+O(\frac{1}{N^6}) \quad ,
\textrm{with} \quad \theta \in ]\frac{\ell-2}{N},\frac{\ell+2}{N}[ 
\end{eqnarray*} 
Therefore, the pointwise error of quasi-interpolation is order 4, for a sufficiently regular function. 
 
\subsubsection{By using the Discrete Fourier Transform} 
 
  Since they are highly accurate, spectral methods are often considered as a reference technique for  
 simulating incompressible turbulent flows. For periodic boundary conditions on the cube $[0,1]^2$, the  
 Discrete Fourier Transform is used to decompose the velocity $\u$. \\  
If $\hat{\u}_{\k}$ means the Discrete Fourier coefficients of $\u$  on a $N^2$ regular grid,  
 $$  
 \hat{\u}_{\k}=\frac{1}{N^2}\sum_{\n\in\{0, 1, \dots, N-1\}^2} \u(\frac{\n}{N})~  
 e^{-2i\pi\frac{\k.\n}{N}}  
 $$   
 the velocity expansion in the Fourier exponential basis is:  
 \begin{equation}  
 \label{serieF} 
 {\u}{(\x)}=\sum_{\k\in\{0, 1, \dots, N-1\}^2}\hat{ \u}_{\k}~e^{2i\pi\k.\x}  
\end{equation} 
  In this context, the divergence-free condition $\mbox{div}~ \u = 0$ 
writes:  
 \begin{equation}  
 \label{Fdiv}  
 \k.\hat{\u}_{\k}=0~,~~\forall \k\in\{0, 1, \dots, N-1\}^2  
 \end{equation}  
   
\  
  
Assume now that the velocity field $\u$ we have to analyze (supposed to be 1-periodic in both  
directions), is obtained from a spectral method and verifies the incompressibility condition in Fourier  
domain (\ref{Fdiv}). To compute its decomposition in a divergence-free wavelet basis of $\R^2$, we have  
first to approximate $\u= (u_1,u_2)$ in the suitable space 
$(V^1_J\otimes V^0_J)\times (V^0_J\otimes V^1_J)$  
which has been introduced in section \ref{isotropic2D}, where $J$ corresponds to $N=2^J$. Then we search for an approximate  
function $\u_J=(u_{J~1},u_{J~2})$ such that:  
$$  
\left\{  
\begin{array}{lll}  
\u_{J~1}&=&\sum_{n_1=0}^{2^J-1}\sum_{n_2=0}^{2^J-1} c^1_{J,n_1,n_2}~\phi_{1,J,n_1}\phi_{0,J,n_2}\\  
&&\\  
\u_{J~2}&=&\sum_{n_1=0}^{2^J-1}\sum_{n_2=0}^{2^J-1} c^2_{J,n_1,n_2}~\phi_{0,J,n_1}\phi_{1,J,n_2}  
\end{array}  
\right.  
$$  
  
For the choice of functions $\phi_{0}$ and $\phi_{1}$ defined above (see equation (\ref{derb})),  
the incompressibility condition $\mbox{div}~\u_J=0$ takes the discrete form on the coefficients $c^i_{J,n_1,n_2}$:  
\begin{equation}  
\label{Cdiv}  
c^1_{J,n_1,n_2}-c^1_{J,n_1+1,n_2} + c^2_{J,n_1,n_2} - c^2_{J,n_1,n_2+1}=0~,~~\forall (n_1,n_2)  
\end{equation}  
  
\  
  
To conserve the incompressibility condition verified by $\u$, a solution consists in considering $\u_J$ as  
the biorthogonal projection onto the space $(V^1_J\otimes V^0_J)\times (V^0_J\otimes V^1_J)$, since we know that this  
projector commutes with the partial derivatives \cite{L92}. This is equivalent to consider that:  
  
$$  
\left\{  
\begin{array}{lll}  
c^1_{J,n_1,n_2}&=&<\u~|~\phi^*_{1,J,n_1}\phi^*_{0,J,n_2}>$$\\  
&&\\  
c^2_{J,n_1,n_2}&=&<\u~|~\phi^*_{0,J,n_1}\phi^*_{1,J,n_2}>$$  
\end{array}  
\right.  
$$  
Replacing $\u$ by its Fourier expansion (\ref{serieF}), it follows:  
\begin{eqnarray*}  
c^1_{J,n_1,n_2}&=&\sum_{\k\in\{0, 1, \dots, N-1\}^2}\hat{ \u}_{\k}~\iint_{\R^2}  
e^{2i\pi\k.\x}~\phi^*_{1,J,n_1}(x_1)~\phi^*_{0,J,n_2}(x_2) ~dx_1dx_2\\  
&=& 2^{-J}\sum_{\k\in\{0, 1, \dots, N-1\}^2}\hat{ \u}_{\k}~  
\overline{\hat{\phi}^*_{1}(2^{-J} 2\pi k_1)} ~~\overline{\hat{\phi}^*_{0}(2^{-J} 2\pi k_2)}~e^{2i\pi\frac{\k.\n}{2^J}}  
\end{eqnarray*}  
where $\hat{\phi}^*_{1}$, $\hat{\phi}^*_{2}$ denote the (continuous) Fourier transforms of the dual  
scaling functions ${\phi}^*_{1}$, ${\phi}^*_{2}$. Finally, we obtain an explicit form for the Discrete   
Fourier Transform $\mbox{DFT}$ of the coefficients $c^1_{J,n_1,n_2}$ (and in the same way for $c^2_{J,n_1,n_2}$):  
\begin{equation}  
\label{ck}  
\left\{  
\begin{array}{lll}  
\mbox{DFT}(c^1_{J,\n})_{\k}&=&\hat{ \u}_{\k}~2^{-J}~  
\overline{\hat{\phi}^*_{1}(2^{-J}(2\pi k_1))} ~~\overline{\hat{\phi}^*_{0}(2^{-J}(2\pi k_2))}\\  
&&\\  
\mbox{DFT}(c^2_{J,\n})_{\k}&=&\hat{ \u}_{\k}~2^{-J}~  
\overline{\hat{\phi}^*_{0}(2^{-J}(2\pi k_1))} ~~\overline{\hat{\phi}^*_{1}(2^{-J}(2\pi k_2))}  
\end{array}  
\right.  
\end{equation}  
It means that the discrete Fourier transform of coefficients $c^i_{J,n_1,n_2}$ is given by the 
discrete Fourier transform of $\u$, multiplied by tabulate values on $[0, 2\pi]$ of the Fourier 
transform of the duals $\hat{\phi}^*_{1}$, $\hat{\phi}^*_{2}$. In practice, we don't know the explicit 
forms of these functions, except by the infinite product: 
  
\begin{eqnarray*} 
\hat{\phi}^*_{0}(\xi)&=&\hat{\phi}_{0}(\xi)~\Pi_{j\geq 1}(2-\cos(\xi 2^{-j}))~= 
\left(\frac{\sin(\xi/2)}{\xi/2}\right)^2~\Pi_{j\geq 1}(2-\cos(\xi 2^{-j}))\\ 
\hat{\phi}^*_{1}(\xi)&=&\hat{\phi}_{0}^*~(\xi)\left(\frac{e^{i\xi}-1}{i\xi}\right)~=~ 
e^{i\xi/2}~\hat{\phi}_{0}^*(\xi)~\left(\frac{\sin(\xi/2)}{\xi/2}\right) 
\end{eqnarray*} 
Nevertheless, the infinite product converges rapidly, which allows to
obtain point values of $\hat{\phi}^*_{0}$ and $\hat{\phi}^*_{1}$,
with sufficiently accurate precision.
  
 \ 
  
In dimension 3, it goes similarly,  by considering the biorthogonal projection of a  3D vector field 
$\u$ onto the space $(V^1_J\otimes V^0_J\otimes V^0_J)\times (V^0_J\otimes V^1_J\otimes V^0_J)\times 
(V^0_J\otimes V^0_J\otimes V^1_J)$.  
  
\subsection{Analysis of 2D incompressible fields} 
 
We focus in this part on the analyses of  two-dimensional  decaying turbulent flows. \\ 
 
The first numerical experiment we present studies the merging of two same 
sign vortices. It concerns free decaying turbulence (no forcing 
term). The experiment was originally designed by M. Farge and N. Kevlahan \cite {SKF97}, and often 
used to test new models \cite{CP96, GK00}. This experiment was here reproduced by using a 
pseudo-spectral/finite-difference method, solving the Navier-Stokes equations in velocity-pressure 
formulation. 
 
\begin{figure}[!hb] 
\begin{center} 
\begin{tabular}{cccc} 
\includegraphics[width=3.5cm,height=3.5cm]{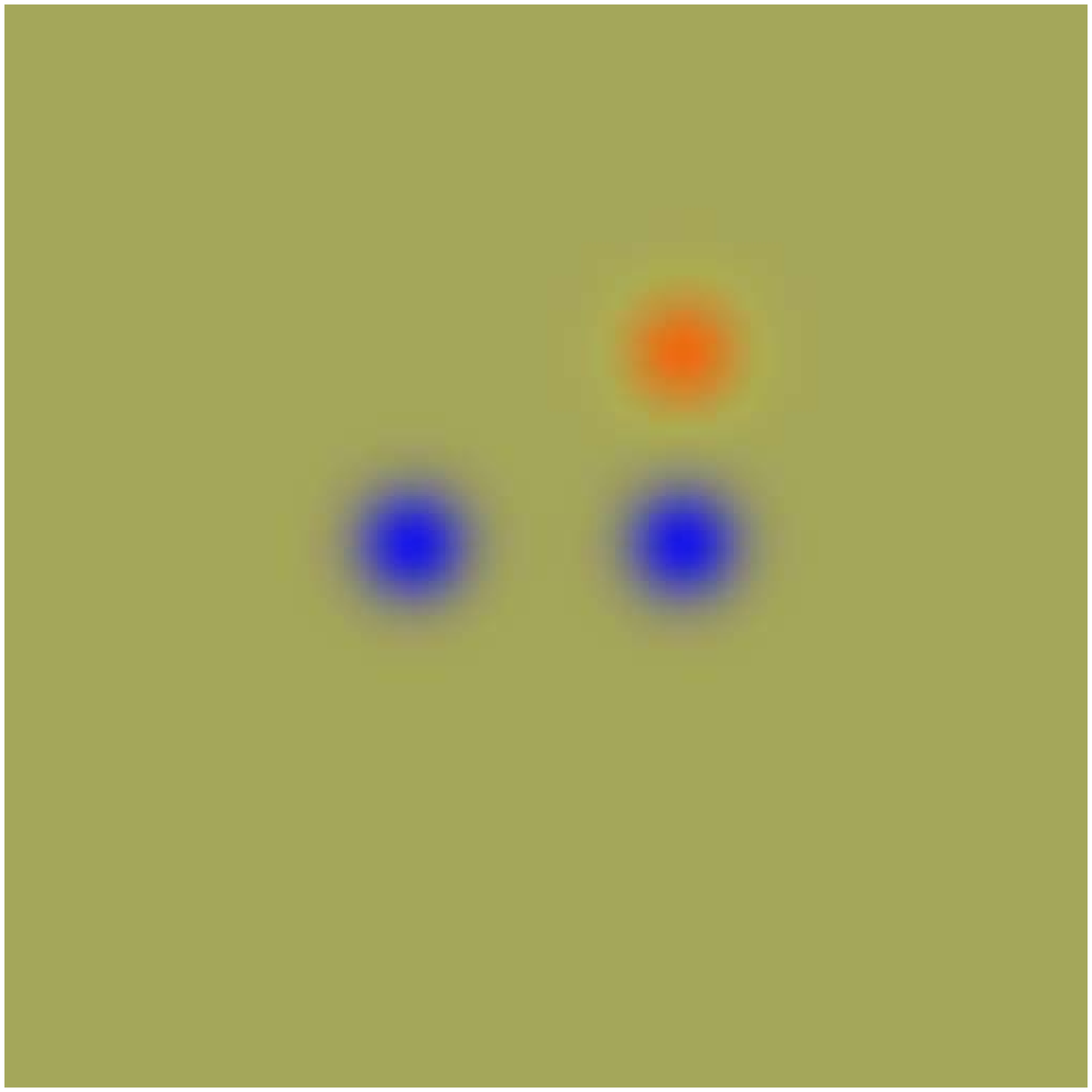}& 
\includegraphics[width=3.5cm,height=3.5cm]{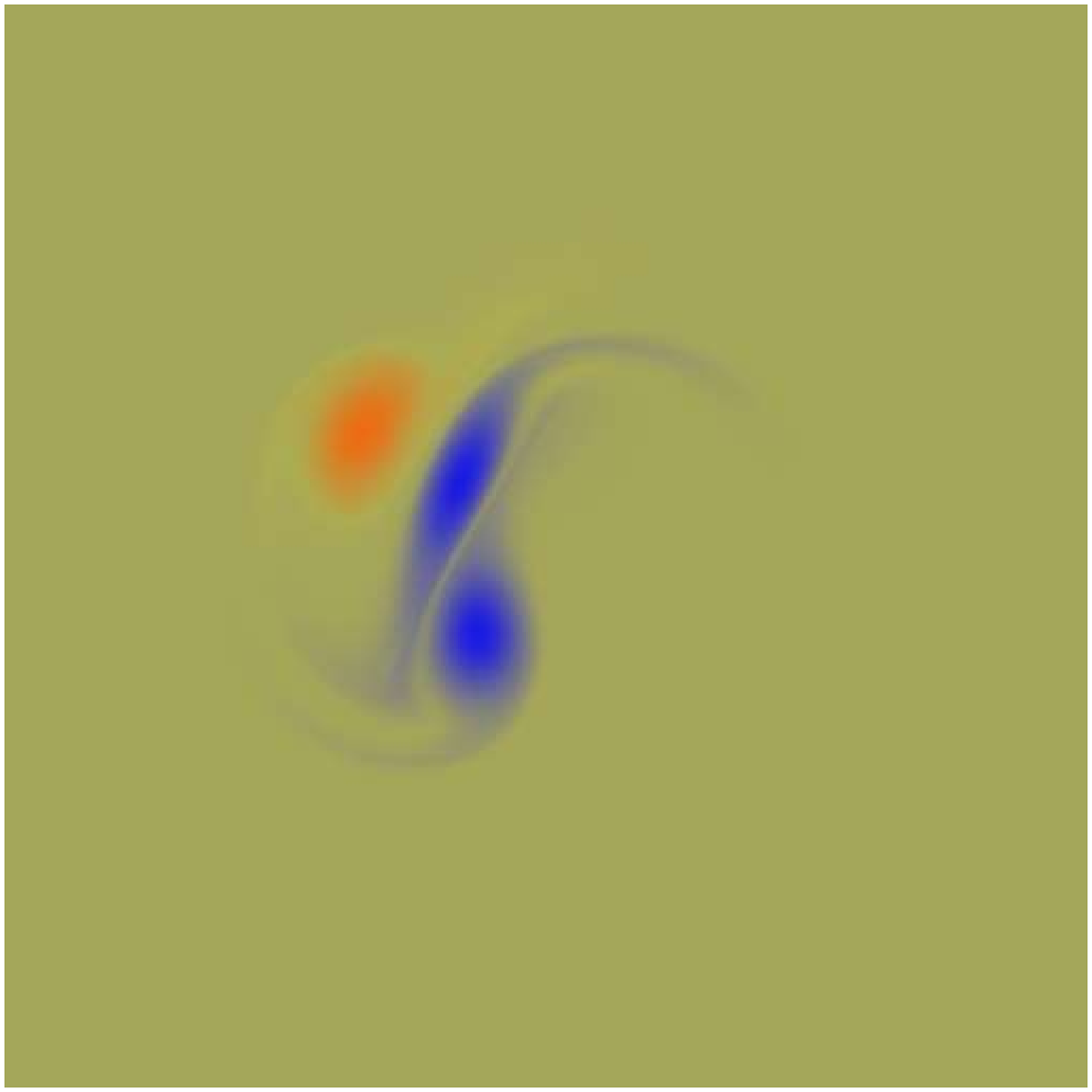}& 
\includegraphics[width=3.5cm,height=3.5cm]{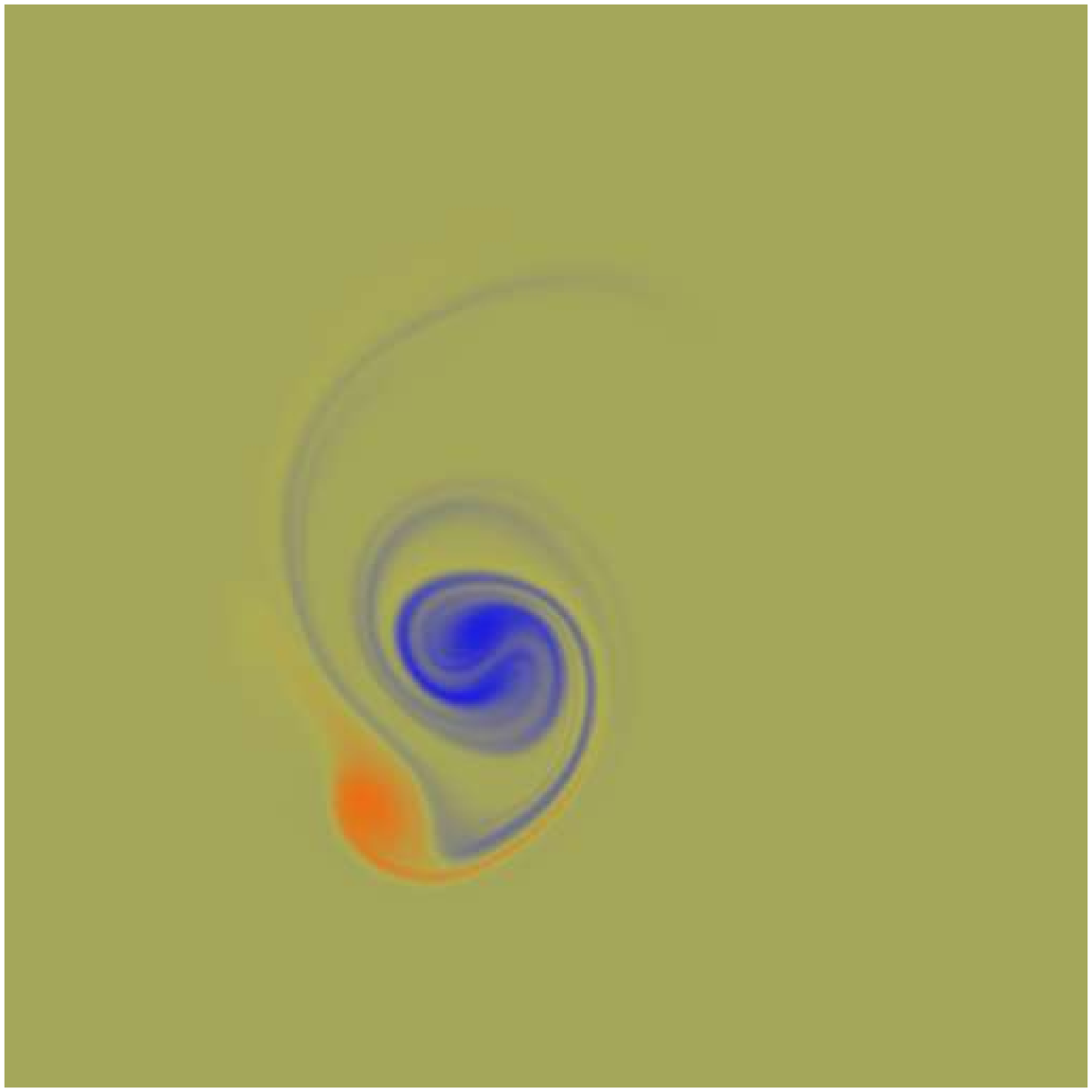}& 
\includegraphics[width=3.5cm,height=3.5cm]{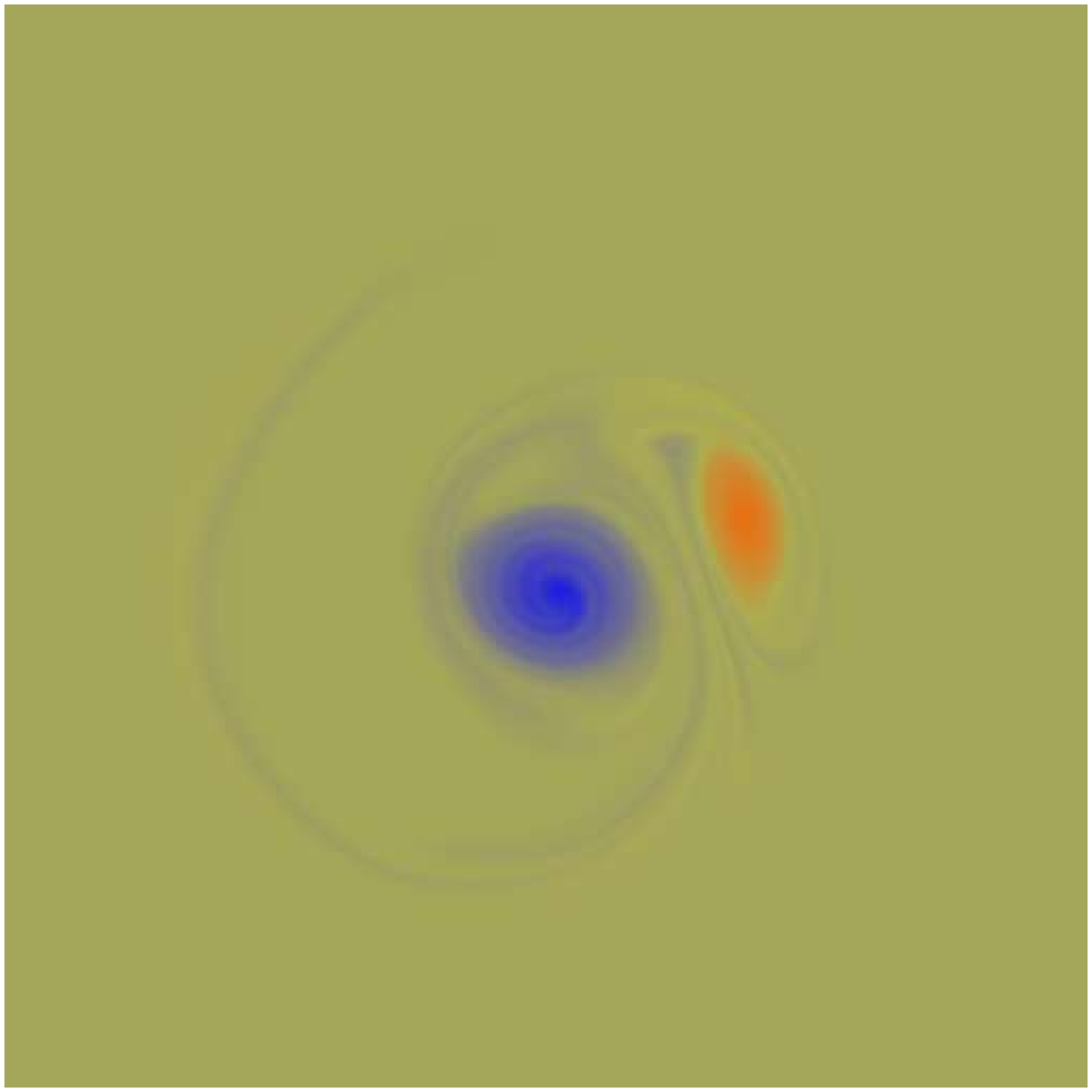}\\ 
\includegraphics[width=3.5cm,height=3.5cm]{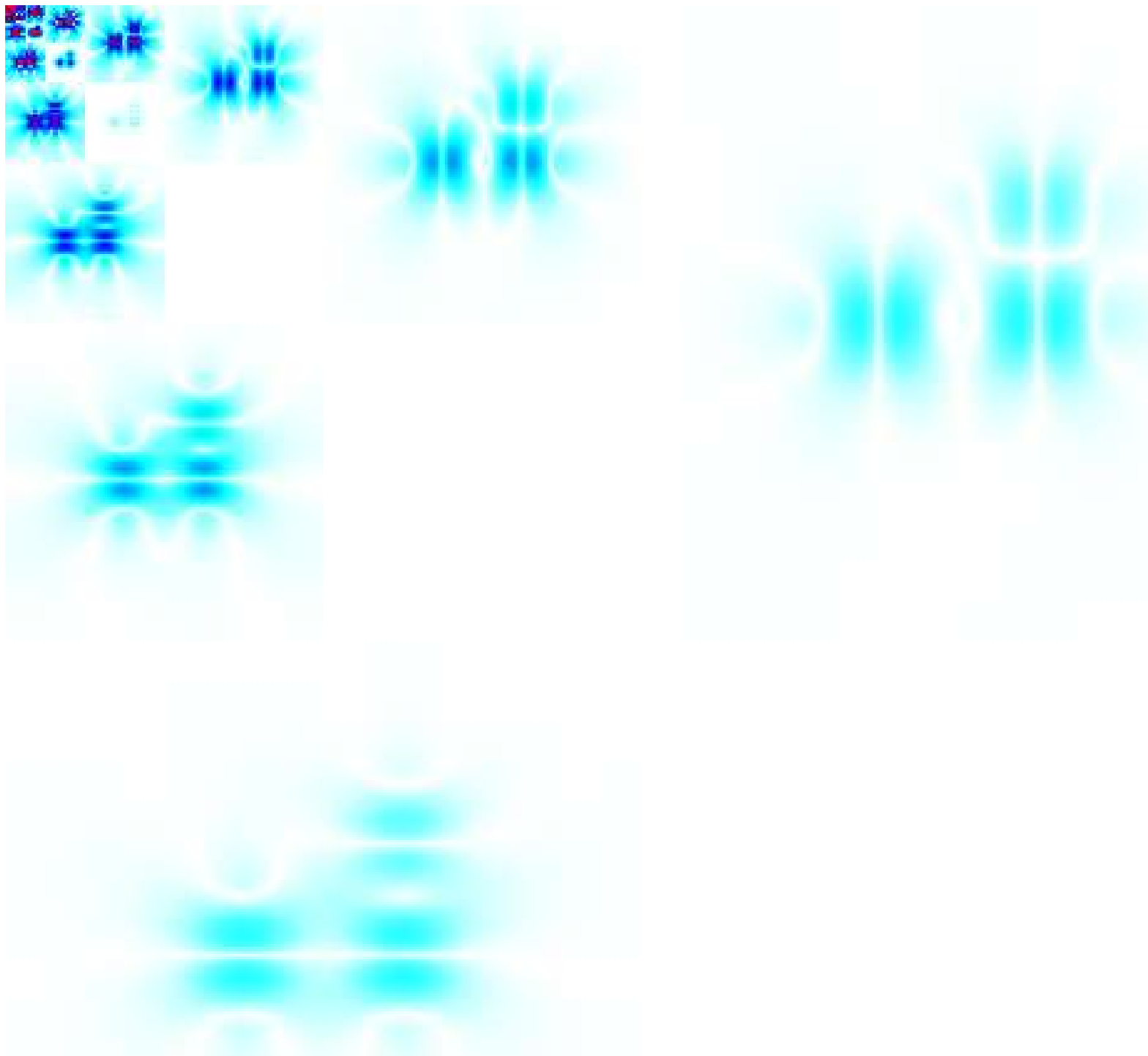}& 
\includegraphics[width=3.5cm,height=3.5cm]{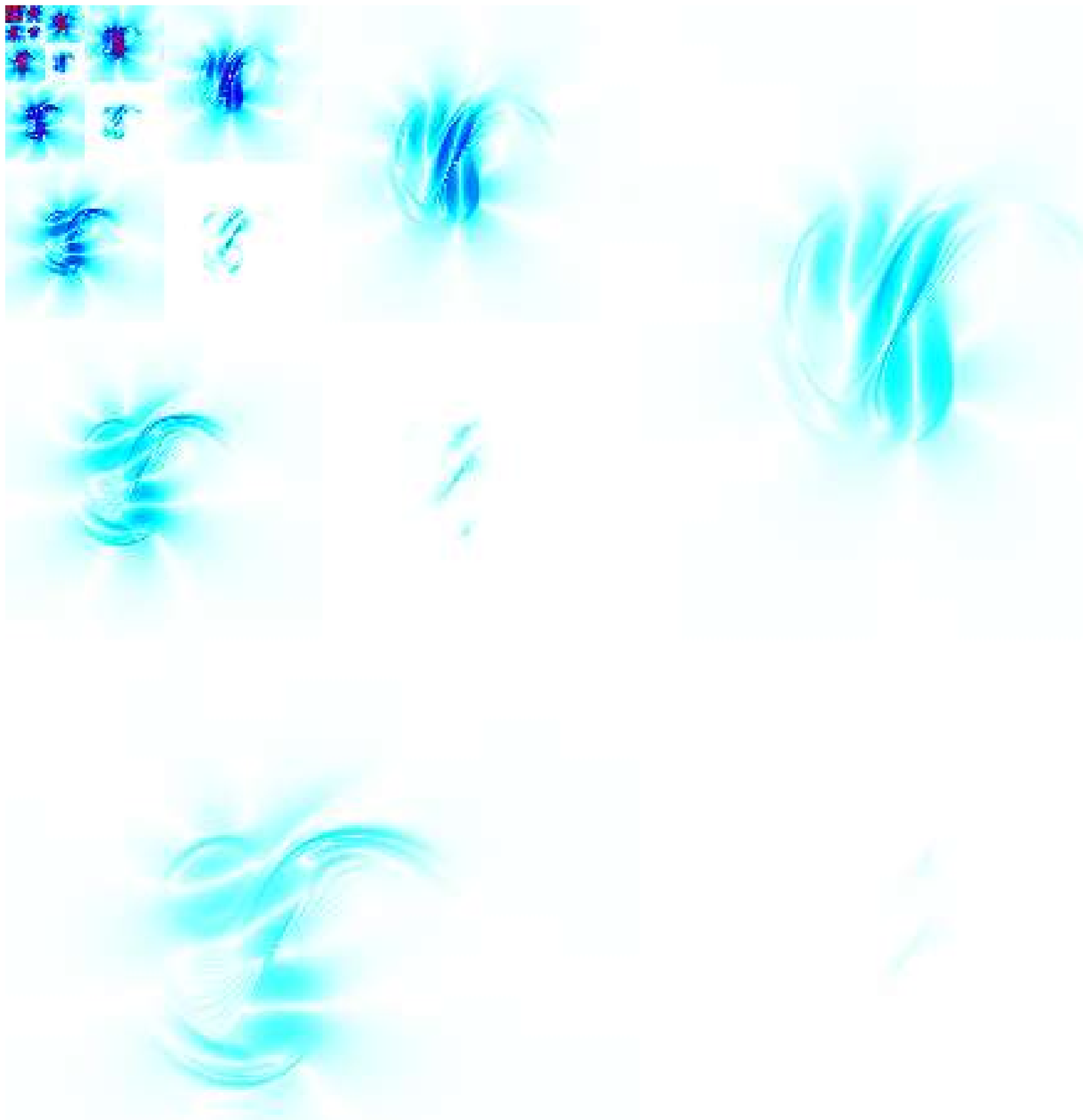}& 
\includegraphics[width=3.5cm,height=3.5cm]{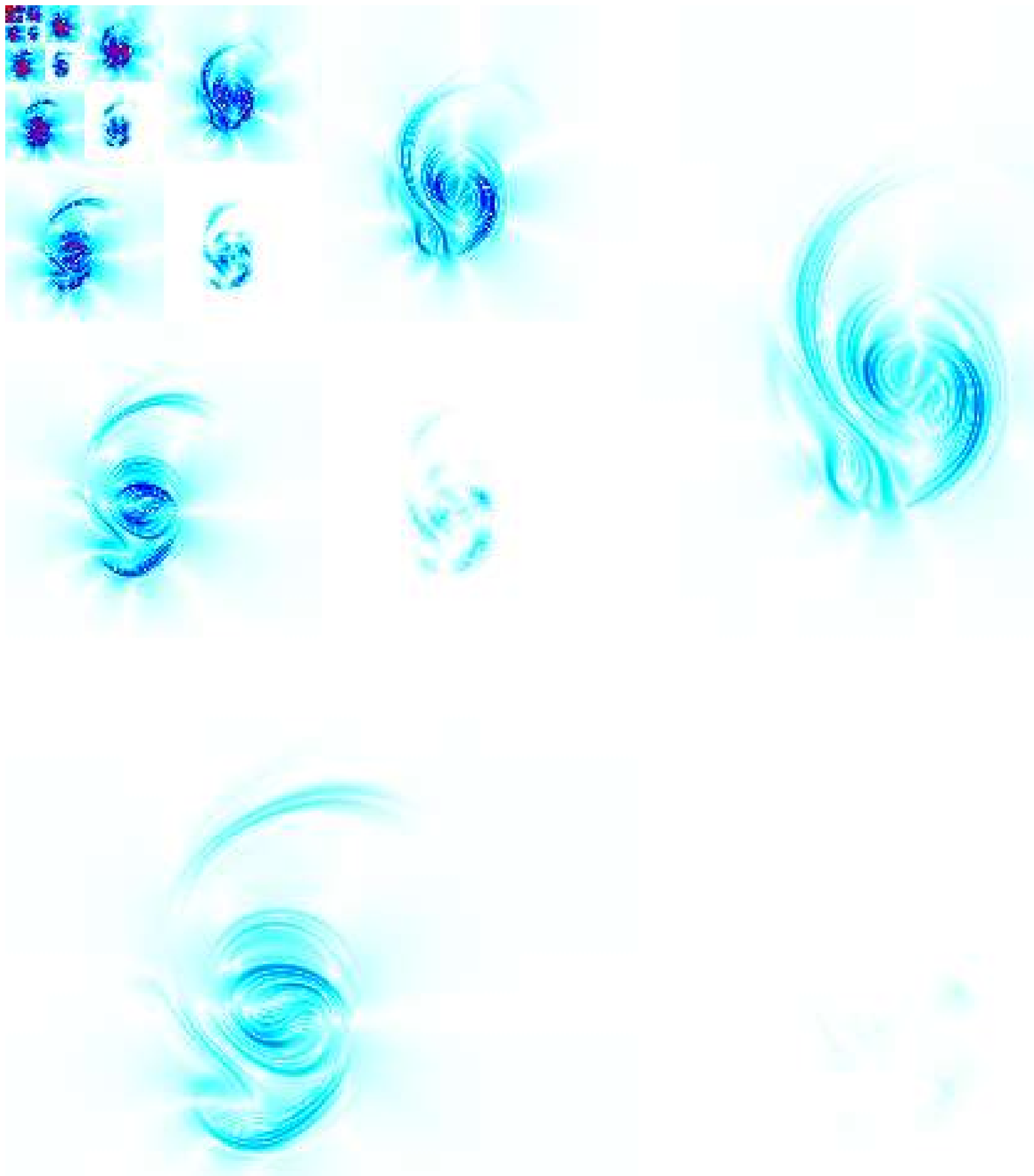}& 
\includegraphics[width=3.5cm,height=3.5cm]{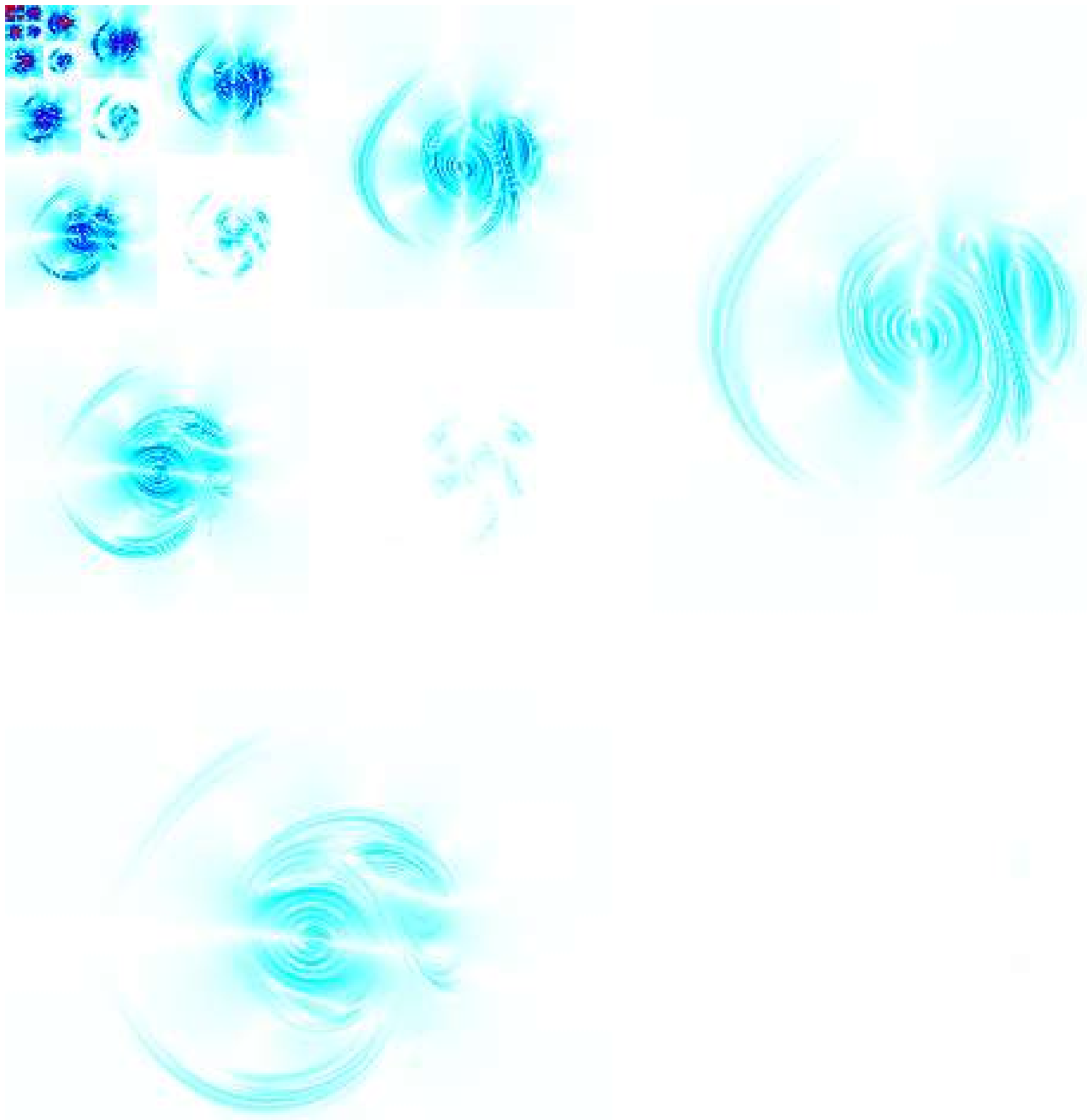} 
\end{tabular} 
\caption{\label{3tourb} Vorticity fields at times $t=0$, $t=10$, $t=20$ and $t=40$, and corresponding 
divergence-free wavelet coefficients of the velocity.} 
\end{center} 
\end{figure}

The initial state is displayed on figure  \ref{3tourb} left.  In a periodic box, three 
vortices with a gaussian vorticity profile are present; two are 
positive with the 
same intensity, one is negative with half the intensity of the 
others. The negative vortex is here to force the merging of the two 
positive ones.  
The time step was $\delta t=10^{-2}$ and the viscosity 
$\nu=5~10^{-5}$. The solution is computed on a
$512\times 512$ point grid. 
 
The vorticity fields at times $t=0$, $t=10$, $t=20$ and $t=40$ are 
displayed on figure \ref{3tourb}.  
The last row of figure \ref{3tourb} displays the absolute values of the isotropic divergence-free 
wavelet coefficients at corresponding times, renormalized by $2^j$ at scale index $j$. As one can see, 
divergence-free wavelet coefficients concentrate on strong change in vorticity zones, 
that is around or in between 
vortices, or along vorticity filaments, as they are equivalent to second derivatives of the 
velocity.\\

\

\begin{center} 
\begin{figure}[!hb]  
\begin{tabular}{cc}  
{\small Turbulent vorticity field} & {\small Divergence-free wavelet decomposition of the velocity} \\  
\includegraphics[scale=0.35]{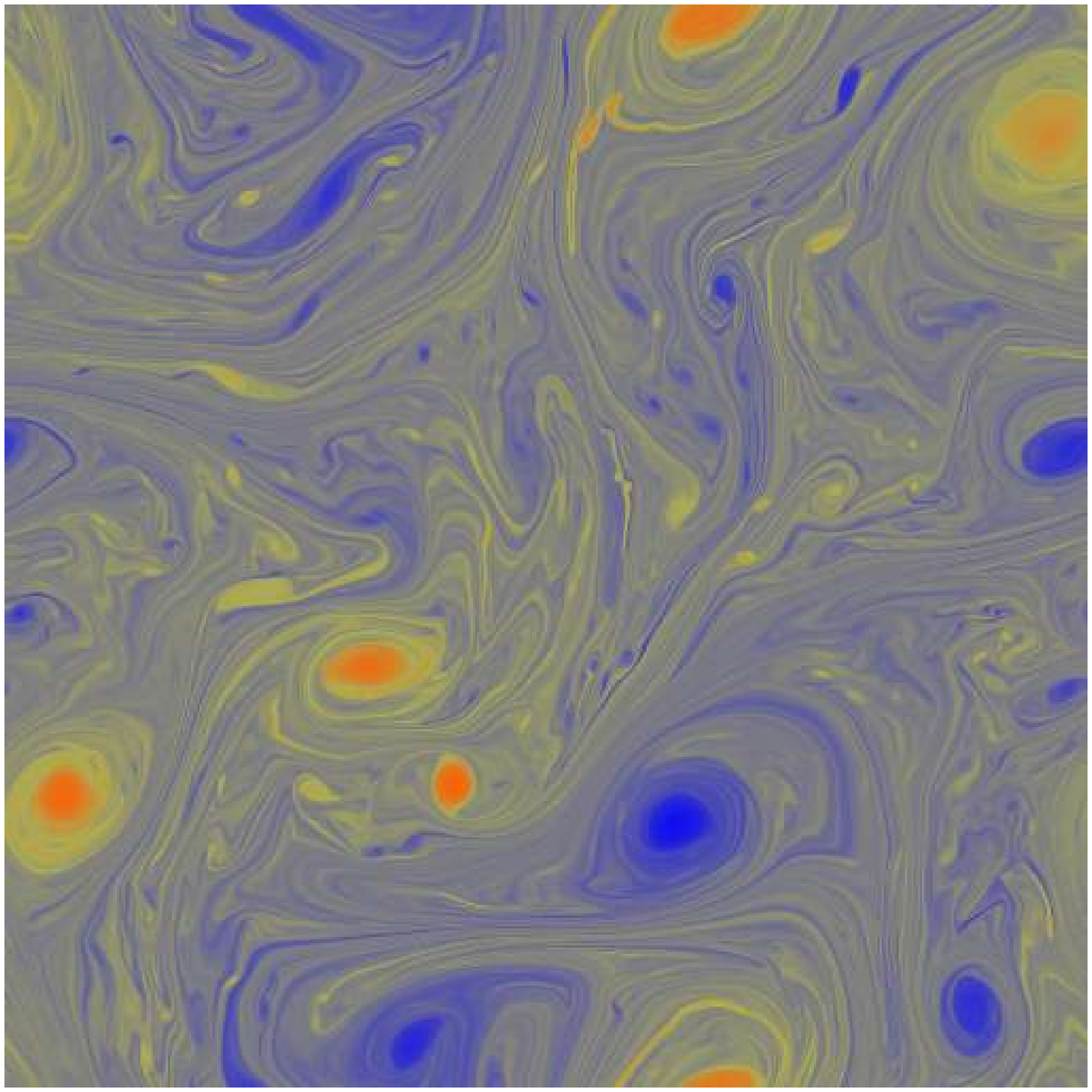} &  
\includegraphics[scale=0.35]{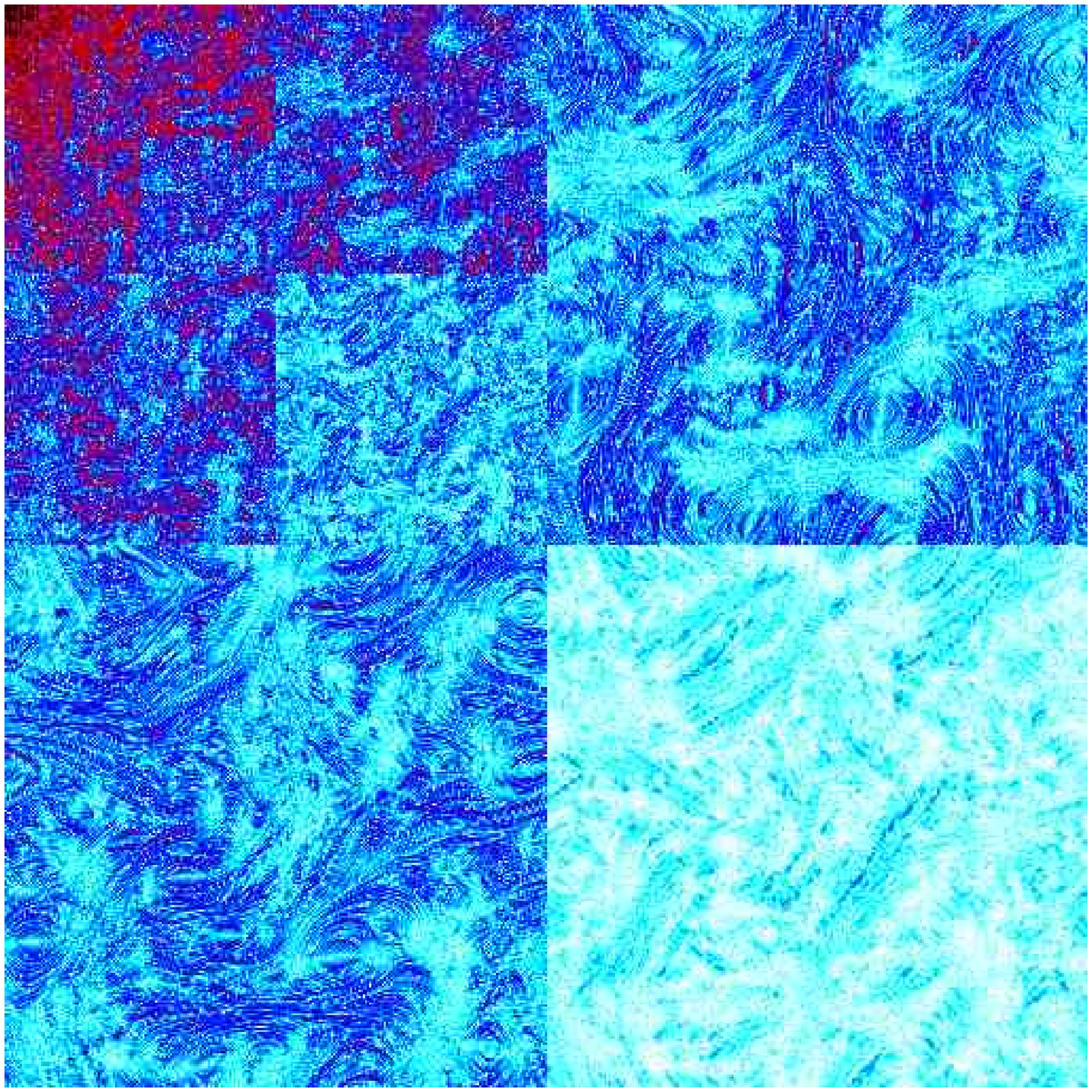} \\  
\end{tabular} 
\caption{\label{turb}Vorticity field for a $1024\times 1024$ simulation of decaying 
turbulence (left), and the corresponding divergence-free wavelet coefficients of the velocity field 
(right). }  
\end{figure} 
\end{center}

The second experiment deals with a decaying two-dimensional turbulent field, obtained with an initial 
state of random phase spectrum. That vorticity field was  computed with the spectral code of \cite{HH86}, 
at a resolution $1024\times 1024$, and a Reynolds number of $3.5\times 10^4$. This field has 
been kindly provided to us by G. Lapeyre \cite{L2000}. Figure \ref{turb} left represents the vorticity, after 
$40$ turnover time-scale, where it exhibits the emergence of coherent structures together with strong 
filamentation of the flow field outside the vortices.\\ 
We show in Figure \ref{turb} right, the isotropic divergence-free wavelet coefficients of the 
corresponding velocity field, in $L^{\inf}$-norm.

As expected, the wavelet coefficients get an insight into the energy distribution over the scales of 
the flow. As one can see on Figure  \ref{turb}, the energy at smallest scale (or highest wavenumbers) 
is localized along the strong deformation lines, and fits the filamentation between vortices, or 
with strong changes in vortices.  The top-right square corresponding to vertical isotropic 
wavelets ($\Psi_{\mbox{div},j,\k}^{(1,0)}$) exhibits vertical structures, whether the bottom-left 
square corresponding to horizontal wavelets ($\Psi_{\mbox{div},j,\k}^{(0,1)}$) exhibits horizontal 
deformation lines. 
 
\ 
\ 
 
\begin{figure}[!h]  
\begin{center}  
\includegraphics[width=7cm,height=4cm]{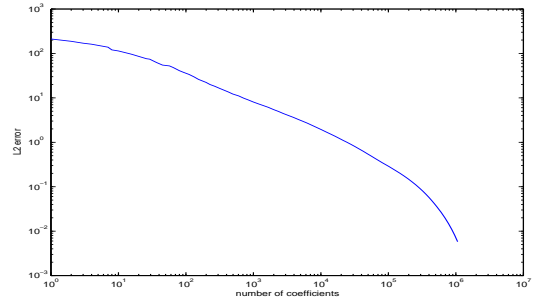}  
\end{center} 
\caption{\label{compTurb} $L^2$-error provided by the nonlinear N-best terms wavelet approximation (\ref{NL}): in 
log-log plot,  $L^2$-error (\ref{errorNL}) versus $N$ for a 2D turbulent flow. } 
\end{figure}

Now we investigate the compression properties of the divergence-free 
wavelet analysis: as predicted by the nonlinear approximation theory 
(see \cite{C00}), the compression ratio in energy-norm is governed by 
the underlying regularity of the velocity field in some Besov space. \\ 
Let $\u$ an incompressible field be given, its divergence-free wavelet expansion writes: 
$$  
\u = \u_0+\sum_{j\geq 0}\sum_{\k\in 
\Z^2}\left(d_{\textrm{div},j,\k}^{(1,0)}~\Psi_{\textrm{div},j,\k}^{(1,0)} 
+d_{\textrm{div},j,\k}^{(0,1)}~\Psi_{\textrm{div},j,\k}^{(0,1)}+  
d_{\textrm{div},j,\k}^{(1,1)}~\Psi_{\textrm{div},j,\k}^{(1,1)}\right)  
$$ 
The nonlinear approximation of $\u$ relies on computing the best N-terms wavelet approximation by 
reordering the wavelet coefficients: 
$$|d_{\textrm{div},j_1,\k_1}^{\varepsilon_1}|>|d_{\textrm{div},j_2,\k_2}^{\varepsilon_2}| > \dots 
>|d_{\textrm{div},j_N,\k_N}^{\varepsilon_N}| > \dots$$ 
and introducing  
\begin{equation} 
\label{NL} 
\Sigma_N(\u)=\u_0+ \sum_{i=1}^N d_{\textrm{div},j_i,\k_i}^{\varepsilon_i} 
~\Psi_{\textrm{div},j_i,\k_i}^{\varepsilon_i} 
\end{equation} 
Then we have 
\begin{equation} 
\label{errorNL} 
\| \u-\Sigma_N(\u)\|_{L^2}\leq C  \left(\frac{1}{N}\right)^{s} \| \u\|_{B^{s,q}_q} 
\end{equation} 
if the quantity 
$\| \u\|_{B^{s,q}_q}^q=\sum_{\varepsilon,j,\k} \left| d_{\textrm{div},j,\k}^{\varepsilon}\right|^q$ is 
finite, with $\frac{1}{q}=\frac{1}{2}+\frac{s}{n}$ (this means that $\u$ belongs to the Besov space 
$B^{s,q}_q$). As stated in \cite{C00}, the evaluated regularity $s$ can't be larger than the order
of polynomial reproduction in scaling spaces plus one (that equals the number of zero moments of the dual wavelet).
In our experiment, the dual spline wavelets $\psi^*_0$ and $\psi^*_1$
introduced in \ref{splines} have respectively two and three zero
moments, which only allows us to evaluate regularities smaller than two.\\ 
Figure \ref{compTurb} shows the nonlinear compression of 
divergence-free wavelets, provided on the $1024^2$ turbulent field. The curve represents the $L^2$-error  $\| 
\u-\Sigma_N(\u)\|_{L^2}$, versus $N$, in log-log plot. The convergence 
rate measured on the curve is $s\thickapprox 1.35$ , which induces that 
the velocity flow belongs to the corresponding Besov space $B^{s,q}_q$ 
with $q=0.85$.\\
When looking at the compression curve  on figure \ref{compTurb}, we observe three zones:

- First, large scale wavelets capture the large scale structure of the flows. Consequently, the compression
progresses slowly and irregularly.

- Then we observe a linear slope that represents the nonlinear
structure of the turbulent flows. In this region, we
are able to evaluate the regularity of the field.

- The last region corresponds to an abrupt decrease, due to the fact that the data are discrete.

One can also remark on Figure \ref{compTurb} that only $1.2$\%  of the
coefficients recover about
$99$\%  of the $L^2$-norm.
 
The same experiment was carried out on the three interacting vortices, but due to the few number of vanishing 
moments of the wavelets we use (2), the slope of the curve saturates at $s=2$, meaning that these fields
are more regular. 
 
\subsection{Analysis of a 3D incompressible field}  
  
\begin{figure}[!h] 
\begin{center}  
\includegraphics[width=6cm,height=6cm]{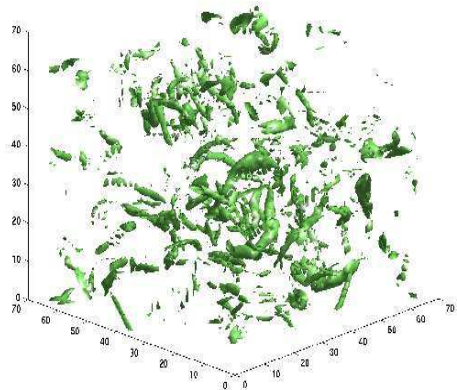}  
\end{center} 
\caption{\label{vort3D} Isosurface of vorticity magnitude after 5 large-eddy turnovers provided by a spectral method \cite{CMOV02}.} 
\end{figure}

In this part we consider a three dimensional periodic field, arising from a freely decaying isotropic 
turbulence, and kindly provided to us by G.-H. Cottet and B. Michaux \cite{CMOV02}. The experiment 
deals with an initial velocity condition of Gaussian distribution, and $128^3$ collocation points. 
Figure \ref{vort3D} displays the vorticity isosurfaces corresponding to about $40\% 
$ of the maximum vorticity at five turnover times.  
  
\begin{figure}[!h] 
\begin{center}  
\begin{tabular}{cc}  
\includegraphics[scale=0.35]{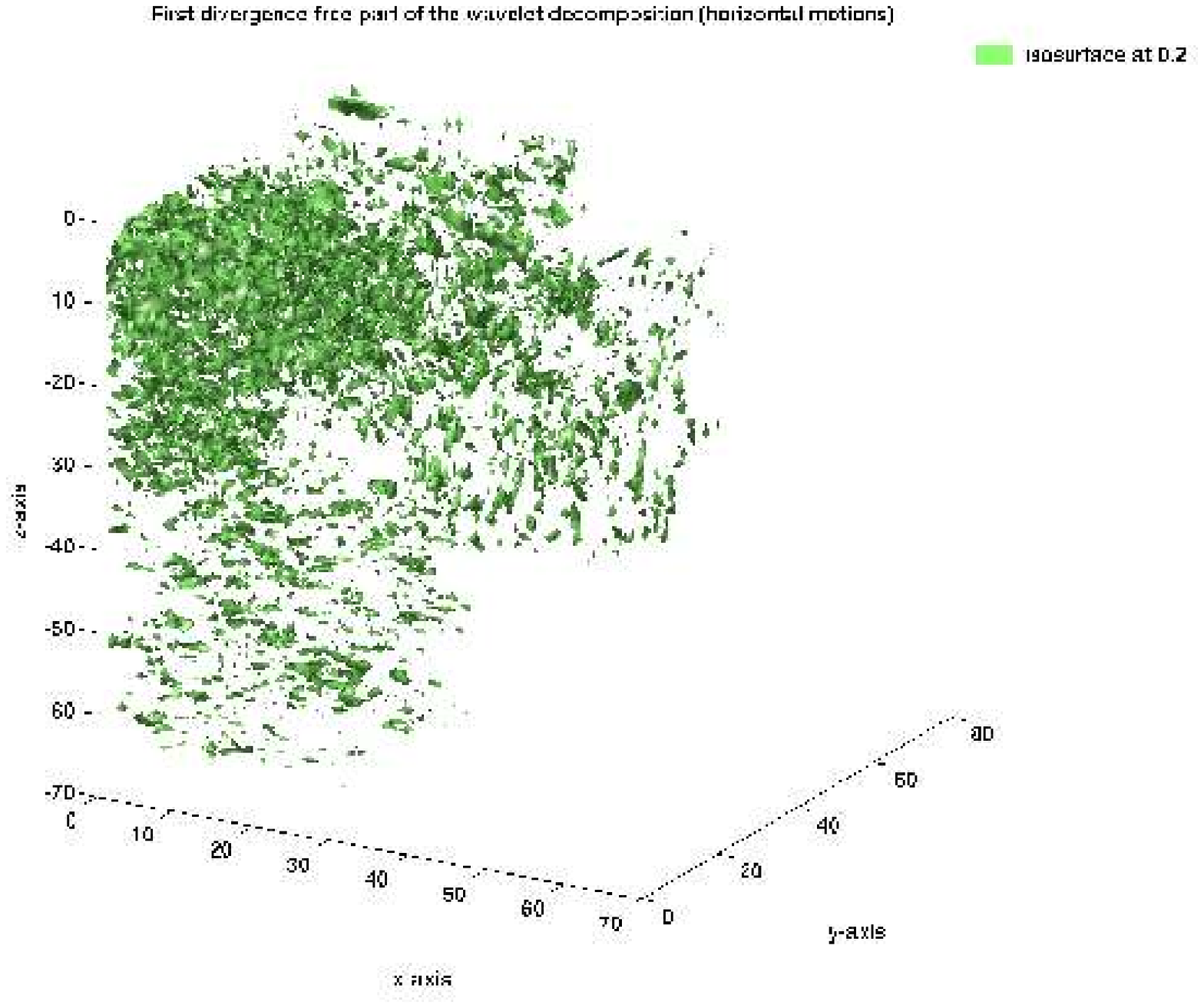}   
\includegraphics[scale=0.35]{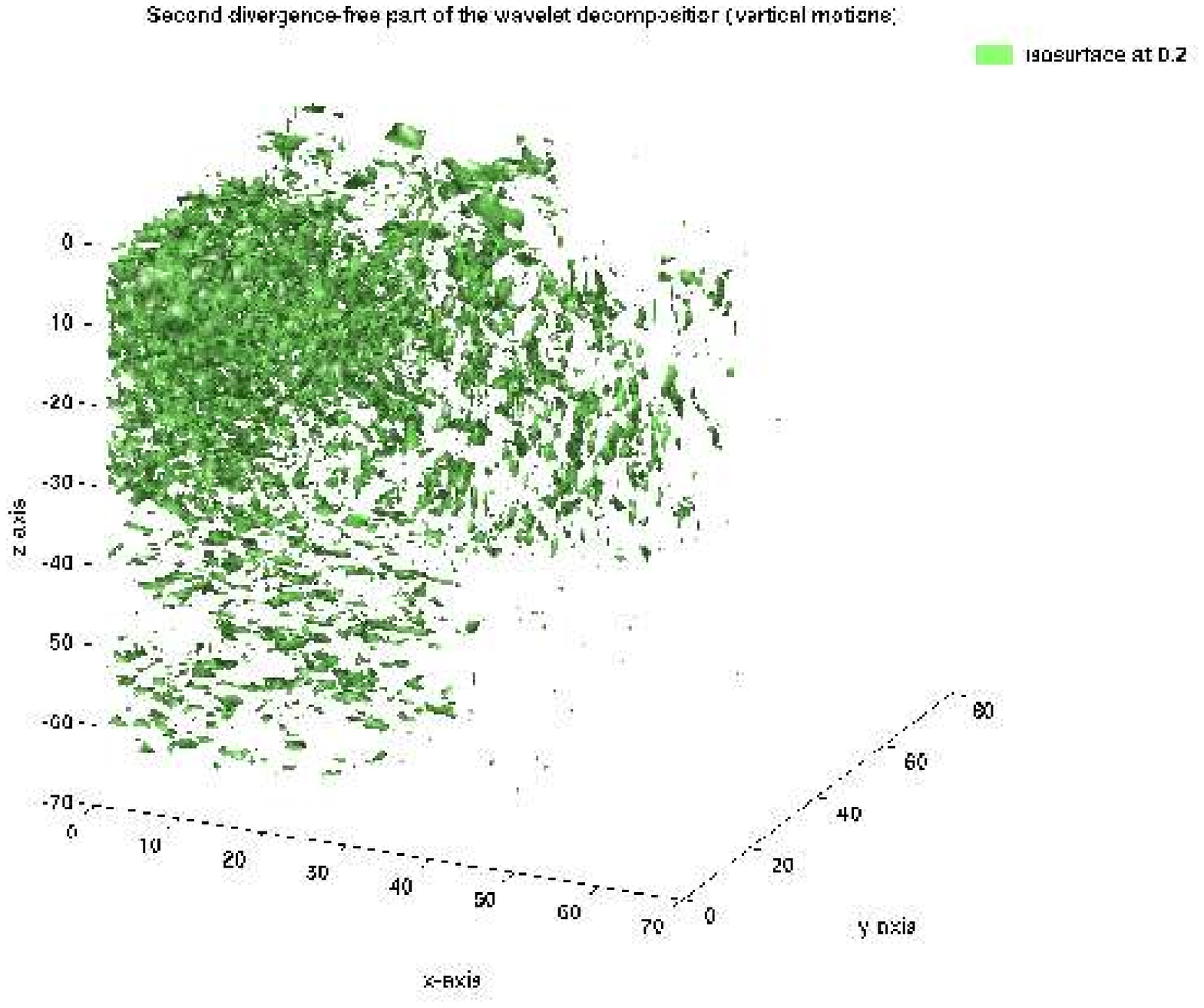}   
\end{tabular}  
\end{center} 
\caption{\label{wave3D} Isosurface $0.2$ of divergence-free wavelet coefficients
associated to $\Psi_{\textrm{div},1,j,\k}^{\varepsilon}$ (left) and to
$\Psi_{\textrm{div},2,j,\k}^{\varepsilon}$ (right),  in absolute value.}
\end{figure}

\begin{figure}[!h] 
\begin{center}  
\begin{tabular}{cc}  
\includegraphics[width=6cm,height=6cm]{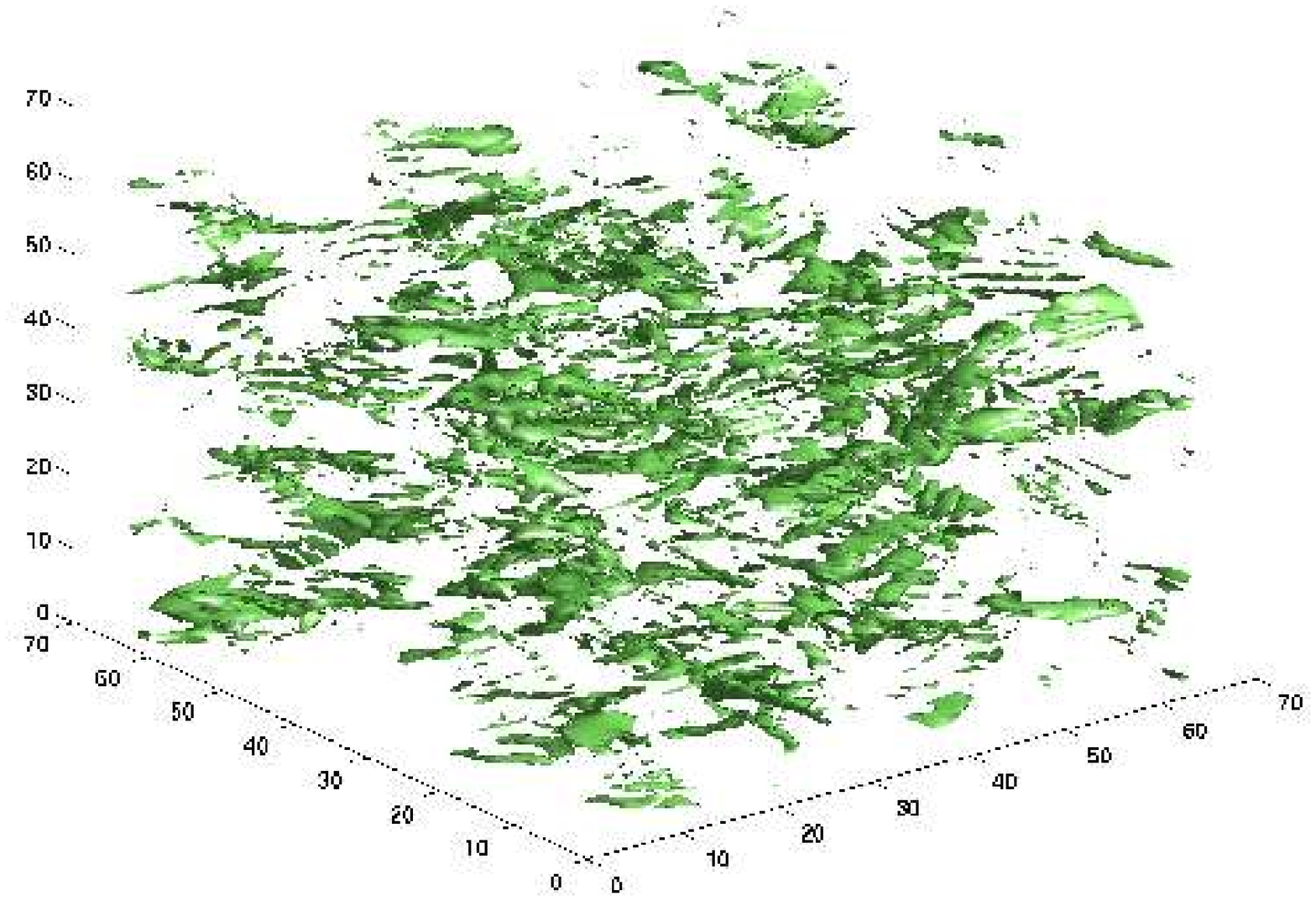}   
\includegraphics[width=6cm,height=6cm]{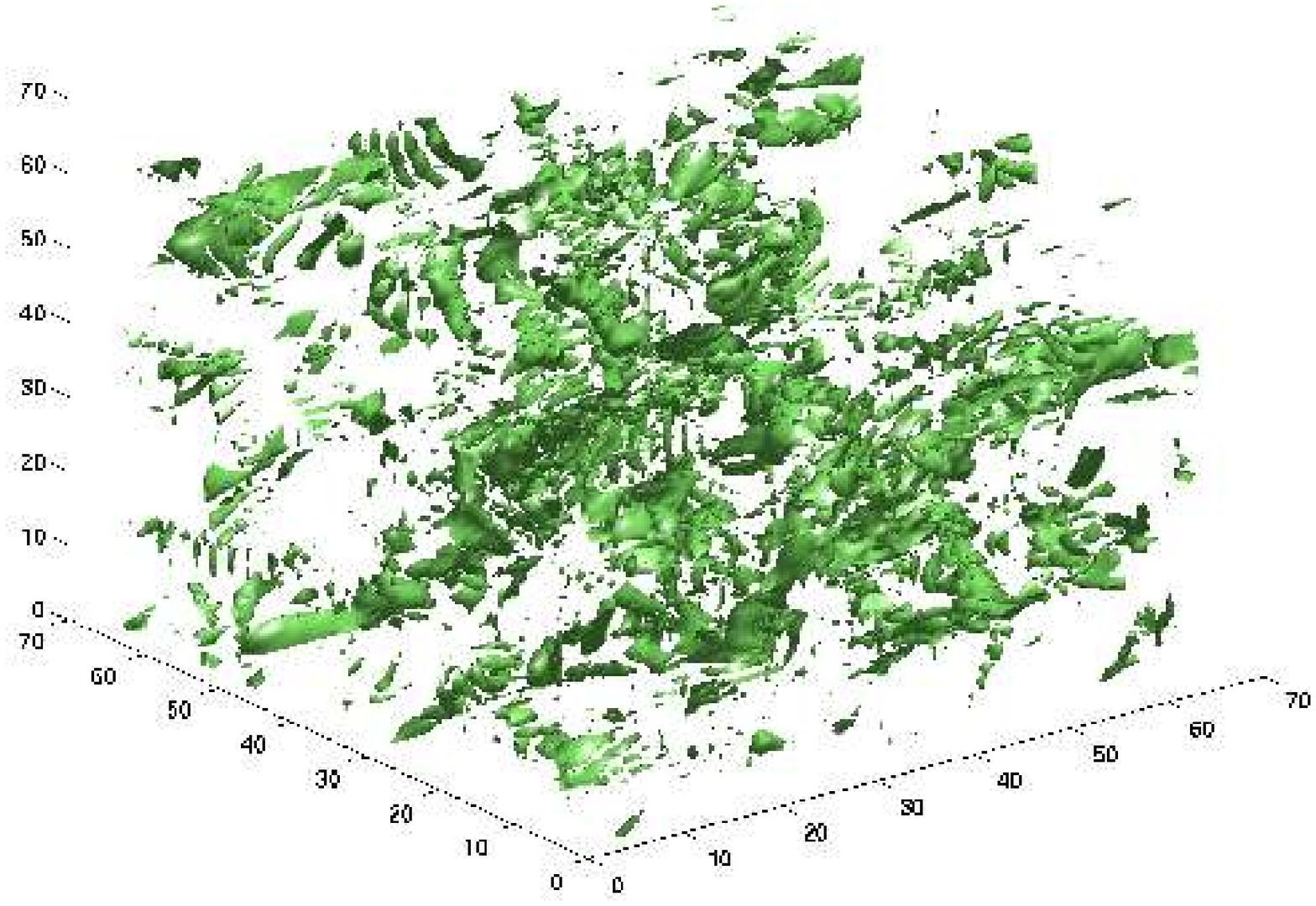}   
\end{tabular}  
\end{center} 
\caption{\label{wave3DJ} Isosurface $0.06$ of divergence-free wavelet coefficients
associated to $\Psi_{\textrm{div},2,J,\k}^{(1,0,0)}$ (left) and to
$\Psi_{\textrm{div},1,J,\k}^{(0,0,1)}$ (right),
in absolute value.} 
\end{figure}

The divergence-free wavelet decomposition of the corresponding velocity field is computed, and displayed on Figures
\ref{wave3D} and \ref{wave3DJ}. As explained in section \ref{iso3D}, the isotropic 3D divergence-free wavelet decomposition provides
with $14$ generating wavelets $\Psi_{\textrm{div},1,j,\k}^{\varepsilon}$,
$\Psi_{\textrm{div},2,j,\k}^{\varepsilon}$. \\ 
\newpage
Figure \ref{wave3D} left shows the corresponding renormalized coefficients $2^j d_{\textrm{div},1,j,\k}$, whereas Figure \ref{wave3D}
right shows the $2^j d_{\textrm{div},2,j,\k}$, until $j=6$.  The smallest scale ($j=7$) wavelet coefficients are displayed on Figure \ref{wave3DJ} below, for two kinds of
generating wavelets: we choose $\Psi_{\textrm{div},2}^{(1,0,0)}$, which corresponds to horizontal structures,
and $\Psi_{\textrm{div},1}^{(0,0,1)}$  which
exhibits vertical ones.  

\begin{figure}[!h] 
\begin{center}  
\includegraphics[width=7cm,height=4cm]{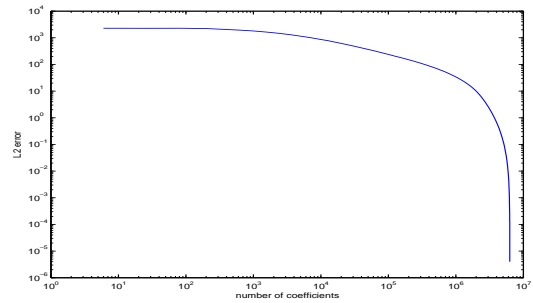}  
\end{center} 
\caption{\label{comp3D} $L^2$-error provided by the nonlinear N-best terms wavelet approximation (\ref{NL}): in log-log plot, 
$L^2$-error (\ref{errorNL}) versus $N$ for a 3D turbulent flow.} 
\end{figure}

Figure \ref{comp3D} displays the nonlinear compression error: we have computed the convergence rate on
the linear part of the graph (which is shorter by comparison with the 2D
case, due the low resolution)  and 
we have found $s\thickapprox 1.45$. 

\subsection{Analysis of 2D compressible fields}  
\label{hodge-exp}

\begin{figure}[h] 
\begin{center}  
\input{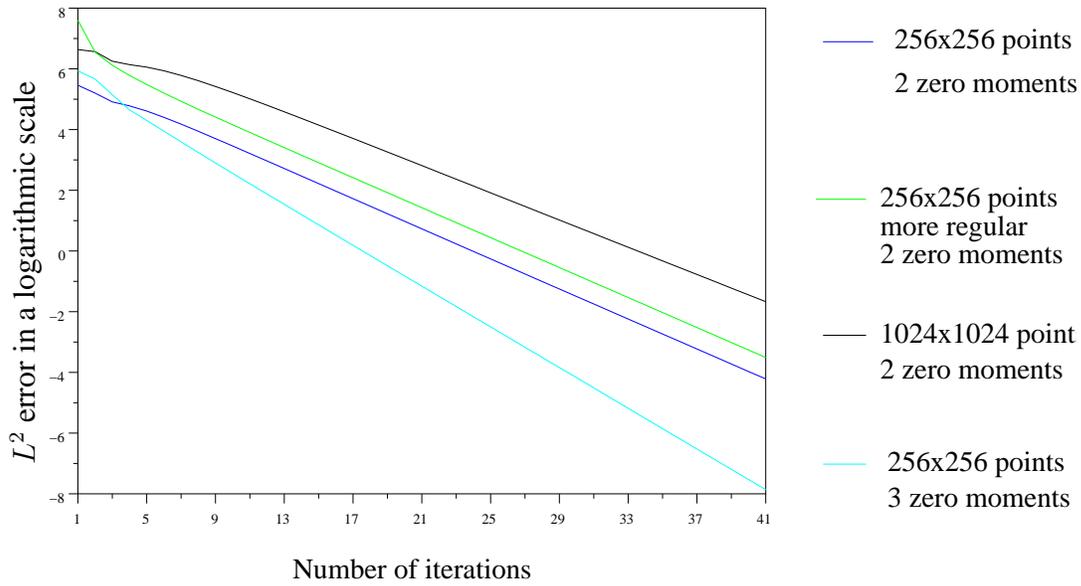}  
\end{center} 
\caption{\label{conv} Convergence curves of the iterative wavelet Hodge algorithm.} 
\end{figure}

We presented in section \ref{Hodge}, an algorithm which gave rise to a
wavelet Hodge decomposition of any flow.  In order to numerically prove that it always converges, we have tested the method on various random two-dimensional 
fields. We constructed some of them by summing random gaussians, and we modified their Fourier 
spectra in order to vary the regularity. Figure \ref{conv}
displays the $L^2$-norm of the residual, in terms of the number of iterations, for four different vector functions. We can infer the following conclusions:

- For all functions we have tested, the method converges, and curves shows that, except at the early beginning,
the convergence is exponential. 

- The slope of the curve do not depend to much on the number of grid-points, but the curve itself, corresponding to $1024^2$ grid points,  is upon these of $256^2$ grid points.
 
 - The convergence rate increase with the number of vanishing moments  of the dual wavelets.
 
In the futur, we will investigate the influence of the wavelet bases and of the interpolating projectors, on the convergence rate. 
 
 \
 
\begin{figure}[h] 
\begin{center} 
\begin{tabular}{cc} 
\includegraphics[scale=0.4]{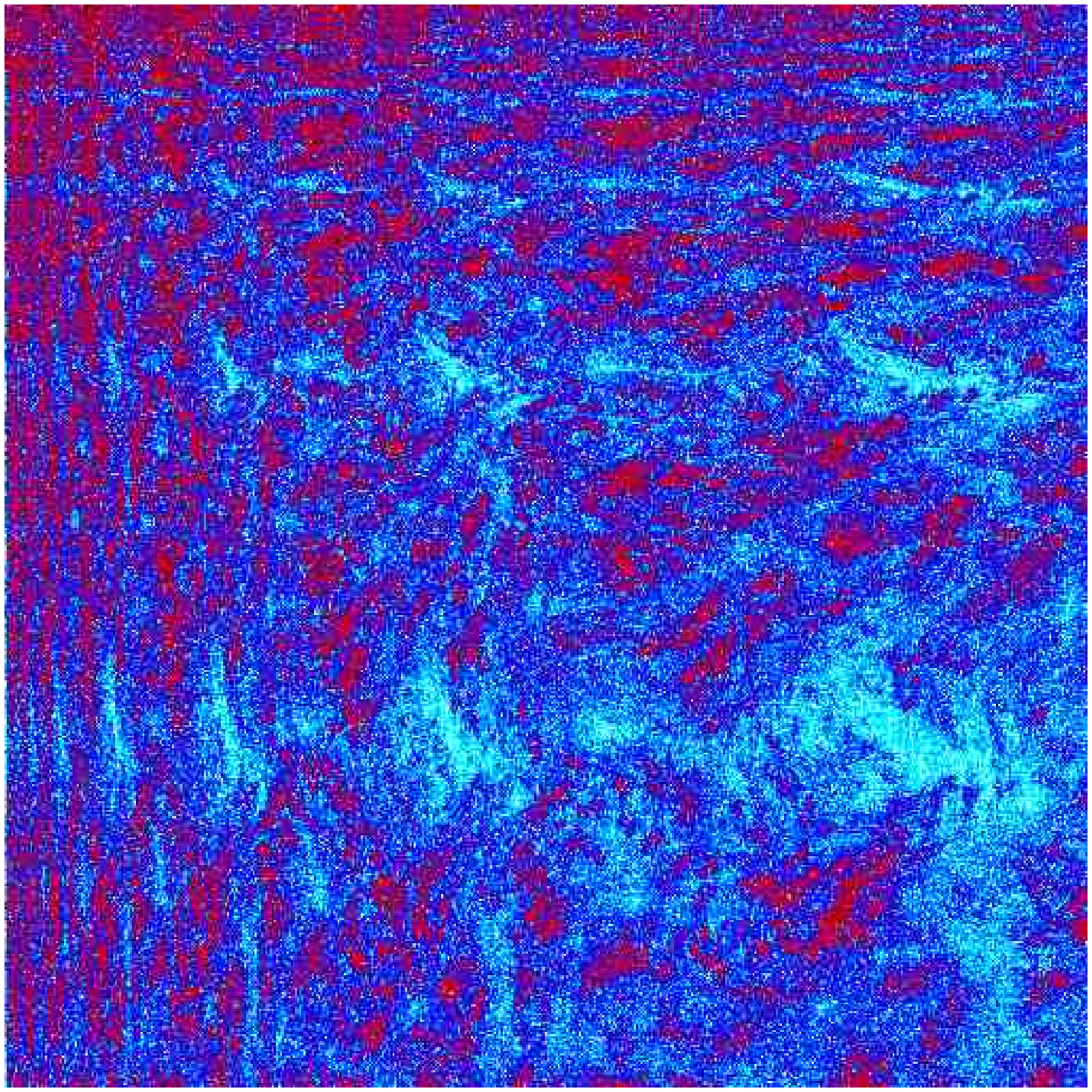} & 
\includegraphics[scale=0.4]{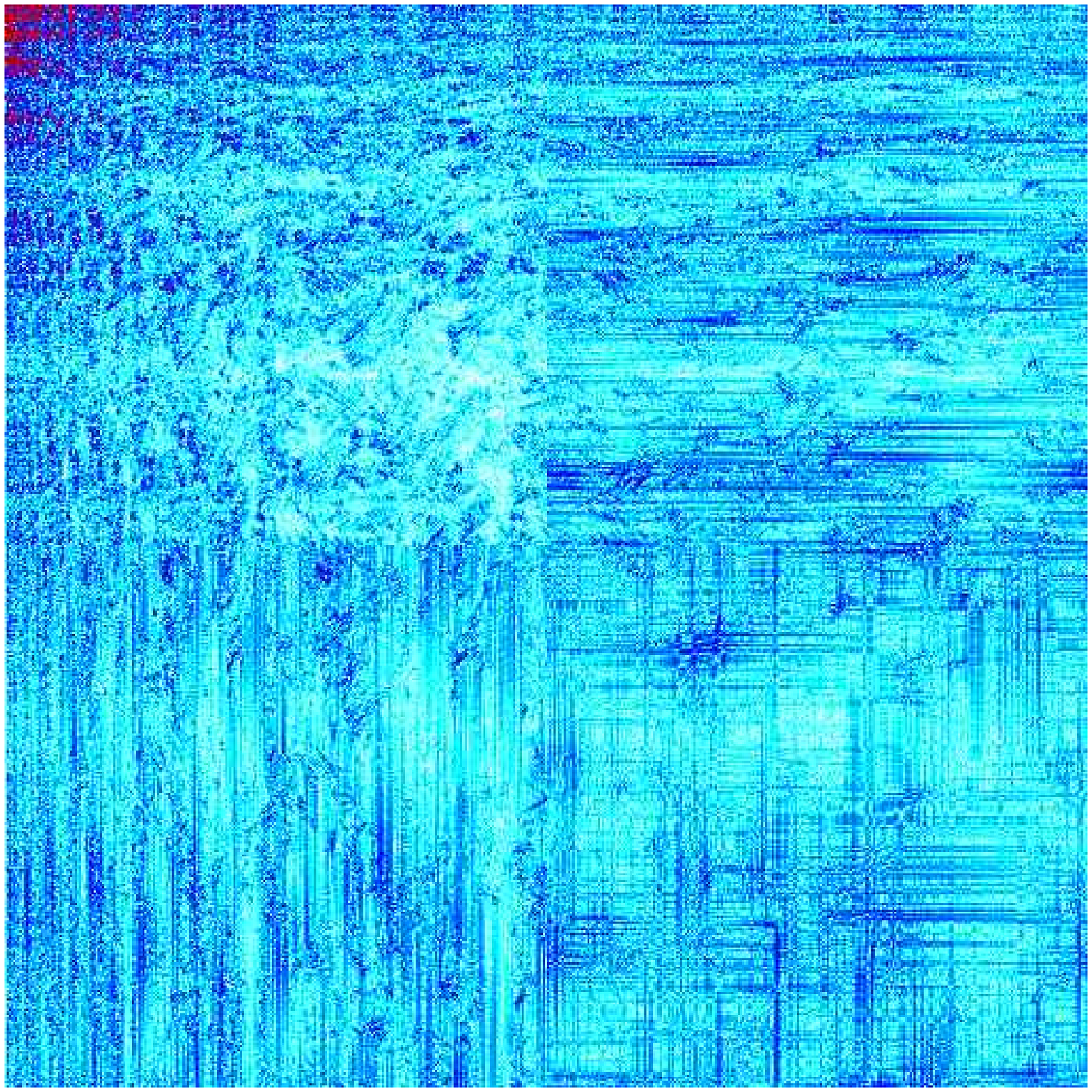} 
\end{tabular} 
\end{center} 
\caption{\label{decompoan} Anisotropic wavelet coefficients corresponding to the wavelet Hodge decomposition of $(\u.\nabla)\u$:  divergence-free coefficients (left), and curl-free coefficients (right).}
\end{figure}

Since our main objective for further research is to use divergence-free wavelets for solving the Navier-Stokes equations, we have to provide  with the wavelet Hodge decomposition of the
nonlinear term $(u.\nabla) u$. Indeed,
although $\u$ should be incompressible, this term breaks the 
divergence-free condition and yields  a compressible part.  As an illustration, 
we consider as $\u$ the 2D turbulent field displayed on Figure \ref{turb},
and we compute, by mean of our wavelet Hodge decomposition,
 the div-free and curl-free wavelet components of $(\u.\nabla) \u$. Figure \ref{decompoan} shows the anisotropic wavelet coefficients of the divergence-free part (left) and of the curl-free part (right) of the $(\u.\nabla) \u$ arising from this decomposition.

\begin{figure}[h] 
\begin{center} 
\begin{tabular}{cc} 
\includegraphics[scale=0.4]{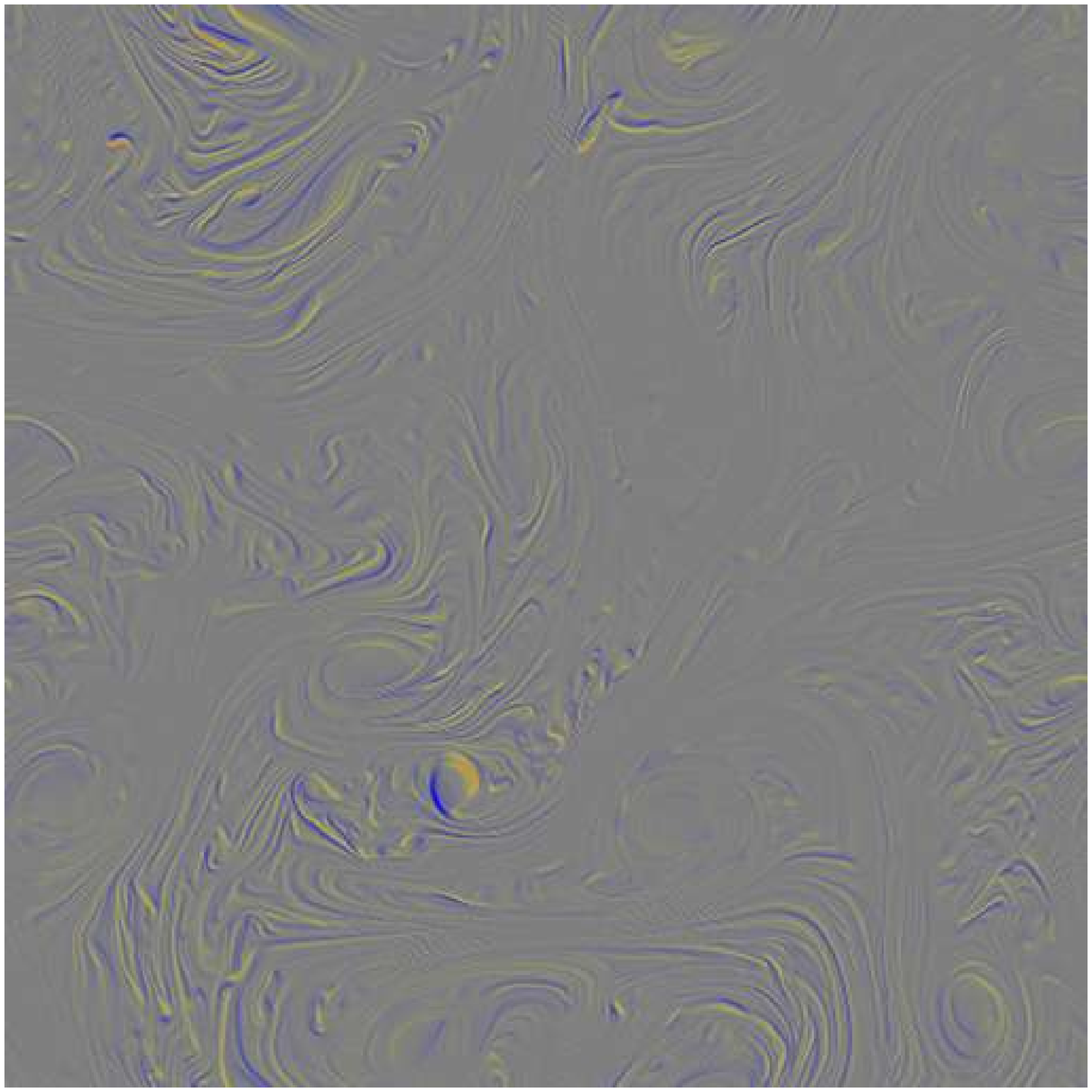} & 
\includegraphics[scale=0.4]{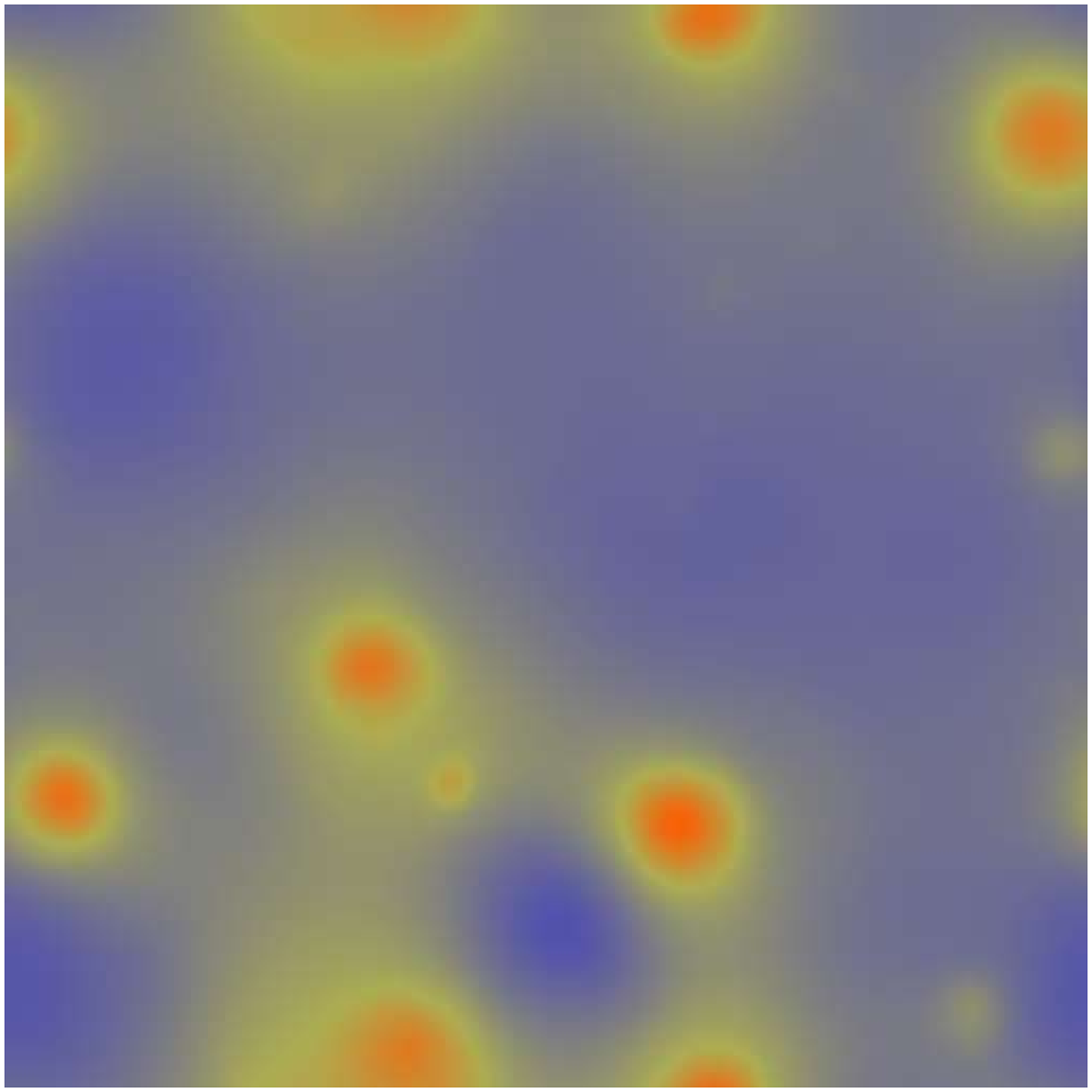} 
\end{tabular} 
\end{center} 
\caption{\label{ugradu} Vorticity (on the left) and pressure (on the right)
derived from the wavelet Hodge decomposition of the nonlinear term $\u\nabla \u$, with $\u$ displayed on Figure \ref{turb}.}
\end{figure}

 Figure \ref{ugradu} (left) displays the vorticity field associated to the divergence-free part of $(\u.\nabla) \u$, while Figure \ref{ugradu} (right) represents the pressure issued from the curl-free term, that is  easily reconstructed in wavelet domain, as it will be explain below.
 
The extracted div-free part of $(\u.\nabla) \u$ is all what is needed to compute the time-evolution of the velocity in the incompressible Navier-Stokes equations. Thanks to the curl-free wavelet definition, we are also able to directly reconstruct the pressure from the curl-free coefficients of $\nabla p$:\\
Indeed, with periodic boundary conditions, the curl-free part of $(\u.\nabla) \u$ writes:
\begin{eqnarray*}
\hspace{-2cm}
[(\u.\nabla) \u]_{\textrm{curl}~1}=\frac{\partial p}{\partial x}=&\sum_{j_1,j_2=0}^{J-1}\sum_{k_1=0}^{2^{j_1}-1}\sum_{k_2=0}^{2^{j_2}-1}
d_{\textrm{curl},k_1,k_2}~2^{j_1}~\psi_0(2^{j_1}~x-k_1)~\psi_1(2^{j_2}~y-k_2)\\
 & +\sum_{j_1=0}^{J-1}\sum_{k_1=0}^{2^{j_1}-1}d_{\textrm{curl},k_1}^1~2^{j_1}~\psi_0(2^{j_1}~x-k_1)
\end{eqnarray*}

\begin{eqnarray*}
\hspace{-2cm}
[(\u.\nabla) \u]_{\textrm{curl}~2}=\frac{\partial p}{\partial y}=&\sum_{j_1,j_2=0}^{J-1}\sum_{k_1=0}^{2^{j_1}-1}\sum_{k_2=0}^{2^{j_2}-1}
d_{\textrm{curl},k_1,k_2}~2^{j_2}~\psi_1(2^{j_1}~x-k_1)~\psi_0(2^{j_2}~y-k_2)\\
 & +\sum_{j_2=0}^{J-1}\sum_{k_2=0}^{2^{j_2}-1}d_{\textrm{curl},k_2}^2~2^{j_2}~\psi_0(2^{j_2}~y-k_2)
\end{eqnarray*}

\begin{figure}[h] 
\begin{center} 
\includegraphics[scale=0.4]{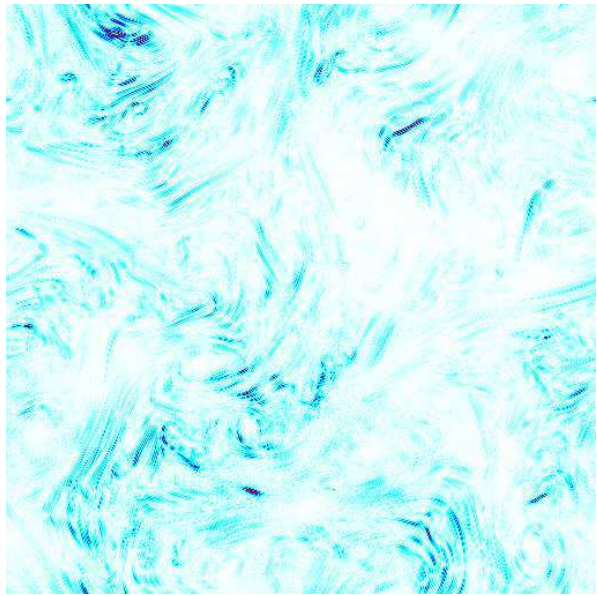} 
\end{center} 
\caption{\label{errugradu} Error between the
divergence-free part of $(\u.\nabla) \u$ (obtained
through Fourier transform) and the one provided by the wavelet Hodge decomposition.}
\end{figure}

Be integrating the system (we recall that $\psi_1'=4~\psi_0$, see the definition of gradient wavelets), we obtain (up to a constant):
\begin{eqnarray*}
\hspace{-1cm}4~p(x,y)=&\sum_{j_1,j_2=0}^{J-1}\sum_{k_1=0}^{2^{j_1}-1}\sum_{k_2=0}^{2^{j_2}-1}
d_{\textrm{curl},k_1,k_2}\psi_1(2^{j_1}~x-k_1)~\psi_1(2^{j_2}~y-k_2)\\
 & +\sum_{j_1=0}^{J-1}\sum_{k_1=0}^{2^{j_1}-1}~d_{\textrm{curl},k_1}^1\psi_1(2^{j_1}~x-k_1)\\
 & +\sum_{j_2=0}^{J-1}\sum_{k_2=0}^{2^{j_2}-1}~d_{\textrm{curl},k_2}^2~\psi_1(2^{j_2}~y-k_2)
\end{eqnarray*}
Thus the computation of the pressure is no more than a standard anisotropic wavelet reconstruction in $V^1_J\times V^1_J$, from the curl-free coefficients. 
By comparison to the pressure computed in Fourier domain, we have found a
relative error of $2.5 ~10^{-4}$
in the $L^2$-norm, which probably arises from the interpolating process.
On the other hand the difference between the Leray projection (in Fourier 
space) and the wavelet projection onto the divergence-free 
space represents $1$\% of the $L^2$-norm, that is to say 
$0.01$\% of the energy. Figure \ref{errugradu} displays the localisation of this error. 

\section*{Conclusion and perspectives}  
  
We have presented in detail the construction of 2D and 3D
divergence-free wavelet bases, and a practical way to compute the associated
coefficients. We have introduced {\it anisotropic} div-free and curl-free wavelet bases, 
which are more easy to handle. We have shown that these bases make possible an iterative 
algorithm to compute the wavelet Hodge decomposition of any flow. Thus, numerical tests
prove the feasibility of divergence-free wavelets for simulating turbulent flows in two and three dimensions.
 A divergence-free wavelet 
based solver for 2D Navier-Stokes equations is underway and will be reported in a forthcoming paper.

An important issue that must be addressed is the great ability  of the method: 
although all numerical tests have been presented in the periodic case, the method extends readily to
non-periodic problems, by using wavelets adapted to the proper boundary conditions \cite{U00, MP98}, in the
div-free 
construction. Another point is since we consider the $(\u,p)$-formulation for the Navier-Stokes equation, and
since we are able to compute the Leray projector in the wavelet domain, the method extends easily to the 3D
case.
At last, this method should be competitive by comparison to a classical Fourier method in the non-periodic case: indeed,
the periodic case corresponds to boundary conditions for which spectral methods are obviously fast, while
it is clear that wavelet methods take advantage both of the compression properties of the wavelet bases for
functions and for operators, in any case.

\ack 
The authors would like to thank G.H. Cottet and G. Lapeyre for helpfully 
providing to them numerical turbulent flows for analyses. 
This work has been supported in part by the European Community's Human Potential 
Programme under contract HPRN-CT-2002-00286, "Breaking Complexity".

 \section*{References} 
\bibliographystyle{plain}

\begin{thebibliography}{99}  
{\small  
\bibitem{AUDRL02} C.-M. Albukrek, K. Urban, W. Dahmen, D. Rempfer, and J.-L. Lumley, \emph{Divergence-Free Wavelet  
Analysis of Turbulent Flows}, J. of Scientific Computing {\bf 17}(1):  
49-66, 2002.  
  
\bibitem{CP96} P. Charton, V. Perrier,  
\emph{A pseudo-wavelet scheme for the two-dimensional Navier-Stokes  
equations},  
Comp. Appl. Math. {\bf 15}(2): 139-160, 1996.  
  
\bibitem{C00} A.Cohen, \emph{Wavelet methods in numerical analysis},  
Handbook of Numerical Analysis, vol. VII, P.G.Ciarlet and J.L.Lions eds.,  
Elsevier, Amsterdam, 2000.  
  
\bibitem{CMOV02} G.-H. Cottet, B. Michaux, S.Ossia and G. Vanderlinden, \emph{A comparison  
of spectral and vortex methods  
in three-dimensional incompressible flows},  J. Comp. Phys, 175, 2002.  
  
\bibitem{DKU96} W. Dahmen, A. Kunoth and K. Urban, \emph{A wavelet-Galerkin method for the  
Stokes problem}, Computing 56, 1996, 259-302.  
  
\bibitem{Dau92} I. Daubechies, \emph{Ten lectures on Wavelets}, SIAM book,  
Philadelphia, Pennsylvania, 1992.   
 
\bibitem{Bo} C. De Boor, \emph{A Practical Guide to Splines}, book, Springer-Verlag New York Inc., 
2001. 
 
\bibitem{F92} M. Farge,  
\emph{Wavelet transforms and their applications to turbulence},  
Ann. Rev. Flu. Mech. {\bf}:395-457, 1992.  
  
\bibitem{FK96} M. Farge, N. Kevlahan, V. Perrier \& E. Goirand,  \emph{Wavelets and turbulence},  
Proc. IEEE 84(4), 639-669,1996.  
  
\bibitem{FS01} M. Farge and K. Schneider, \emph{Coherent Vortex  
Simulation (CVS), A Semi-Deterministic Turbulence Model Using Wavelets}, Flow,  
Turbulence and Combution, {\bf 66}: 393-426, 2001.  
  
\bibitem{FS96} J. Fr\"ohlich  and K. Schneider,  
\emph{ Numerical simulation of decaying turbulence in an adaptive wavelet basis},  
   Appl. Comput. Harmon. Anal., {\bf 3}: 393-397,  
  1996.  
  
  \bibitem{GR79} V. Girault, P.A. Raviart,  
 \emph{ Finite element approximations of the 
		  Navier-Stokes equations},  
Lecture Notes  in Mathematics, 
  {Springer-Verlag},  
  {1979}.  
   
 \bibitem{GK00} M.~Griebel and F.~Koster,  
 \emph{ Adaptive wavelet solvers for the unsteady incompressible  
		  Navier-Stokes equations},  
Advances in Mathematical Fluid Mechanics,  
  {J.~Malek and J.~Necas and M.~Rokyta eds},  
  {Springer-Verlag},  
  {2000}.  
   
  \bibitem{HH86} B.L. Hua and D. Haidvogel, 
  \emph{Numerical simulations of the vertical structure of quasi-geostrophic turbulence}, 
  J. Atmos. Sci.   {\bf 3}:2923-2936, 1986.  
  
\bibitem{KKR00} 
 J. Ko,  A.J. Kurdila and  O.K. Rediniotis, 
\emph{ divergence-free Bases and Multiresolution Methods for 
Reduced-Order Flow Modeling}, 
 AIAA Journal, {\bf 38}(2): 2219-2232, 2000. 
   
 \bibitem{KGKFS98} F.~Koster and M.~Griebel and N.~Kevlahan and M.~Farge and  
		  K.~Schneider,  
  \emph{Towards an adaptive wavelet-based {3D} {N}avier-{S}tokes  
		  solver},  
Numerical flow simulation I, Notes on Numerical Fluid Mechanics, Vol.  
66, 339-364, E.H. Hirschel eds, Vieweg-Verlag, Braunschweig, 1998.  
  
  
\bibitem{KL98} J.-P. Kahane and P.-G. Lemari\'e-Rieusset, \emph{Fourier  
series and wavelets}, book, Gordon \& Breach, London, 1995.  
  
\bibitem{L2000} G. Lapeyre, \emph{Topologie de m\'elange dans un fluide turbulent  
g\'eophysique} (in french), Th\`ese de doctorat de l'Universit\'e Paris VI, 2000.  
  
\bibitem{L92} P.-G. Lemari\'e-Rieusset, \emph{Analyses  
multi-r\'esolutions non orthogonales, commutation entre projecteurs et  
d\'erivation et ondelettes vecteurs \`a divergence nulle} (in french), Revista Matem\'atica  
Iberoamericana, {\bf 8}(2): 221-236, 1992.  
  
\bibitem{L94} P.- G. Lemari\'e-Rieusset, \emph{Un th\'eor\`eme d'inexistance pour les  
ondelettes vecteurs \`a divergence nulle} (in french), C. R. Acad. Sci. Paris, t. 319,  
S\'erie I, p. 811-813, 1994.  
  
\bibitem{L93} J. Lewalle,  
\emph{Wavelet transform of the Navier-Stokes equations and the  
generalized dimensions of turbulence},   
Appl. Sci. Res. {\bf 51}(1-2):109-113, 1993.  
  
\bibitem{M98} S. Mallat, \emph{A Wavelet Tour of Signal Processing}, book,  
Academic Press, 1999.  
  
  
\bibitem{M91} C. Meneveau,  
\emph{Analysis of turbulence in the orthonormal wavelet  
representation},  
Journal of Fluid Mechanics {\bf 232}: 469-520, 1991.  
 
\bibitem{MP98} P. Monasse and V. Perrier,
 \emph{Orthonormal wavelet bases adapted for partial differential
equations with boundary conditions}, SIAM J. on Math. Analysis {\bf
29}(4):1040-1065, 1998.  

\bibitem{SKF97} K. Schneider, N. Kevlahan, M. Farge,  
 \emph{Comparison of an adaptive wavelet method and nonlinearly filtered pseudo-spectral methods for two-dimensional turbulence}, 
 Theor. Comput. Fluid Dyn. {\bf 9}: 191-206, 1997. 
  
\bibitem{U94} K. Urban, \emph{A Wavelet-Galerkin Algorithm for the  
Driven--Cavity--Stokes--Problem in Two Space Dimensions}, RWTH Aachen,  
Preprint 1994.  
  
\bibitem{U96}  
K. Urban,  
   \emph{Using divergence free wavelets for the numerical solution of the {S}tokes problem},  
   {AMLI}'96: Proceedings of the Conference on Algebraic Multilevel Iteration Methods with Applications,  
   {\bf 2}: 261--277,   
University of Nijmegan, The Netherlands, 1996.  
  
\bibitem{U00} K. Urban, \emph{Wavelet Bases in H(div) and H(curl)}, Mathematics  
of Computation {\bf 70}(234): 739-766, 2000.  
  
\bibitem{U02} K. Urban, \emph{Wavelets in Numerical Simulation}, Springer,  
2002.  
  
}  
  
\end{thebibliography}

 
\Figures 
  
\Figure{From left to right: the scaling function $\phi$ with its associated symmetric wavelet 
with shortest support, and their duals: the dual scaling function $\phi^*$ and the dual
wavelet $\psi^*$.}
\Figure{Scaling functions and associated wavelets with shortest support, for splines of degree
1 (left) and 2 (right).}
\Figure{Anisotropic 2D wavelet transform.}
\Figure{Isotropic 2D wavelet transform.}
\Figure{Isotropic 2D generating divergence free wavelets $\Psi_{\textrm{div}}^{(1,0)}$ (left), 
$\Psi_{\textrm{div}}^{(0,1)}$ (center) and  $\Psi_{\textrm{div}}^{(1,1)}$ (right).}
\Figure{The two scaling functions $\phi_0$ and $\phi_1$, and their symmetry centers.}

\Figure{Vorticity fields at times $t=0$, $t=10$, $t=20$ and $t=40$, and corresponding 
divergence-free wavelet coefficients of the velocity.} 
\Figure{Vorticity field for a $1024\times 1024$ simulation of decaying turbulence (left),
and the corresponding divergence-free wavelet coefficients of the velocity field (right).}
\Figure{$L^2$-error provided by the nonlinear N-best terms wavelet approximation (\ref{NL}):
in log-log plot,  $L^2$-error (\ref{errorNL}) versus $N$ for a 2D turbulent flow.}
\Figure{Isosurface of vorticity magnitude after 5 large-eddy turnovers provided by a spectral
method \cite{CMOV02}.}
\Figure{Isosurface $0.2$ of divergence-free wavelet coefficients
associated to $\Psi_{\textrm{div},1,j,\k}^{\varepsilon}$ (left) and to
$\Psi_{\textrm{div},2,j,\k}^{\varepsilon}$ (right) in absolute value.}
\Figure{Isosurface $0.06$ of divergence-free wavelet coefficients
associated to $\Psi_{\textrm{div},2,J,\k}^{(1,0,0)}$ (left) and to
$\Psi_{\textrm{div},1,J,\k}^{(0,0,1)}$ (right),
in absolute value.}

\Figure{$L^2$-error provided by the nonlinear N-best terms of wavelet approximation (\ref{NL}):
in log-log plot,  $L^2$-error (\ref{errorNL}) versus $N$ for a 3D turbulent flow.}
\Figure{Convergence curves of the iterative wavelet Hodge algorithm.}
\Figure{Anisotropic wavelet coefficients corresponding to the wavelet Hodge decomposition of $(\u.\nabla)\u$:  divergence-free coefficients (left), and curl-free coefficients (right).}

\Figure{Vorticity (on the left) and pressure (on the right)
derived from the wavelet Hodge decomposition of the nonlinear term $\u\nabla \u$, with $\u$ displayed on Figure \ref{turb}.}

\Figure{Error between the
divergence-free part of $(\u.\nabla) \u$ (obtained
through Fourier transform) and the one provided by the wavelet Hodge decomposition.}

\end{document}